\documentclass[final,3p,times]{elsarticle}
\usepackage{graphicx}
\usepackage[latin1]{inputenc}
\usepackage[svgnames]{xcolor}
\usepackage[french,english]{babel}
\usepackage{amsmath,amssymb,amsthm,epsfig,graphics,float,subfigure}
\usepackage{thmtools}
\usepackage{subfigure,color,colortbl}
\usepackage{natbib}
\usepackage[
  bookmarks,
  linkcolor= blue,
  citecolor= blue,
  filecolor= blue,
  urlcolor=  blue,
  colorlinks=true,
  hyperindex=true,
  pdftitle={Density-dependent},
  pdfauthor={R. Fekih-Salem, T. Sari}
]{hyperref}
\usepackage{etoolbox}
\addto\captionsfrench{}
\input rotate.tex

\newtheorem{theorem}{Theorem}
\AfterEndEnvironment{theorem}{\noindent\ignorespaces}

\newtheorem{Proposition}{Proposition}

\theoremstyle{remark}
\newtheorem{Remark}{Remark}
\newcommand{{\sLP}}{\ensuremath{ {\mbox{\tiny{LP}}} }}
\newcommand{{\sBP}}{\ensuremath{ {\mbox{\tiny{BP}}} }}
\newcommand{{\sH}}{\ensuremath{ {\mbox{\tiny{H}}} }}

\def\sc{\scriptsize}
\def\sl{\small}
\def\ds{\displaystyle}
\def\scM{\scriptstyle}

\graphicspath{{Images/}}

\usepackage[mathlines]{lineno}
\usepackage{lineno}
\biboptions{numbers,sort&compress}
\begin{document}
\begin{frontmatter}

\cortext[cor1]{Corresponding author}
\title{The operating diagram of a flocculation model in the chemostat and its dependence on the biological parameters}
\author[a,c]{{\bf Radhouane Fekih-Salem}}  \corref{cor1}  
\ead{radhouane.fekih-salem@enit.utm.tn}
\author[b]{{\bf Tewfik Sari}}
\ead{tewfik.sari@inrae.fr}
\address[a]{University of Tunis El Manar, National Engineering School of Tunis, LAMSIN, 1002, Tunis, Tunisia}
\address[b]{ITAP, Univ Montpellier, INRAE, Institut Agro, Montpellier, France}
\address[c]{University of Monastir, Higher Institute of Computer Science of Mahdia, 5147, Mahdia, Tunisia}
\begin{abstract}
In this paper, we consider a flocculation model in a chemostat where one species is present in two forms: planktonic and aggregated bacteria with the presence of a single resource. The removal rates of isolated and attached bacteria are distinct and include the specific death rates. Considering distinct yield coefficients with a large class of growth rates, we present a mathematical analysis of the model by establishing the necessary and sufficient conditions of the existence and local asymptotic stability of all steady states according to the two operating parameters which are the dilution rate and the inflowing concentration of the substrate. Using these conditions, we determine first theoretically the operating diagram of the flocculation process describing the asymptotic behavior of the system with respect to two control parameters. The bifurcations analysis shows a rich set of possible types of bifurcations: transcritical bifurcation or branch points of steady states, saddle-node bifurcation or limit points of steady states, Hopf, and homoclinic bifurcations. Using the numerical method with MATCONT software based on a continuation and correction algorithm, we find the same operating diagram obtained theoretically. However, MATCONT detects other types of two-parameter bifurcations such as Bogdanov-Takens and Cusp bifurcations.
\end{abstract}
\end{frontmatter}
\section{Introduction}
The chemostat is an important laboratory apparatus used for experiments on the controlled growth of microorganisms in microbiology and ecology.
It has played an important role in many fields, such as the wastewater treatment process, biomass energy recovery, and biotechnologies in a broad sense.
Mathematical models of competition on a single limiting nutrient in a chemostat have played a central role in microbial ecology, microbiology, and evolutionary and applied biology.
The mathematical study of the classical chemostat model of several species competing on the same limiting resource can be found in the monograph by Smith and Waltman \cite{SmithBook1995}. They have shown that only the most competitive species that consumes less substrate to reach its steady state (or that has the lowest {\em break-even concentration}) survives the competition of several species on a single nutrient while all other species are excluded.
This result is well known as the {\em Competitive Exclusion Principle} (CEP) which states that two species competing for identical limited resource cannot coexist indefinitely.
However, the CEP contradicts the biodiversity observed in nature and microbial ecosystems.

In order to reconcile the mathematical results of the classical chemostat model asserting the CEP and the experimental results and nature showing the biodiversity of microbial species, various recent studies have revised the mathematical modeling of the competition of several microbial species competing on a single resource.
More specifically, a lot of research has tried to understand and explain the biodiversity in microbial ecosystems by analyzing the various types of interactions favoring the coexistence of microbial species.
In \cite{RobledoNARWA2012}, the constant input of some species in a chemostat of $n$ species competing on a single nutrient can lead to coexistence.
In the literature, we can cite these various mechanisms of coexistence:
flocculation \cite{FekihJMAA2013,FekihAMM2016,FekihSIADS2019,FekihArima2020,HaegemanJBD2008},
intra- and interspecific interference \cite{AbdellatifMBE2016,DeLeenheerJMAA2006},
density-dependence \cite{FekihMB2017,LobryCRBio2006EN,LobryCRMA2005,LobryCRBio2006,MtarIJB2020,MtarDcdsB2022},
presence of internal or external inhibitors \cite{BarDCDSB2020,DellalDcDsB2021,DellalDcDsB2022,DellalMB2018,HsuJapanJIAM1998},
predator-prey interaction \cite{BoerMB1998,ButlerJMB1983},
simple or complex food web \cite{ButlerJMB1986,HsuDcDsB2009,WolkowiczMB1989},
and the references therein.


Allelopathy and bacteriocin represent another mechanism of coexistence between a wild-type organism and a single mutant in the chemostat \cite{ZouNARWA2011}.
Using the specific growth rates of Monod-type, the authors study the existence and local stability of the steady states.
The analysis of the bifurcations shows that there can be either a transcritical or a pitchfork bifurcation \cite{ZouNARWA2011}.

This paper is a follow-up to a previous work \cite{FekihSIADS2019,FekihArima2020} where we considered a flocculation model of one microbial species that is decomposed into isolated (or planktonic)
bacteria and attached (or aggregate) bacteria with a single nutrient $S$ in a chemostat. Moreover, isolated bacteria can aggregate with isolated bacteria or flocs to form new flocs, with a rate $a(u + v)u$,
while flocs can split and liberate isolated bacteria, with a rate $bv$. This model was introduced in \cite{FekihJMAA2013} and was also considered in \cite{HarmandBook2017,RapaportEJAM2018}.
The model is given by the following three-dimensional system of ordinary differential equations
\begin{equation}                                   \label{ModFlocGen}
\left\{
\begin{array}{lll}
\dot{S}  &=& D(S_{in}-S)-\frac{1}{y_u}f(S)u-\frac{1}{y_v}g(S)v \\[1mm]
\dot{u}  &=& [f(S)-D_u]u-a(u+v)u+b v \\[1mm]
\dot{v}  &=& [g(S)-D_v]v+a(u+v)u-b v
\end{array}
\right.
\end{equation}
where $S(t)$ is the concentration of the substrate at time $t$;
$f(S)$ and $g(S)$ represent, respectively, the growth rates of isolated and attached bacteria;
$D$ and $S_{in}$ are, respectively, the dilution rate and the concentration of the substrate in the feed device;
$D_u$ and $D_v$ represent, respectively, the disappearance rates of planktonic and attached bacteria;
$y_u$ and $y_v$ are ``yield'' constants reflecting the conversion of nutrient to planktonic and aggregated bacteria, respectively.
In \cite{FekihSIADS2019}, we have considered model (\ref{ModFlocGen}) in the case where the yields coefficients $y_u$ and $y_v$ are equal.
In this case, they can be normalized to 1. However, because of the structure of model (\ref{ModFlocGen}), when these coefficients are distinct, they cannot be normalized to 1 by the usual change of variable where $u$ and $v$ are replaced by $u/y_u$ and $v/y_v$, respectively.
Consideration of these yield coefficients is very important in the mathematical models of the chemostat to model reproduction by nutrient uptake as mentioned in \cite{HarmandBook2017,SmithBook1995}.
In this work, we study model (\ref{ModFlocGen}) where $D_u$ and $D_v$ can be modeled as in \cite{Marsili1996,Shen2007} by:
\begin{equation}                                      \label{ExpDuDv}
D_u = \alpha D + m_u, \quad D_v = \beta D + m_v
\end{equation}
where the non-negative parameters $m_u$ and $m_v$ representing mortality rate are taken into consideration.


In \cite{FekihSIADS2019}, we have determined the existence and local stability of all steady states of system (\ref{ModFlocGen}) with the same yields coefficients $y_u$ and $y_v$.
The model presents a multiplicity of positive steady states that can only appear or disappear through saddle-node or transcritical bifurcations.
Under the joined effect of flocculation and mortality, the coexistence steady states may destabilize via a supercritical Hopf bifurcations with emergence of a stable limit cycle that can disappear through a homoclinic bifurcation.
However, the study of bifurcations is limited to the one-parameter diagrams by fixing $D$ and varying $S_{in}$.


In \cite{FekihArima2020}, the theoretical study of the operating diagram of model (\ref{ModFlocGen}) with the same yields coefficients shows that the system can exhibit bistability between the washout steady state $E_0$ and the coexistence steady state $E_1$.
There may also only be coexistence around the positive steady state $E_1$.
The construction of the operating diagram of model (\ref{ModFlocGen}) in \cite{FekihArima2020} has omitted the existence of the region of destabilization of a positive steady state where there can be the emergence of a stable limit cycle via a Hopf bifurcation for very small values of $D$ as demonstrated in \ref{Sec-AppedixCasArima}.


Indeed, the operating diagram is a very useful tool to visualize and summarize the asymptotic behavior of a process according to the operating parameters which are the most easily manipulated parameters in a chemostat as explained in \cite{HarmandBook2017,SmithBook1995}. In the existing literature, the study of the operating diagram can be purely numerical.
By exploring the set of operating parameters $D$ and $S_{in}$ with a certain discretization step, a significant steady state (i.e. with nonnegative components) is determined by solving numerically algebraic equations giving steady states. Their asymptotic behaviors are established by solving the characteristic polynomial and the sign of its roots \cite{WadePloS2017}.
This method can be applied to complex processes with a large number of state variables and parameters \cite{HanakiProc2021,KhedimAMM2018,WadeJTB2016,WeedermannNonlinDyn2015,XuJTB2011}.


Another numerical alternative consists in constructing the boundaries of the various regions of the operating diagram using a continuation and correction algorithm.
Various software packages have been developed in order to determine the values of the critical parameters corresponding to the different types of bifurcations for autonomous dynamic systems.
The most used software packages are MATCONT, CONTENT, AUTO, and XPPAUT (see \cite{DhoogeMCMDS2008} and the reference therein).


However, the theoretical determination of the operating diagram consists in constructing the boundaries of the different regions from the theoretical analysis of the dynamic system.
More precisely, using a scientific numerical platform (like MAPLE \cite{MAPLE18}), these boundaries are drawn from the conditions of existence and stability of all steady states according to the operating parameters when all biological parameters are fixed \cite{AbdellatifMBE2016,BarDCDSB2020,DaliYoucefMBE2020,DaoudMMNP2018,DellalDcDsB2021,DellalDcDsB2022,DellalMB2018,DhoogeMCMDS2008,FekihArima2020,FekihMB2017,FekihSIADS2021,
MtarIJB2020,MtarDcdsB2022,NouaouraJMB2022,SariNonLinDyn2021,SariMB2016,SariMB2017}.
Note that the single-parameter or two-parameter bifurcation diagrams obtained with MATCONT \cite{MATCONT} allow additional phenomena to be detected (such as homoclinic, Cusp, and Bogdanov-Takens bifurcations) compared with those obtained theoretically from the existence and stability conditions.

Our main objective in this paper is to extend our mathematical study in \cite{FekihSIADS2019} by considering distinct yields and
to describe theoretically and numerically the operating diagram of model (\ref{ModFlocGen}).
Moreover,
this work is an extension of our study presented in \cite{FekihSIADS2019}, which is limited to the numerical analysis of the bifurcation diagram according to the single parameter $S_{in}$.
Thus, this study of the operating diagram provides a more general analysis of the asymptotic behavior of solutions of the system according to the two operating parameters $S_{in}$ and $D$.
In addition, our aim is to use bifurcation theory to complement previous studies.
Moreover, our in-depth theoretical study of the operating diagram shows the emergence of a region of destabilization of the positive steady state via a Hopf bifurcation with coexistence around a stable limit cycle.
The one- and two-parameter diagrams are also obtained by the numerical continuation method using MATCONT software \cite{MATCONT}, which allowed us to detect other types of bifurcations according to two parameters.
In addition, the effect of attachment and detachment on the operating diagram is analyzed theoretically to show the importance of considering the phenomenon of flocculation as a coexistence mechanism in the classic chemostat model.


This paper is organized as follows.
First, we present in Section \ref{Sec-HypAnalMod} a general hypothesis about the growth functions of flocculation model (\ref{ModFlocGen}).
Then, we determine the existence and the local stability conditions of all steady states according to the dilution rate and the input concentration of the substrate.
In Section \ref{Sec-DO}, we analyze theoretically the operating diagram.
First, in Section \ref{SubSec-DOCasFig13}, a simple case is considered where there is only a Branch Point (BP) and no Limit Point (LP) or Hopf bifurcation.
In Section \ref{SubSec-DOCasFig12Siam}, a case with LP and Hopf bifurcations is considered.
In Section \ref{SubSec-DOCasFig6}, another case with LP and Hopf bifurcations is considered but a new region of instability of two positive steady states emerges in the operating diagram.
In Section \ref{Sec-NumDO}, we study numerically the operating diagram and the bifurcation diagram according to one parameter using the software MATCONT for the two cases in Sections \ref{SubSec-DOCasFig12Siam} and \ref{SubSec-DOCasFig6}.
In Section \ref{Sec-EffectFloc}, we study the effect of flocculation on the operating diagram for the set of parameters considered in Sections \ref{SubSec-DOCasFig6} but where $a$ and $b$ are variable.
Finally, conclusions are drawn in the last Section \ref{Sec-Conc}.
In \ref{Sec-AppedixCasFig13}, we show that a stability condition of the positive steady state holds for the set of parameters considered in Section \ref{SubSec-DOCasFig13}.
In \ref{Sec-AppedixCasFig12}, we illustrate that this stability condition does not hold for the set of parameters considered in Section \ref{SubSec-DOCasFig12Siam}.
In \ref{Sec-AppedixCasFig6}, we show the destabilization of a positive steady state and then illustrate the stable limit cycles in the three-dimensional space $(S,u,v)$ and their disappear via a homoclinic bifurcation
for the set of parameters considered in Section \ref{SubSec-DOCasFig6}.
In \ref{Sec-AppedixCasArima}, we show that the region of destabilization of the positive steady state is omitted in the construction of the operating diagram in \cite{FekihArima2020}.
All the values of parameters used throughout this paper are provided in \ref{Sec-AppedixParamVal}.
\section{Hypothesis and model analysis}         \label{Sec-HypAnalMod}
In this paper, we make the following general assumption on the growth functions $f(S)$ and $g(S)$ which are continuously differentiable ($\mathcal{C}^1$).\\
\noindent
(H1) $f(0)=g(0)=0$ and $f'(S)>0$ and $g'(S)>0$ for all $S>0$.  \\
Assumption (H1) means that no growth can occur for isolated bacteria $u$ and attached bacteria $v$ without the presence of the substrate $S$.
Moreover, the growth rates of isolated and attached bacteria increase with the concentration of the substrate $S$.

In this section, we summarize the main results of the existence and stability of all steady states of system (\ref{ModFlocGen}).
A steady state exists if and only if all its components are nonnegative. This predicts two types of steady states, labeled as follows:
\begin{itemize}
\item $E_0$ ($u=0$, $v=0$): the washout of planktonic and attached bacteria.
\item $E_1$ ($u>0$, $v>0$): both planktonic and attached bacteria are present.
\end{itemize}

To determine these steady states, we define the following auxiliary functions
\begin{equation}                                           \label{eqH}
H(S):=\frac{1}{y_u}f(S)U(S)+\frac{1}{y_v}g(S)V(S)
\end{equation}
where
\begin{equation}                                          \label{eqUV}
U(S):=\frac{\varphi(S)(\psi(S)-b)}{a[\psi(S)-\varphi(S)]}  \quad \mbox{and} \quad V(S):=-\frac{\varphi^2(S)(\psi(S)-b)}{a[\psi(S)-\varphi(S)]\psi(S)}
\end{equation}
and
\begin{equation}                                      \label{eqphipsi}
\varphi(S):=f(S)-D_u \quad \mbox{and} \quad \psi(S):=g(S)-D_v.
\end{equation}
In addition, we need to define the following interval of existence of the positive steady states:
$$
 I=]\lambda_u,\lambda_v[ \quad \mbox{if} \quad \lambda_u<\lambda_v, \quad \mbox{else} \quad  I= ]\lambda_v,\min(\lambda_u,\lambda_b)[
$$
For convenience, we shall use the abbreviation LES for Locally Exponentially Stable.
Any reference to steady state stability should be considered as local exponential stability, that is to say, the real parts of the eigenvalues of the Jacobian matrix are negative.
To determine the stability of the positive steady state $E_1=(S^*,u^*,v^*)$, we define the Routh--Hurwitz coefficients by
\begin{equation}                                        \label{ExpRHC}
\begin{array}{c}
c_1=m_{11}+m_{22}+m_{33},\\[1mm]
c_2=m_{12}m_{21}+m_{13}m_{31}-m_{32}a_{23}+m_{11}m_{22}+m_{11}m_{33}+m_{22}m_{33},\\[1mm]
c_3=m_{11}(m_{22}m_{33}-m_{32}a_{23})+m_{21}(m_{12}m_{33}+m_{32}m_{13})+m_{31}(m_{12}a_{23}+m_{13}m_{22})\\[1mm]
c_4=c_1c_2-c_3.
\end{array}
\end{equation}
where
\begin{equation}                                        \label{Exp-mij}
\left\{
\begin{array}{lll}
m_{11}=D+\frac{1}{y_u}f'(S^*)u^*+\frac{1}{y_v}g'(S^*)v^*, \quad m_{12}=\frac{1}{y_u}f(S^*), \quad m_{13}=\frac{1}{y_v}g(S^*),   \\[1mm]
m_{21}=f'(S^*)u^*, \quad m_{22}=a(2u^*+v^*)-\varphi(S^*), \quad a_{23}=b-au^*,\\[1mm]
m_{31}=g'(S^*)v^*, \quad m_{32}=a(2u^*+v^*) \quad \mbox{and} \quad m_{33}=b-au^*-\psi(S^*).
\end{array}
\right.
\end{equation}
Now, we can state the main result which establishes the components of all steady states of (\ref{ModFlocGen}) and their existence and local asymptotic stability conditions.
\begin{theorem}                                    \label{TheoExisGen}
Assume that Hypothesis (H1) holds. The steady states of (\ref{ModFlocGen}) and the necessary and sufficient conditions of existence and local asymptotic stability are given in Tables \ref{TabSSi} and \ref{TabExisStabGen}, respectively.
\end{theorem}
\begin{table}[ht]
\caption{Steady states of (\ref{ModFlocGen}). The functions $H(S)$, $U(S)$ and $V(S)$ are defined by (\ref{eqH}) and (\ref{eqUV}).} \label{TabSSi}
\vspace{-0.3cm}
\begin{center}
\begin{tabular}{ @{\hspace{0mm}}l@{\hspace{1mm}} @{\hspace{1mm}}l@{\hspace{-2mm}} }
              &  $S$, $u$, $v$ components
\\ \hline
$E_0$        &  $S=S_{in}$, $u=0$, $v=0$
\\ \hline
$E_1$        &  $S^*$ solution of equation $D(S_{in}-S)=H(S)$, $u^*=U(S^*)$ and $v^*=V(S^*)$
\end{tabular}
\end{center}
\vspace{-0.3cm}
\end{table}
\begin{table}[ht]
\caption{Necessary and sufficient existence and local stability conditions of steady states of (\ref{ModFlocGen}) where $c_4$ is defined by (\ref{ExpRHC}).} \label{TabExisStabGen}
\vspace{-0.3cm}
\begin{center}	
\begin{tabular}{ @{\hspace{1mm}}l@{\hspace{1mm}} @{\hspace{1mm}}l@{\hspace{1mm}} @{\hspace{1mm}}l@{\hspace{1mm}}}			
                        &   Existence conditions  &   Stability conditions
\\ \hline
$E_0$                   & always exists                &  $S_{in}<\min(\lambda_u,\lambda_b)$.
\\ \hline
$E_1$                  &
\begin{tabular}{l}		
equation $D(S_{in}-S)=H(S)$ \\
has a solution $S^*\in I$
\end{tabular}
 &
\begin{tabular}{l}		
$c_3=\varphi(S^*)(b-\psi(S^*))(D+H'(S^*))>0$ and \\
$c_4>0$
\end{tabular}
\\ \hline
\end{tabular}
\end{center}
\vspace{-0.3cm}
\end{table}
\begin{proof}
The proof for the components of the steady states given in Table \ref{TabSSi} and their existence conditions given in Table \ref{TabExisStabGen} is the same as the proof of \cite[Lemma 2.4]{FekihSIADS2019}.
The proof for the stability condition of $E_0$ is the same as the proof of \cite[Proposition 3.1]{FekihSIADS2019}.
The proof for the stability condition of $E_1$ is the same as the proof of \cite[Proposition 3.3]{FekihSIADS2019}.
\end{proof}
It was shown in \cite{HarmandBook2017} (see also \cite{RapaportEJAM2018}) that when $D_u=D_v=D$, then the positive steady state $E_1$ exists and is unique and LES if and only if $S_{in}>\lambda_u$.
\section{Operating diagrams}                            \label{Sec-DO}
In this section, we study theoretically the operating diagrams of model (\ref{ModFlocGen}) to determine the various qualitative asymptotic behaviors of the process according to the operating parameters which are the concentration of substrate in the feed bottle $S_{in}$ and the dilution rate $D$.
Each region of the diagram is characterized by a different color according to the number of existing steady states and their various asymptotic behaviors.
Except of the operating parameters $S_{in}$ and $D$, all the biological parameters are fixed since they cannot be easily manipulated by the biologist as they depend on the nature of the organisms and the substrate introduced into the bioreactor.

From definition of $\lambda_v(D)$ and $\lambda_b(D)$ in Table \ref{TabFuncDO}, we have $\lambda_v(D)<\lambda_b(D)$ for all $D\in[0,\overline{D}_v[$.
To construct theoretically the operating diagram of system (\ref{ModFlocGen}) by determining the various curves,
we define the auxiliary functions according to the dilution rate $D$ in Table \ref{TabFuncDO} and the set of curves $\Gamma_i$,
$i= \left\{ u, v, {\rm{\scM{BP}}},{\rm{\scM{LP}}},{\rm{\scM{H}}} \right\}$
in Table \ref{TabCurvDO}, which are the boundaries of different regions of the $(S_{in}, D)$-plane.
As in this work, the construction of the operating diagram will be done with the specific growth rates of Monod-type (\ref{SpeciFunc}) satisfying hypothesis (H1) and we know that in this case the function $S\mapsto H(S)$ is convex, we can define the functions $D \mapsto S_{\sLP}(D)$ and $D\mapsto \lambda_{\sLP}(D)$ in Table \ref{TabFuncDO} by following \cite{FekihArima2020}.
\begin{table}[ht]
\caption{Notations, auxiliary functions, and their domains of definition.} \label{TabFuncDO}
\vspace{-0.2cm}
\begin{center}
\begin{tabular}{ @{\hspace{1mm}}l@{\hspace{1mm}} @{\hspace{1mm}}l@{\hspace{1mm}}  }	
                               & \qquad \qquad Definition
\\ \hline
$\lambda_u(D)$                 &
\begin{tabular}{l}
$\lambda_u(D) = f^{-1}(\alpha D + m_u)$.   \\
It is defined for $0 \leq D < \overline{D}_u:=(f(+\infty)-m_u)/\alpha$.
\end{tabular}
\\ \hline
$\lambda_v(D)$                 &
\begin{tabular}{l}
$\lambda_v(D) = g^{-1}(\beta D + m_v)$.   \\
It is defined for $0 \leq D < \overline{D}_v:=(g(+\infty)-m_v)/\beta$.
\end{tabular}
\\ \hline
$\lambda_b(D)$                 &
\begin{tabular}{l}
$\lambda_b(D) = g^{-1}(\beta D + m_v + b)$.   \\
It is defined for $0 \leq D < \overline{D}_b:=(g(+\infty)-m_v-b)/\beta$.
\end{tabular}
\\ \hline
$\lambda_{\sBP}(D)$                 &
\begin{tabular}{l}
$\lambda_{\sBP}(D) = \min(\lambda_u(D),\lambda_b(D))$.   \\
It is defined for $0 \leq D < \max\left(\overline{D}_u,\overline{D}_b\right)$.
\end{tabular}
\\ \hline
$S=S_{\sLP}(D)$                 &
\begin{tabular}{l}
$S=S_{\sLP}(D)$ is the unique solution of equation $H'(S)=-D$ on $]\lambda_v(D),\lambda_{\sBP}(D)]$.  \\
It is defined for $\overline{D} \leq D < \overline{D}_v$ where $\overline{D}$ is the unique solution in $\left]0,\max\left(\overline{D}_u,\overline{D}_b\right)\right[$\\
    of equation $H'(\lambda_{\sBP}(D))=-D$.
\end{tabular}
\\ \hline
$\lambda_{\sLP}(D)$                 &
\begin{tabular}{l}
$\lambda_{\sLP}(D) = H(S_{\sLP}(D)) / D + S_{\sLP}(D)$.   \\
It is defined for $\overline{D} \leq D < \overline{D}_v$.
\end{tabular}
\end{tabular}
\end{center}
\vspace{-0.1cm}
\end{table}

The passage through the $\Gamma_{\rm{\scM{BP}}}$ curve corresponds to a transcritical bifurcation or Branch Point (BP) between $E_0$ and $E_1^1$ or between $E_0$ and $E_1^2$ as we will see in the following section.
As we shall see later, passing through the curve $\Gamma_{\sLP}$ in the operating plan $(S_{in},D)$ gives rise to the two positive steady states $E_1^1$ and $E_1^2$ via a Limit Points (LP) or saddle-node bifurcation.
In addition, the passage through the curve $\Gamma_H$ corresponds to Hopf bifurcation with the appearance or disappearance of a stable limit cycle.
\begin{table}[ht]
\caption{Definitions of the curves $\Gamma_i$, $i= \left\{ u, v, {\rm{\scM{BP}}},{\rm{\scM{LP}}},{\rm{\scM{H}}} \right\}$, and the corresponding colors where all functions $\lambda_i$ and $c_4$ are defined in Table \ref{TabFuncDO} and (\ref{ExpRHC}), resp. The abbreviations BP, LP, and H mean a Branch Point, Limit Point, and Hopf bifurcations, respectively.}  \label{TabCurvDO}
\vspace{-0.2cm}
\begin{center}
\begin{tabular}{ @{\hspace{1mm}}l@{\hspace{2mm}}  @{\hspace{2mm}}l@{\hspace{1mm}} @{\hspace{2mm}}l@{\hspace{1mm}} }	
Curves                                                                    &  Color    &   Bifurcation  \\ \hline
$\Gamma_u = \left\{(S_{in},D): S_{in} = \lambda_u(D) \right\}$            &  Red        &  BP           \\
$\Gamma_b = \left\{(S_{in},D): S_{in} = \lambda_b(D) \right\}$            &  Blue       &  BP           \\
$\Gamma_{\sBP} = \left\{(S_{in},D): S_{in} = \lambda_{\sBP}(D) \right\}$  &  Red or Blue&  BP           \\
$\Gamma_{\sLP} =\left\{(S_{in},D): S_{in} = \lambda_{\sLP}(D) \right\}$   &  Green      &  LP            \\
$\Gamma_{\sH}=     \left\{(S_{in},D): c_4(S_{in}, D) = 0 \right\}$        &  Magenta    &  H
\end{tabular}\end{center}
\end{table}
To illustrate the operating diagram of model (\ref{ModFlocGen}), we choose the following specific growth rates of Monod-type satisfying hypothesis (H1):
\begin{equation}                                     \label{SpeciFunc}
f(S)=\frac{m_1S}{k_1+S} \quad \mbox{and} \quad g(S)=\frac{m_2S}{k_2+S},
\end{equation}
where $m_i$ is the maximum growth rate and $k_i$ is the Michaelis-Menten constant, $i=1,2$.
The values of these biological parameters are provided in Table \ref{TabParamVal}.
In addition, the construction of the operating diagram is similar for any other specific growth rate satisfying hypothesis (H1).

In the next section, we start with the simplest case where the operating diagram does not present the regions of destabilization of the positive steady state and the emergence of two positive steady states.
Then, we study the general case with the emergence of the limit cycle and BP bifurcation. Then, we find these results using the numerical continuation method with the MATCONT software.
Finally, we determine the effect of flocculation on the appearance and disappearance of various regions.
\subsection{A case where the positive steady state is unique and stable if it exists}  \label{SubSec-DOCasFig13}
In this section, we consider a case where there is only BP of two steady states and no LP or Hopf bifurcation. Therefore, the positive steady state $E_1$ is unique and stable if it exists and can only bifurcate with the washout steady state $E_0$. For this purpose, we consider the biological parameter values that were used in \cite[Fig. 13]{FekihSIADS2019} (see Table \ref{TabParamVal}, line 1).
However, $S_{in}$ and $D$ are variable and not fixed as in \cite{FekihSIADS2019} where $S_{in}=5$ and $D=3.5$.
With this set of parameters, we have $\lambda_u(D)<\lambda_v(D)$ for all $D$ in their definition domain so that the function $S\mapsto H(S)$ is defined and increasing on $\left[\lambda_u(D),\lambda_v(D)\right)$ as shown in Fig. \ref{Fig-CasFig13SiamHC4}(a). Moreover, \ref{Sec-AppedixCasFig13} shows that the stability condition $c_4>0$ holds for all $S_{in}$ and $D$ in the existence domain of $E_1$.

\begin{figure}[!h]
\setlength{\unitlength}{1.0cm}
\begin{center}
\begin{picture}(5.8,5.5)(0,0)
\put(-5.1,0){\rotatebox{0}{\includegraphics[width=5.6cm,height=7cm]{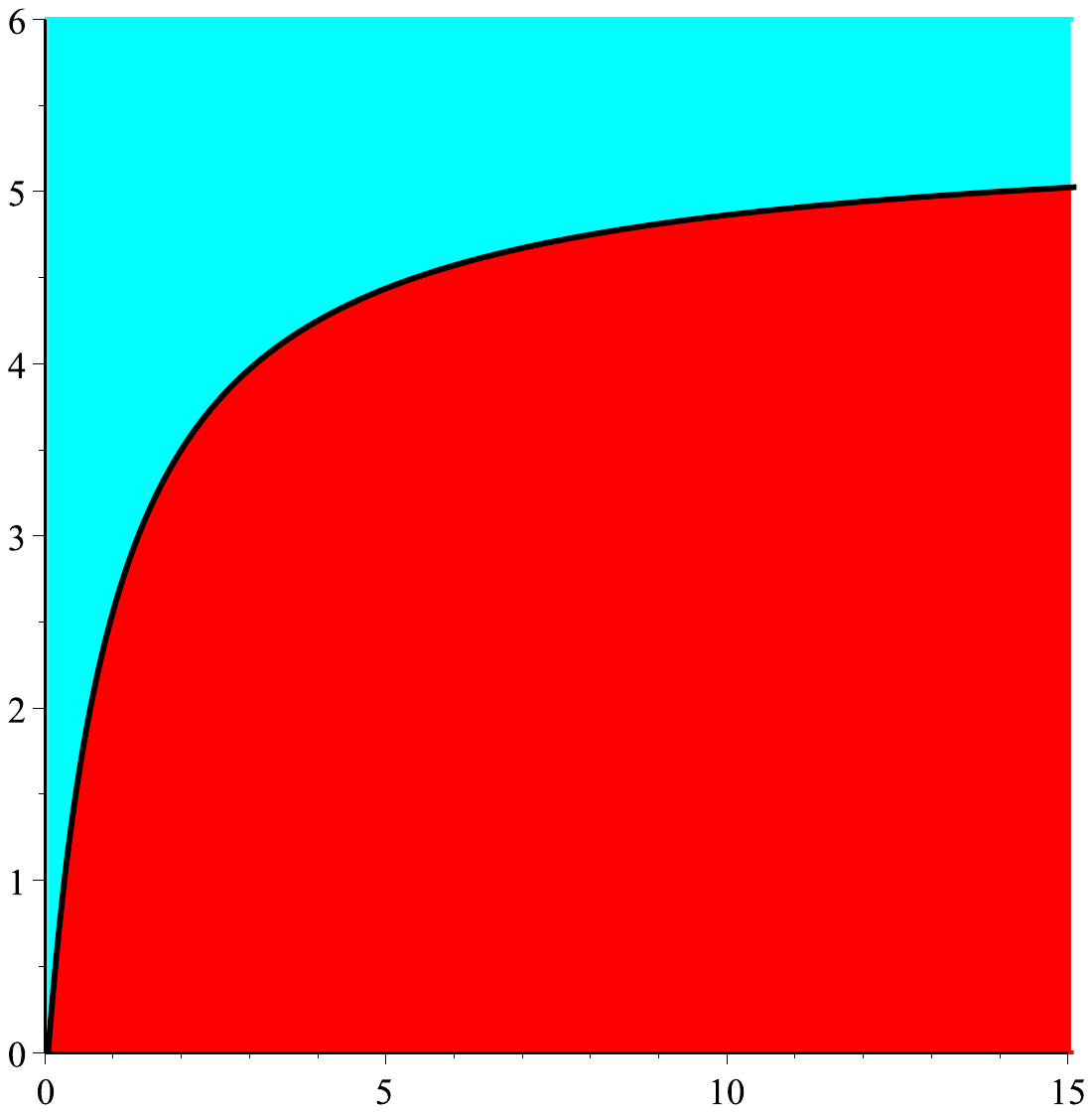}}}
\put(0,0){\rotatebox{0}{\includegraphics[width=5.6cm,height=7cm]{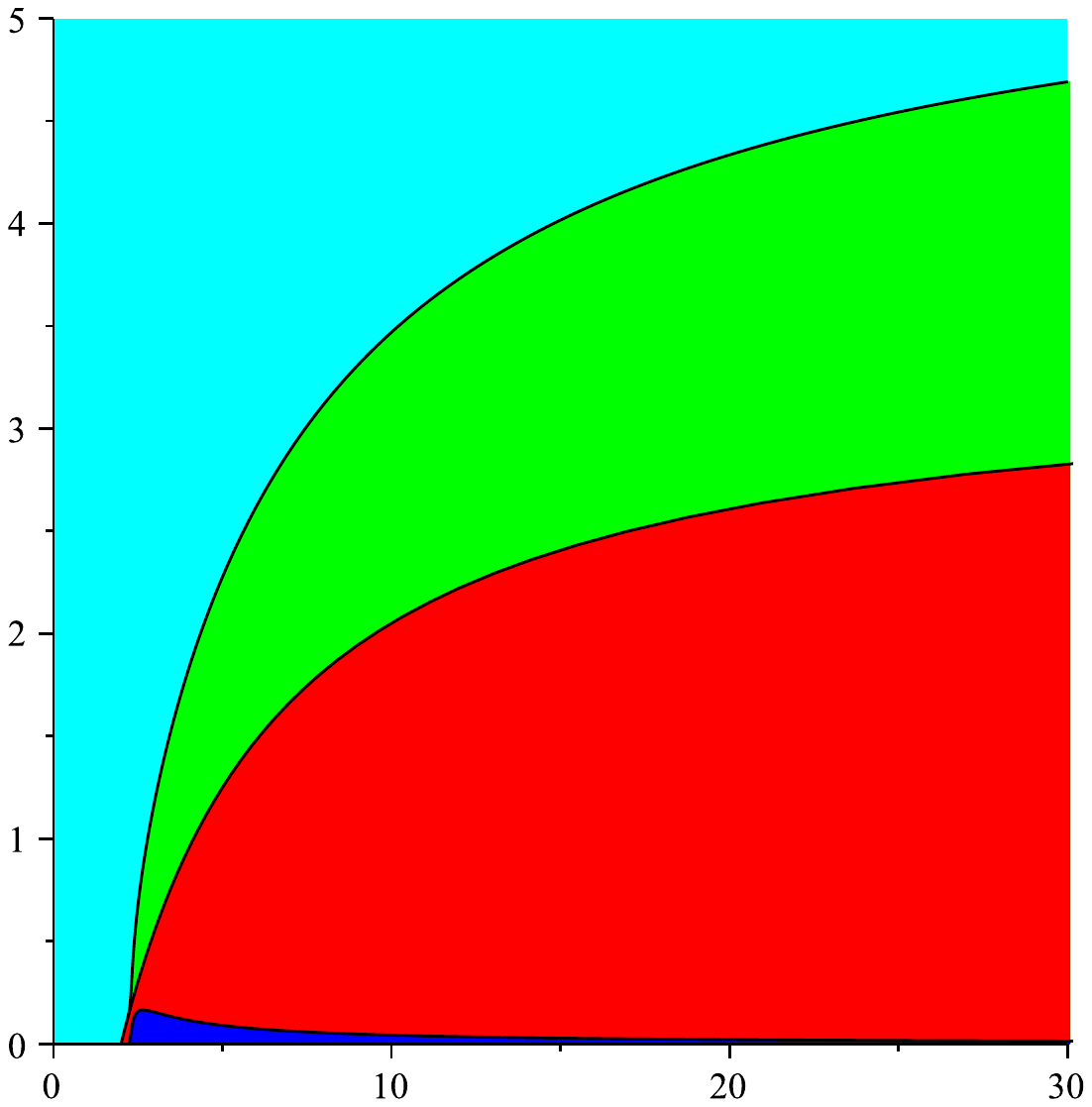}}}
\put(5.1,0){\rotatebox{0}{\includegraphics[width=5.6cm,height=7cm]{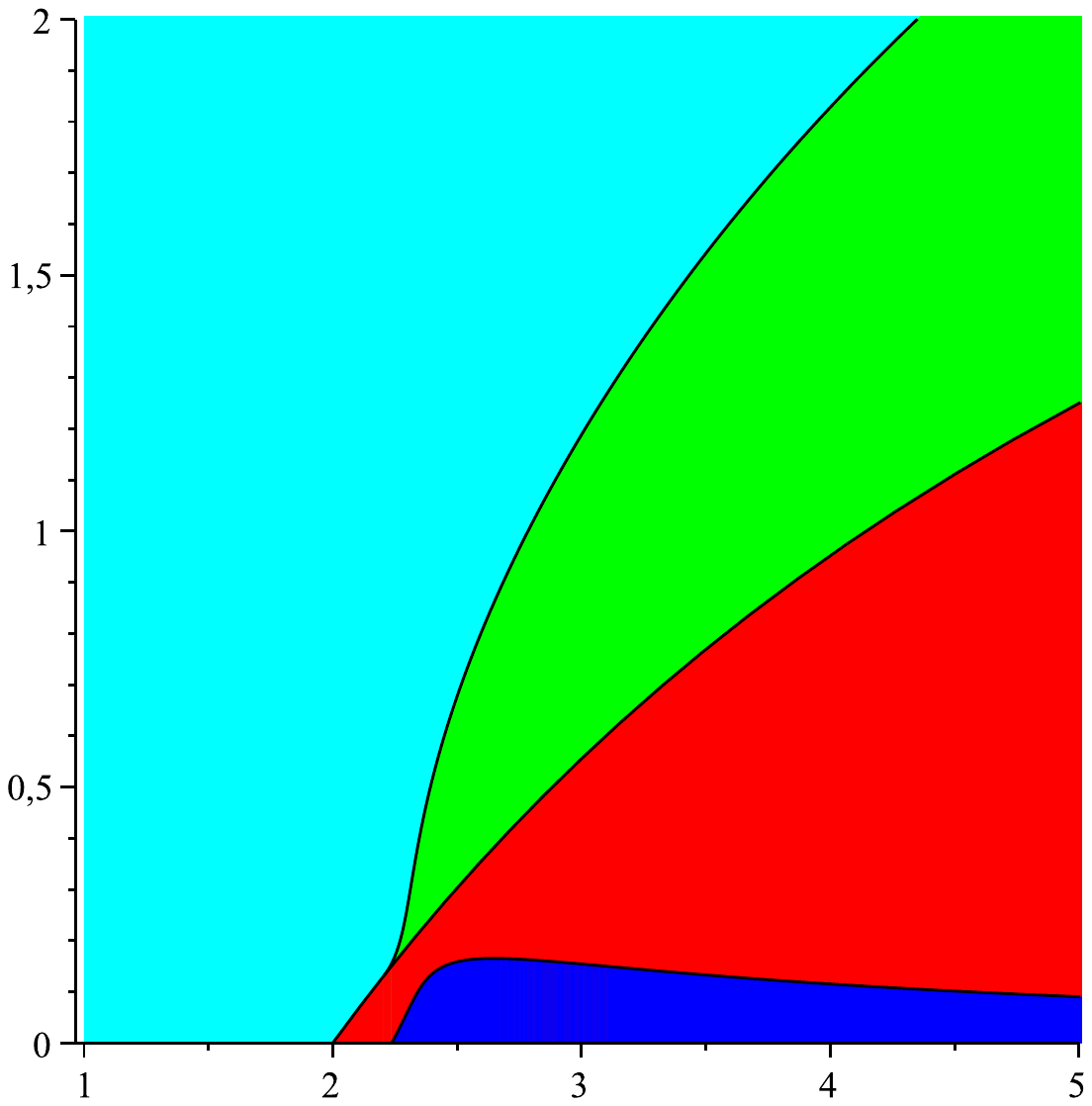}}}
\put(-2.4,5.95){{\sc $(a)$}}
\put(-4.5,5.85){{\sc $D$}}
\put(-3.3,5.1){{\sc $\mathcal{I}_0$}}
\put(-2,3.3){{\sc ${\color{white}\mathcal{I}_1}$}}
\put(0.05,5){{\sc $\Gamma_u$}}
\put(0.05,1.35){{\sc $S_{in}$}}
\put(2.7,5.95){{\sc $(b)$}}
\put(0.6,5.85){{\sc $D$}}
\put(5.15,5.45){{\sc $\Gamma_{\sLP}$}}
\put(5.15,3.8){{\sc $\Gamma_{\sBP}$}}
\put(1.4,4.9){{\sc $\mathcal{I}_0$}}
\put(3.6,4.4){{\sc $\mathcal{I}_2$}}
\put(3.4,2.4){{\sc ${\color{white}\mathcal{I}_1}$}}
\put(5.15,1.35){{\sc $S_{in}$}}
\put(7.7,5.95){{\sc $(c)$}}
\put(5.8,5.85){{\sc $D$}}
\put(9.4,5.88){{\sc $\Gamma_{\sLP}$}}
\put(10.3,4.1){{\sc $\Gamma_{\sBP}$}}
\put(6.8,4.5){{\sc $\mathcal{I}_0$}}
\put(9,4.2){{\sc $\mathcal{I}_2$}}
\put(9.3,2.6){{\sc ${\color{white}\mathcal{I}_1}$}}
\put(7.7,1.5){{\sc ${\color{white}\mathcal{I}_3}$}}
\put(7.15,1.5){{\sc ${\color{white}\Gamma_{\sH}^1}$}}
\put(8.5,1.7){{\sc ${\color{white}\Gamma_{\sH}^2}$}}
\put(10.3,1.4){{\sc $S_{in}$}}
\end{picture}
\end{center}
\vspace{-1.8cm}
\caption{Operating diagram of (\ref{ModFlocGen}). (a) The case considered in Section \ref{SubSec-DOCasFig13}. (b) The case considered in Section \ref{SubSec-DOCasFig12Siam}. (c) Magnification of (b) showing the curve $\Gamma_{\sH}$.}\label{Fig-CasFig13Siam}
\end{figure}

From Table \ref{TabExisStabGen} providing the existence and local stability conditions of steady states, we can state the next result determining theoretically the operating diagram in the case of Table \ref{TabParamVal} (line 1) where the various functions and the corresponding curves are defined in Tables \ref{TabFuncDO} and \ref{TabCurvDO}, respectively.
\begin{Proposition}                                 \label{PropRegionDO}
For the specific growth rates $f_1$ and $f_2$ defined in (\ref{SpeciFunc}) and the set of the biological parameter values in Table \ref{TabParamVal} (line 1),
we have $\lambda_u(D)<\lambda_v(D)$ for all $D\in\left[0,\overline{D}_u\right[$.
In addition, the existence and the local stability of the steady states $E_0$ and $E_1$ of model (\ref{ModFlocGen}) in the two regions $\mathcal{I}_0$ and $\mathcal{I}_1$ of the operating diagram shown in Fig. \ref{Fig-CasFig13Siam}(a) are described in Table \ref{Tab-DO_CasFig13Siam}.
\end{Proposition}
\begin{table}[!h]
\caption{Existence and local stability of steady states according to the regions in the operating diagram of Fig. \ref{Fig-CasFig13Siam}(a). The letter S [resp. U] means stable [resp. unstable]. No letter means that the steady state does not exist.} \label{Tab-DO_CasFig13Siam}
\begin{center}
\begin{tabular}{ @{\hspace{1mm}}l@{\hspace{2mm}} @{\hspace{2mm}}l@{\hspace{2mm}} @{\hspace{2mm}}c@{\hspace{2mm}}
                 @{\hspace{2mm}}l@{\hspace{2mm}} @{\hspace{2mm}}l@{\hspace{2mm}}  }
\hline\hline\\[-3mm]
 Condition                 &     Region       &       Color     &  $E_0$  &  $E_1$   \\[0.5mm] \hline \\[-3mm]
 $S_{in}< \lambda_u(D)$    & $\mathcal{I}_0$  &       Cyan      &    S    &              \\
 $S_{in}> \lambda_u(D)$    & $\mathcal{I}_1$  &       Red       &    I    &    S          \\
\end{tabular}
\end{center}
\end{table}

Note that the cyan region $\mathcal{I}_0$ in the operating diagram in Fig. \ref{Fig-CasFig13Siam}(a) corresponds to the washout of isolated and attached bacteria
while the red region $\mathcal{I}_1$ corresponds to the coexistence of both species around a steady state.
The theoretical study of the operating diagram determines the asymptotic behavior of the solutions for the set of biological parameters in \cite[Fig. 13]{FekihSIADS2019}.
Note that this figure illustrates the trajectories over time which converge towards the positive steady state $E_1$ for various initial conditions when the point of the plane $(S_{in},D)=(5,3.5)$ belongs to region $\mathcal{I}_1$.
Thus, this operating diagram presents a global vision of the behavior of the solutions according to the two operating parameters.
\subsection{A case with multiplicity of positive steady states and Hopf bifurcations}  \label{SubSec-DOCasFig12Siam}
In this section, we consider a case where there are LP and Hopf bifurcations. In this example, we can have two positive steady states, one being stable and the other unstable, and in addition, the stable one can be destabilized through a Hopf bifurcation.
For this purpose, we consider the biological parameter values that were used in \cite[Fig. 12]{FekihSIADS2019} (see Table \ref{TabParamVal}, line 2).
However, $S_{in}$ and $D$ are variable and not fixed as in \cite{FekihSIADS2019} where $D=0.1$.
With this set of parameters, we have $\lambda_v(D)<\lambda_{\sBP}(D)=\lambda_b(D)$ for all $D\in \left[0,\overline{D}_v\right[$ so that the function $S\mapsto H(S)$ is defined and decreasing on $\left]\lambda_v(D),\lambda_b(D)\right]$.
Using Tables \ref{TabExisStabGen}, \ref{TabFuncDO} and \ref{TabCurvDO}, we can state the next proposition determining theoretically the operating diagram. 
\begin{Proposition}
For the specific growth rates $f_1$ and $f_2$ defined in (\ref{SpeciFunc}) and the set of the biological parameter values in Table \ref{TabParamVal} (line 2),
we have $\lambda_v(D)<\lambda_{\sBP}(D)$ for all $D\in [0,\overline{D}_v[$.
In addition, the existence and the local stability of the steady states $E_0$, $E_1^1$ and $E_1^2$ of model (\ref{ModFlocGen}) in the four regions $\mathcal{I}_i$, $i=0,\ldots,3$ of the operating diagram shown in Fig. \ref{Fig-CasFig13Siam}(b-c) are described in Table \ref{Tab-DO-LP-Fig12Siam}.
\end{Proposition}
\begin{table}[!h]
\caption{Existence and local stability of steady states according to the regions in the operating diagram of Fig. \ref{Fig-CasFig13Siam}(b-c).} \label{Tab-DO-LP-Fig12Siam}
\begin{center}
\begin{tabular}{ @{\hspace{1mm}}l@{\hspace{2mm}} @{\hspace{2mm}}l@{\hspace{2mm}} @{\hspace{2mm}}c@{\hspace{2mm}} @{\hspace{2mm}}l@{\hspace{2mm}}
                 @{\hspace{2mm}}l@{\hspace{2mm}} @{\hspace{2mm}}l@{\hspace{2mm}} @{\hspace{2mm}}c@{\hspace{1mm}} }
\hline\hline\\[-3mm]
 Condition 1                                        &     Condition 2      &     Region       &       Color     &  $E_0$  &  $E_1^1$  &  $E_1^2$ \\[0.5mm] \hline \\[-3mm]
 $S_{in}<\min(\lambda_{\sLP}(D),\lambda_{\sBP}(D))$ &                      & $\mathcal{I}_0$  &       Cyan      &    S    &         &            \\
 $\lambda_{\sBP}(D)<S_{in}$                         & $c_4(S_{in}, D) > 0$ & $\mathcal{I}_1$  &       Red       &    I    &   S     &              \\
 $\lambda_{\sLP}(D)<S_{in}<\lambda_{\sBP}(D)$       & $c_4(S_{in}, D) > 0$ & $\mathcal{I}_2$  &       Green     &    S    &   S     &   I         \\
 $\lambda_{\sBP}(D)<S_{in}$                         & $c_4(S_{in}, D) < 0$ & $\mathcal{I}_3$  &       Blue      &    I    &   I     &                \\
\hline\hline
\end{tabular}
\end{center}
\end{table}

Note that the construction of the operating diagram in Fig. \ref{Fig-CasFig13Siam}(b-c) is obtained by plotting the various curves $\Gamma_i$, $i=\left\{{\rm{\scM{BP}}},{\rm{\scM{LP}}},{\rm{\scM{H}}}\right\}$ defined in Tables \ref{TabCurvDO}.
They correspond to the existence and stability conditions of all steady states provided in Table \ref{TabExisStabGen}.
The green region $\mathcal{I}_2$ corresponds to the bistability with either the coexistence around a steady state or the washout of the isolated and attached bacteria according to the initial condition.
The blue region $\mathcal{I}_3$ corresponds to
the instability of the positive steady state $E_1^1$ where there can be coexistence around a stable limit cycle.

Note that $c_4$ is a function of $S$, that is, $c_4=c_4(S)$ because it depends on the three state variables $S$, $u=U(S)$ and $v=V (S)$ defined in (\ref{eqUV}).
Moreover, $c_4$ is a function of $(S_{in},D)$, that is, $c_4=c_4(S_{in},D)$ because we can determine $S$ from the equation $D (S_{in}-S)=H(S)$.
With the set of parameters in Table \ref{TabParamVal} (line 2), we provide numerical evidence in \ref{Sec-AppedixCasFig12} of the change of sign of the function $c_4(S)$ on the existence interval of the positive steady state $E_1^1$ according to $D$.
Indeed, for $D<D_{\sH}^{max} \approx 0.165$ and fixed, the function $c_4(S)$ changes sign in $]\lambda_v(D),\lambda_{\sBP}(D)[$ so that the equation $c_4(S)=0$ has two solutions noted by
$$
S_{\sH}^2(D) < S_{\sH}^1(D) < \lambda_{\sBP}(D).
$$
For all $D<D_{\sH}^{max}$, we define the two solutions of the equation $c_4(S_{in},D)=0$ by the following critical values of $S_{in}$ which corresponds to a Hopf bifurcation
$$
S_{in} =  S_{in}^{Hi}(D) = \frac{1}{D} H\left(S_{\sH}^i(D)\right) + S_{\sH}^i(D), \quad  i=1,2.
$$
Consequently, the $\Gamma_{\sH}$ curve of the equation $c_4(S_{in},D)=0$ is given by the union of the two curves $\Gamma_{\sH}^1$ (on the left of the maximum) of equation $S_{in} = S_{in}^{H1}(D)$ and $\Gamma_{\sH}^2$ (on the right of the maximum) of equation $S_{in} = S_{in}^{H2}(D)$, see Fig. \ref{Fig-CasFig13Siam}(c).
\begin{Remark}
In \ref{Sec-AppedixCasArima}, we establish the operating diagram with the parameter set in \cite{FekihArima2020}. It is similar to that in Fig. \ref{Fig-CasFig13Siam}(b-c) where we find the same regions in Table \ref{Tab-DO-LP-Fig12Siam} (just $\lambda_{\sBP}(D)$ is equal to $\lambda_u(D)$ instead of $\lambda_b(D)$). However, the region $\mathcal{I}_3$ of destabilization of the positive steady state $E_1^1$ with the appearance of a stable limit cycle was not detected in \cite{FekihArima2020} because of the order of magnitude of $D_{\sH}^{max}$ as demonstrated in \ref{Sec-AppedixCasArima}.
\end{Remark}
\subsection{Another case with multiplicity of positive steady states and Hopf bifurcations}   \label{SubSec-DOCasFig6}
\begin{figure}[!h]
\setlength{\unitlength}{1.0cm}
\begin{center}
\begin{picture}(6.2,6.2)(0,0)
\put(-5.4,0){\rotatebox{0}{\includegraphics[width=6cm,height=7cm]{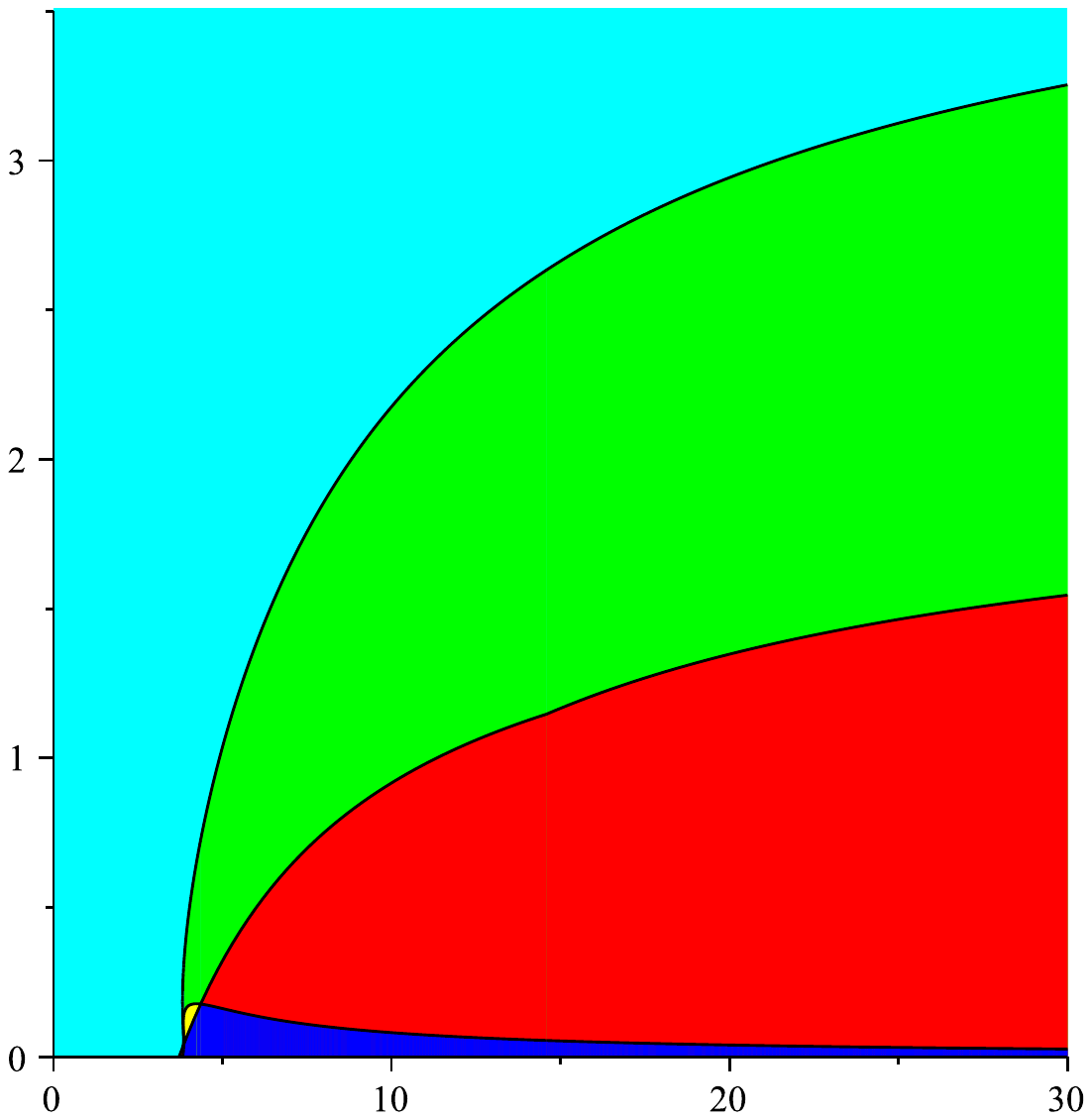}}}
\put(0,0){\rotatebox{0}{\includegraphics[width=6cm,height=7.1cm]{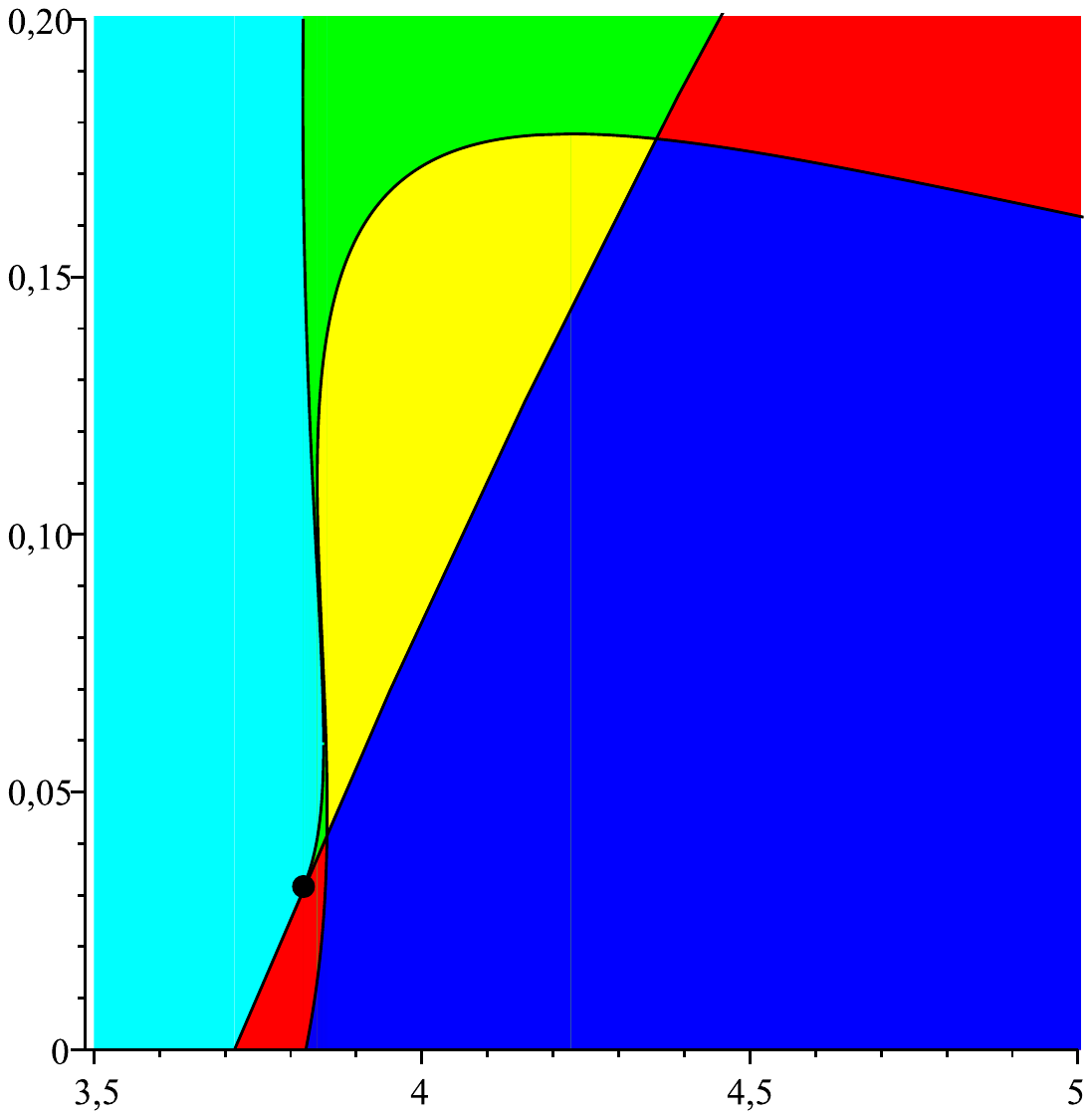}}}
\put(5.4,0){\rotatebox{0}{\includegraphics[width=6cm,height=7.1cm]{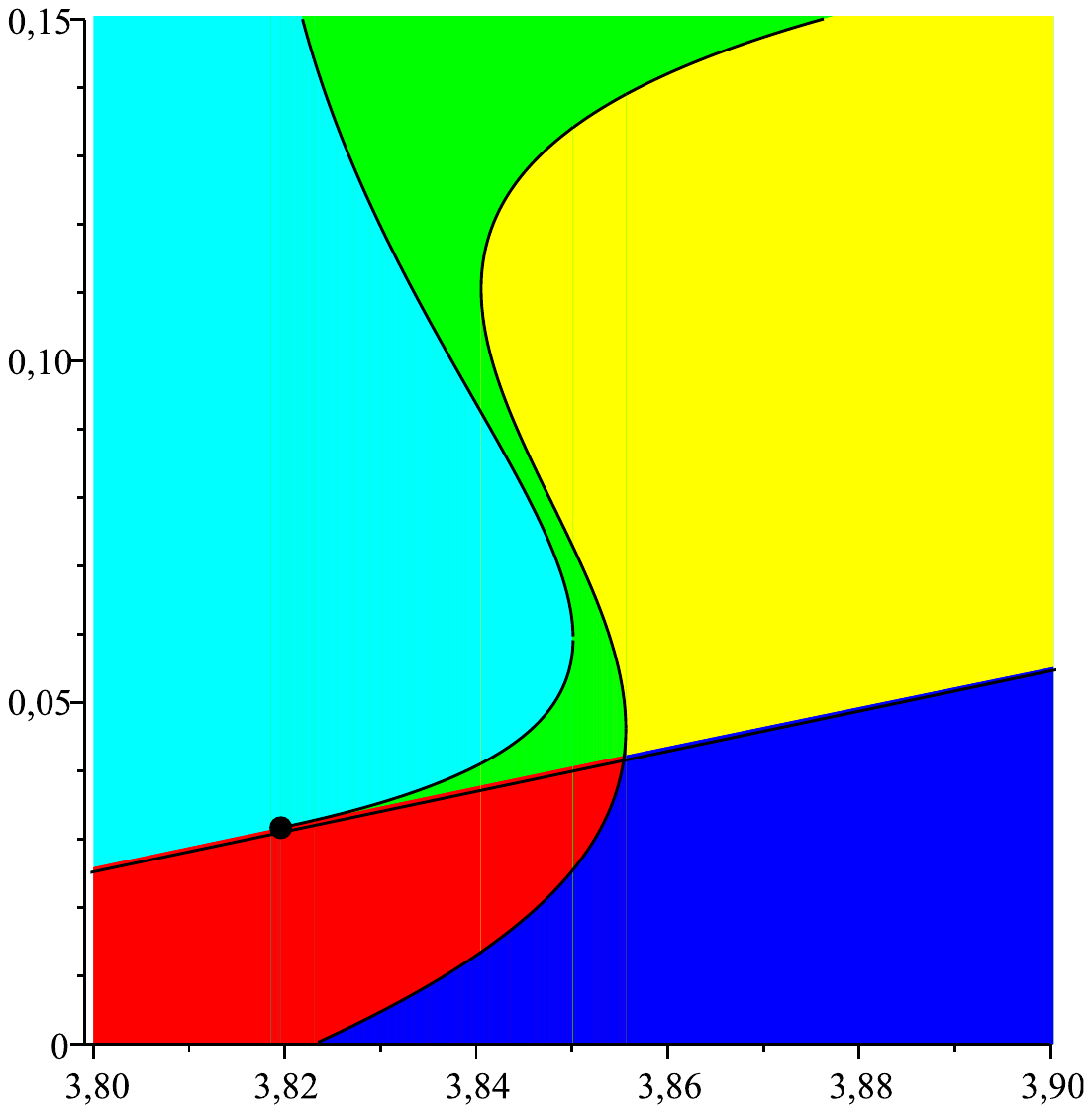}}}
\put(-2.4,6){{\sc $(a)$}}
\put(-4.7,5.95){{\sc $D$}}
\put(0.1,5.45){{\sc $\Gamma_{\sLP}$}}
\put(0.1,3.3){{\sc $\Gamma_b$}}
\put(-3.75,2.35){{\sc $\Gamma_u$}}
\put(-3.8,1.6){{\sc ${\color{white}\Gamma_{\sH}}$}}
\put(-3.7,5){{\sc $\mathcal{I}_0$}}
\put(-1.5,4){{\sc $\mathcal{I}_2$}}
\put(-1.8,2.3){{\sc ${\color{white}\mathcal{I}_1}$}}
\put(-2.9,1.9){{\sc ${\color{white}\mathcal{I}_3}$}}
\put(-2.9,1.86){\tiny {\color{white}\vector(-1,-1){0.4}}}
\put(0.1,1.4){{\sc $S_{in}$}}
\put(3.1,6){{\sc $(b)$}}
\put(0.8,5.95){{\sc $D$}}
\put(1.8,5.95){{\sc $\Gamma_{\sLP}$}}
\put(3.75,5.95){{\sc $\Gamma_u$}}
\put(5.55,4.95){{\sc $\Gamma_{\sH}$}}
\put(1.2,3.8){{\sc $\mathcal{I}_0$}}
\put(2.6,5.55){{\sc $\mathcal{I}_2$}}
\put(2.4,4.3){{\sc $\mathcal{I}_4$}}
\put(4.6,5.45){{\sc ${\color{white}\mathcal{I}_1}$}}
\put(3.9,3.2){{\sc ${\color{white}\mathcal{I}_3}$}}
\put(5.55,1.4){{\sc $S_{in}$}}
\put(8.3,6){{\sc $(c)$}}
\put(6.2,5.95){{\sc $D$}}
\put(7.15,5.95){{\sc $\Gamma_{\sLP}$}}
\put(9.6,5.95){{\sc $\Gamma_{\sH}$}}
\put(10.85,3){{\sc $\Gamma_u$}}
\put(7,3.4){{\sc $\mathcal{I}_0$}}
\put(8,5.45){{\sc $\mathcal{I}_2$}}
\put(9.4,4.3){{\sc $\mathcal{I}_4$}}
\put(7,1.8){{\sc ${\color{white}\mathcal{I}_1}$}}
\put(9.5,1.9){{\sc ${\color{white}\mathcal{I}_3}$}}
\put(10.9,1.4){{\sc $S_{in}$}}
\end{picture}
\end{center}
\vspace{-1.8cm}
\caption{Operating diagram of (\ref{ModFlocGen}). (a) The case considered in Section \ref{SubSec-DOCasFig6}. (b-c) Magnifications of (a) showing the curve $\Gamma_{\sH}$.}\label{Fig-CasFig6Siam}
\end{figure}

In the operating diagram of the case considered in Section \ref{SubSec-DOCasFig12Siam}, we do not have a region where the two positive steady states are both unstable, see Table \ref{Tab-DO-LP-Fig12Siam}.
Thus, the aim of this section is to provide an example where there is a new region (labeled $\mathcal{I}_4$, see Table \ref{Tab-DO1})
of instability of the two positive steady states. For this purpose, we consider the biological parameter values that were used in \cite[Fig. 6]{FekihSIADS2019} (see Table \ref{TabParamVal}, line 3).
However, $S_{in}$ and $D$ are variable and not fixed as in \cite{FekihSIADS2019} where $D=0.1$. Indeed,
the study in \cite{FekihSIADS2019} was limited to one parameter bifurcation diagrams according to $S_{in}$.
With this set of parameters, we have $\lambda_v(D)<\lambda_{\sBP}(D)$ so that the function $S\mapsto H(S)$ is defined and decreasing on $(\lambda_v(D),\lambda_{\sBP}(D)]$ as shown in Fig. \ref{Fig-Hc4-CasFig6Siam}.
In addition, the two curves $\Gamma_u$ and $\Gamma_b$ intersect when $(S_{in},D) = (14.588,1.147)$ so that $\Gamma_{\sBP}=\Gamma_u$ for all $D\in[0,1.147]$ and $\Gamma_{\sBP}=\Gamma_b$ for all $D\in\left[1.147,\overline{D}_b\right[$.
Similarly to the previous cases, we can state the next result.
\begin{Proposition}                                 \label{PropRegionDO}
For the specific growth rates $f_1$ and $f_2$ defined in (\ref{SpeciFunc}) and the set of the biological parameter values in Table \ref{TabParamVal} (line 3),
we have $\lambda_v(D)<\lambda_{\sBP}(D)$ for all $D\in [0,\overline{D}_v[$.
In addition, the existence and the local stability of the steady states $E_0$, $E_1^1$ and $E_1^2$  of model (\ref{ModFlocGen}) in the five regions $\mathcal{I}_k$, $k = 0,\ldots,4$ of the operating diagram shown in Fig. \ref{Fig-CasFig6Siam} are described in Table \ref{Tab-DO1}.
\end{Proposition}
\begin{table}[!h]
\caption{Existence and local stability of steady states according to the regions in the operating diagram of Fig. \ref{Fig-CasFig6Siam}.} \label{Tab-DO1}
\begin{center}
\begin{tabular}{ @{\hspace{1mm}}l@{\hspace{2mm}} @{\hspace{2mm}}l@{\hspace{2mm}} @{\hspace{2mm}}c@{\hspace{2mm}} @{\hspace{2mm}}l@{\hspace{2mm}}
                 @{\hspace{2mm}}l@{\hspace{2mm}} @{\hspace{2mm}}l@{\hspace{2mm}} @{\hspace{2mm}}c@{\hspace{1mm}} }
\hline\hline\\[-3mm]
 Condition 1                                        &    Condition 2       &     Region       &       Color     &  $E_0$  &  $E_1^1$  &  $E_1^2$ \\[0.5mm] \hline \\[-3mm]
 $S_{in}<\min(\lambda_{\sLP}(D),\lambda_{\sBP}(D))$ &                      & $\mathcal{I}_0$  &       Cyan      &    S    &         &            \\
 $\lambda_{\sBP}(D)<S_{in}$                         & $c_4(S_{in}, D) > 0$ & $\mathcal{I}_1$  &       Red       &    I    &   S     &             \\
 $\lambda_{\sLP}(D)<S_{in}<\lambda_{\sBP}(D)$       & $c_4(S_{in}, D) > 0$ & $\mathcal{I}_2$  &       Green     &    S    &   S     &   I          \\
 $\lambda_{\sBP}(D)<S_{in}$                         & $c_4(S_{in}, D) < 0$ & $\mathcal{I}_3$  &       Blue      &    I    &   I     &                \\
 $\lambda_{\sLP}(D)<S_{in}<\lambda_{\sBP}(D)$       & $c_4(S_{in}, D) < 0$ & $\mathcal{I}_4$  &       Yellow    &    S    &   I     &    I            \\
\hline\hline
\end{tabular}
\end{center}
\end{table}

Fig. \ref{Fig-CasFig6Siam}(a) illustrates the operating diagram of model (\ref{ModFlocGen}) while Figs. \ref{Fig-CasFig6Siam}(b-c) illustrate magnifications of regions $\mathcal{I}_1$, $\mathcal{I}_3$ and $\mathcal{I}_4$.
The operating diagram in Fig. \ref{Fig-CasFig6Siam} is divided into five regions.
The blue region $\mathcal{I}_3$ corresponds to the instability of the positive steady state.
The yellow region $\mathcal{I}_4$ corresponds to the instability of the two positive steady states $E_1^1$ and $E_1^2$ where the system can exhibit bistability
with either coexistence around a stable limit cycle or the washout of the isolated and attached bacteria.

In Fig. \ref{Fig-Hc4-CasFig6Siam}, we give the justification that the operating diagram is the one shown in Fig. \ref{Fig-CasFig6Siam}.
Indeed, it illustrates the functions $H(S)$ and $c_4(S)$ for $D$ fixed at $D^\ast=0.1$ to see the change of the sign of $c_4(S)$.
The solutions $S_{\sH}^1$ and $S_{\sH}^2$ of the equation $c_4(S)=0$ correspond to the critical values $S_{in}^{H1}$ and $S_{in}^{H2}$ which are the intersections of the horizontal line of equation $D=D^\ast$ in the $(S_{in},D)$-plane of the operating diagram in Fig. \ref{Fig-CasFig6Siam}.
\section{Operating diagrams and bifurcations diagrams in MATCONT}            \label{Sec-NumDO}
In this section, we use MATCONT \cite{MATCONT} to numerically analyze the one- and two-parameter diagrams of model (\ref{ModFlocGen}) and to detect two-parameter bifurcations that cannot be established theoretically.
It also allows us to validate our theoretical results.
In fact, MATCONT is a MATLAB numerical continuation package used to analyze the different types of bifurcations of the continuous and discrete parameterized systems of ODEs.
It allows to trace the trajectories over time according to the initial condition and the bifurcation diagrams with a single parameter or two parameters.
More precisely, it allows one to visualize the curves of steady states according to a parameter by determining their local asymptotic behavior thanks to the calculation of the eigenvalues of the Jacobian matrix evaluated at the steady state. Moreover, it also allows one to determine the stable or unstable limit cycles by calculation of the sign of the First Lyapunov coefficient.
Thanks to test functions, MATCONT detects all types of bifurcations such as the transcritical bifurcation or Branch Points (BP), saddle-node or Limit Points (LP) bifurcation, Cusp (CP) bifurcation, Hopf (H) bifurcation, Limit Point of Cycles (LPC) or fold bifurcation points of limit cycles, period doubling bifurcation points of limit cycles.

From these critical bifurcation points, MATCONT can determine the various curves in the operating diagram according to two parameters by numerical continuation.
These curves of objects of a given type (e.g. steady states, limit cycle, Hopf bifurcation points, homoclinic orbits, etc.) are calculated under variation of one or more system parameters.
The reader is addressed to the relevant paper of Dhooge et al. \cite{DhoogeMCMDS2008} for more on this interesting subject.
\subsection{Operating diagram in the case considered in Section \ref{SubSec-DOCasFig12Siam}, obtained with MATCONT}
\begin{figure}[!h]
\setlength{\unitlength}{1.0cm}
\begin{center}
\begin{picture}(5.7,4.4)(0,0)
\put(-5.2,0){\rotatebox{0}{\includegraphics[width=7cm,height=11.6cm]{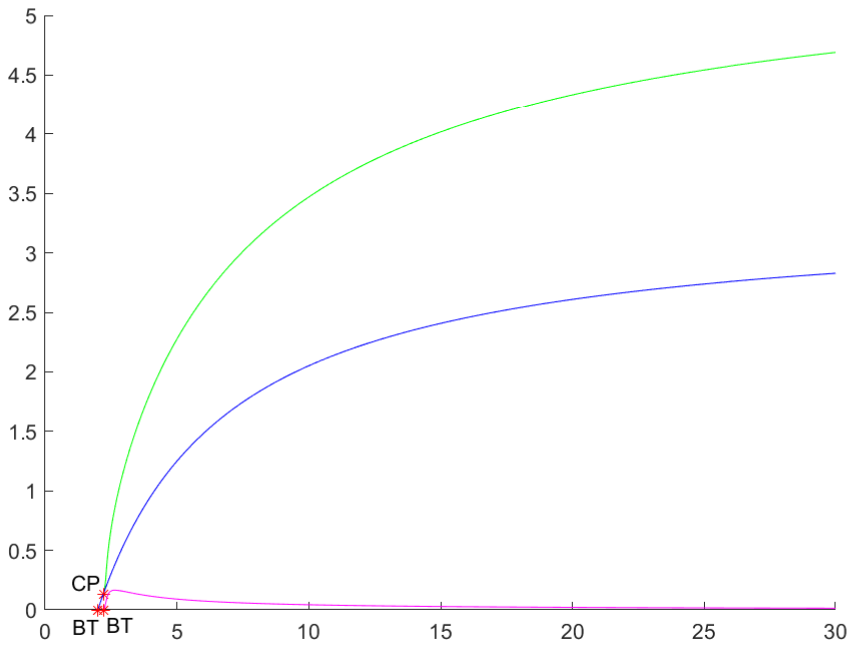}}}
\put(0,0){\rotatebox{0}{\includegraphics[width=7cm,height=11.6cm]{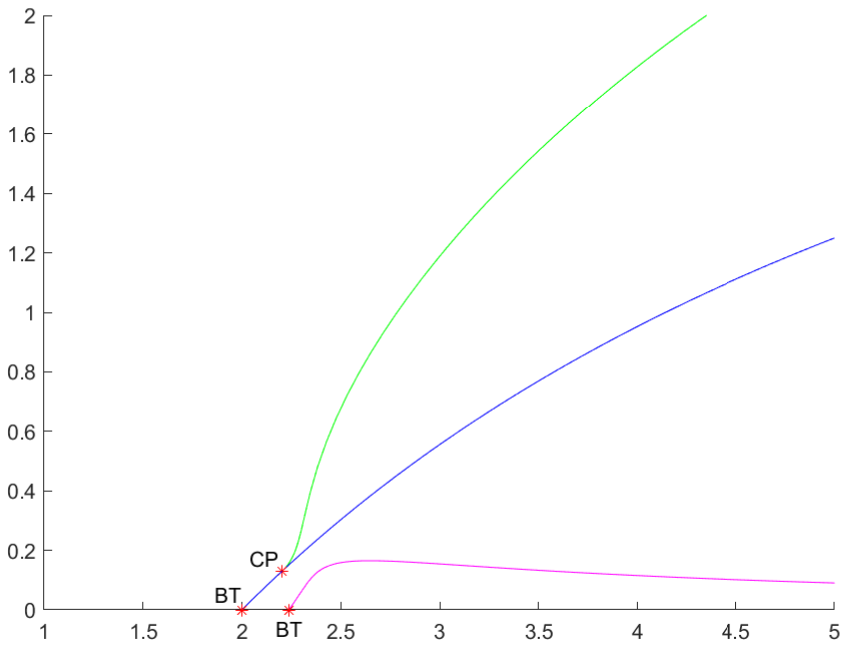}}}
\put(5.2,0){\rotatebox{0}{\includegraphics[width=7.7cm,height=12.8cm]{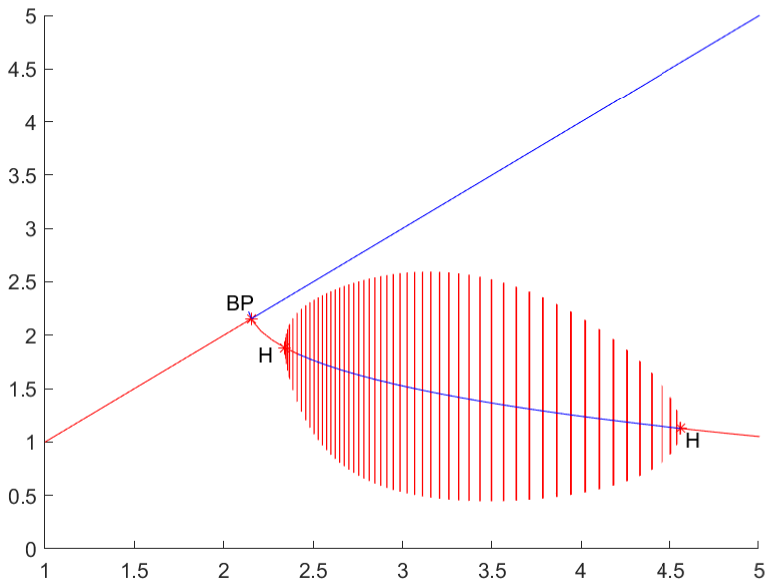}}}
\put(-2.4,4.8){{\sc $(a)$}}
\put(-4.4,4.7){{\sc $D$}}
\put(-0.1,4.4){{\sc ${\color{green}\Gamma_{\sLP}}$}}
\put(-0.1,2.85){{\sc ${\color{blue}\Gamma_b}$}}
\put(-3.4,0.65){{\sc ${\color{magenta}\Gamma_{\sH}}$}}
\put(-3.7,4){{\sc $\mathcal{I}_0$}}
\put(-1.8,3.4){{\sc $\mathcal{I}_2$}}
\put(-1.8,1.6){{\sc $\mathcal{I}_1$}}
\put(-3.6,1.05){{\sc $\mathcal{I}_3$}}
\put(-3.6,1){\tiny \vector(-1,-1){0.4}}
\put(-0.1,0.55){{\sc $S_{in}$}}
\put(2.7,4.8){{\sc $(b)$}}
\put(0.8,4.7){{\sc $D$}}
\put(4.4,4.7){{\sc ${\color{green}\Gamma_{\sLP}}$}}
\put(5.1,3.1){{\sc ${\color{blue}\Gamma_b}$}}
\put(2.7,0.95){{\sc ${\color{magenta}\Gamma_{\sH}}$}}
\put(1.7,3.6){{\sc $\mathcal{I}_0$}}
\put(3.6,3.1){{\sc $\mathcal{I}_2$}}
\put(3.8,1.4){{\sc $\mathcal{I}_1$}}
\put(2.3,0.65){{\sc $\mathcal{I}_3$}}
\put(5.1,0.55){{\sc $S_{in}$}}
\put(7.9,4.8){{\sc $(c)$}}
\put(6,4.7){{\sc $S$}}
\put(6.5,2.05){{\sc ${\color{red}E_0}$}}
\put(8.6,3.65){{\sc ${\color{blue}E_0}$}}
\put(7.22,2.3){{\sc ${\color{red}E_1}$}}
\put(8.3,1.8){{\sc ${\color{blue}E_1}$}}
\put(9.95,1.52){{\sc ${\color{red}E_1}$}}
\put(10.25,0.55){{\sc $S_{in}$}}
\end{picture}
\end{center}
\vspace{-1cm}
\caption{MATCONT: (a) operating diagram of (\ref{ModFlocGen}) in the case considered in Section \ref{SubSec-DOCasFig12Siam}. (b) Magnification of (a) showing the curve $\Gamma_{\sH}$. (c) The corresponding one-parameter bifurcation diagram in variable $S$ when $D=0.1$.}\label{FigDO-Fig12Siam-Num}
\end{figure}

In this section, we determine the one and two-parameter bifurcation diagrams in Fig. \ref{FigDO-Fig12Siam-Num} using MATCONT
for the set of the biological parameter values in \cite[Fig. 12]{FekihSIADS2019} where the one bifurcation diagram is obtained using SCILAB.
The corresponding set of the parameters are provided in Table \ref{TabParamVal} (line 2).
The intersection point between $\Gamma_b$ and $\Gamma_{\sLP}$ is a two parameters bifurcation of type Cusp (CP) while the
intersection points between $\Gamma_{\sH}$ and $\Gamma_b$ with the $D=0$ axis is of type Bogdanov-Takens (BT).
These types of bifurcation are not detected in the theoretical study of the operating diagram obtained in Section \ref{SubSec-DOCasFig12Siam}.
Table \ref{TabCriParmStatNFC2} summarizes the critical operating parameters, the state, and the normal form coefficient for BT and CP bifurcations.
\begin{table}[!h]
\caption{Operating parameters, state, normal form coefficient values $(\theta_1, \theta_2)$ [resp. $\theta_3$] for BT [resp. CP], at the bifurcation points in Fig. \ref{FigDO-Fig12Siam-Num}. The abbreviation BP [resp. CP] means a Bogdanov-Takens point [resp. Cusp] bifurcation.}  \label{TabCriParmStatNFC2}
\vspace{-0.2cm}
\begin{center}\begin{tabular}{ @{\hspace{1mm}}l@{\hspace{2mm}} @{\hspace{2mm}}l@{\hspace{2mm}} @{\hspace{2mm}}l@{\hspace{1mm}} @{\hspace{2mm}}l@{\hspace{1mm}} }			
\hline
Bifurcation       & Parameter $(S_{in},D)$ &  State $(S,u,v)$  &  Normal form coefficient
\\ \hline
BT                & (2,0)                  & (2,0,0)           &   $(\theta_1, \theta_2)=(1.03\,10^{-6}, -0.6)$
\\
CP                & (2.204,0.130)          & (2.204,0,0)       &  $\theta_3=-0.523$
\\
BT                & (2.236,0)              & (2,0,0)           &   $(\theta_1, \theta_2)$ impossible
\\ \hline
\end{tabular}\end{center}
\end{table}
Fig. \ref{FigDO-Fig12Siam-Num}(c) illustrates the one-parameter bifurcation diagram in variable $S$ when $D$ is fixed at $D=0.1$ in the case considered in Section \ref{SubSec-DOCasFig12Siam}. It reveals the appearance and the disappearance of stable limit cycles via two Hopf bifurcations.
\subsection{Operating diagram in the case considered in Section \ref{SubSec-DOCasFig6}, obtained with MATCONT}          \label{Sec-OD-MATCONT}
\begin{figure}[!h]
\setlength{\unitlength}{1.0cm}
\begin{center}
\begin{picture}(6.2,5)(0,0)
\put(-5.5,0){\rotatebox{0}{\includegraphics[width=6.5cm,height=10.8cm]{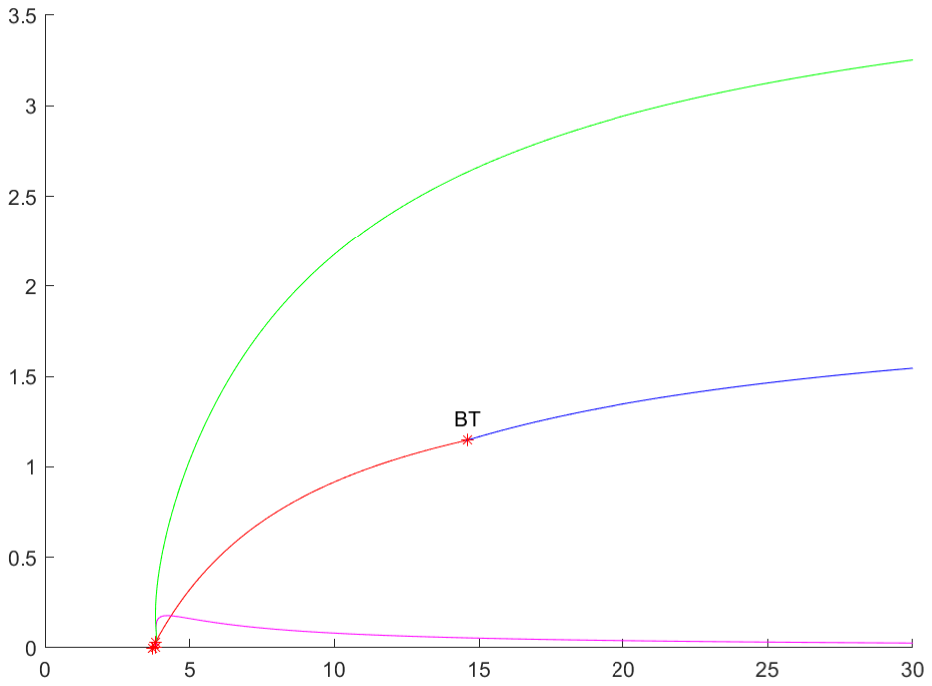}}}
\put(-0.6,-0.65){\rotatebox{0}{\includegraphics[width=7.2cm,height=14.6cm]{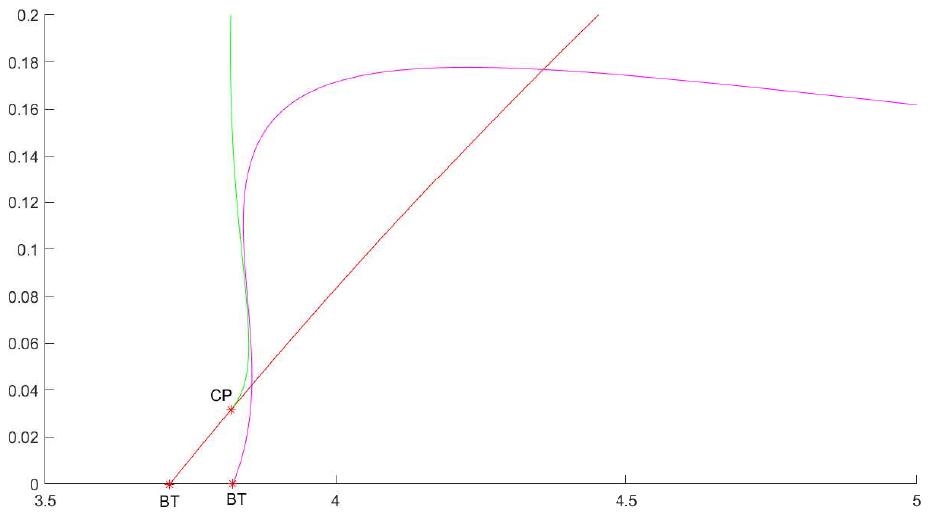}}}
\put(5.8,0){\rotatebox{0}{\includegraphics[width=6.5cm,height=10.8cm]{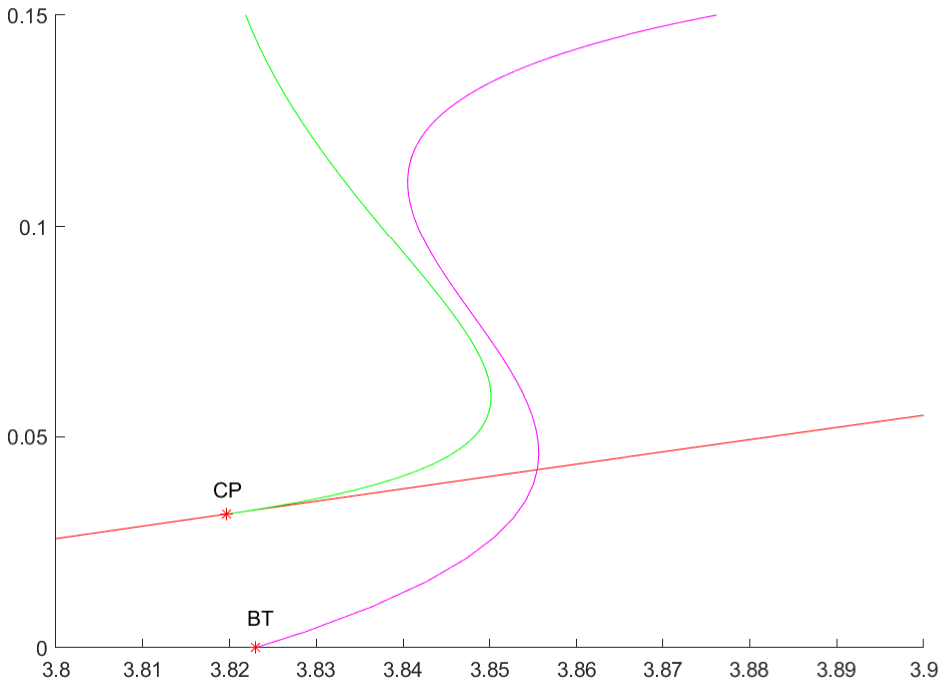}}}
\put(-2.7,4.8){{\sc $(a)$}}
\put(-4.7,4.7){{\sc $D$}}
\put(-0.3,4.35){{\sc ${\color{green}\Gamma_{\sLP}}$}}
\put(-0.3,2.35){{\sc ${\color{blue}\Gamma_b}$}}
\put(-3.5,1.7){{\sc ${\color{red}\Gamma_u}$}}
\put(-2.3,0.7){{\sc ${\color{magenta}\Gamma_{\sH}}$}}
\put(-3.9,4){{\sc $\mathcal{I}_0$}}
\put(-1.8,3){{\sc $\mathcal{I}_2$}}
\put(-1.8,1.4){{\sc $\mathcal{I}_1$}}
\put(-3,1.05){{\sc $\mathcal{I}_3$}}
\put(-3,1){\tiny \vector(-1,-1){0.4}}
\put(-0.3,0.55){{\sc $S_{in}$}}
\put(2.9,4.8){{\sc $(b)$}}
\put(0.75,4.7){{\sc $D$}}
\put(1.7,4.72){{\sc ${\color{green}\Gamma_{\sLP}}$}}
\put(3.7,4.75){{\sc ${\color{red}\Gamma_u}$}}
\put(5.6,3.8){{\sc ${\color{magenta}\Gamma_{\sH}}$}}
\put(1,3){{\sc $\mathcal{I}_0$}}
\put(2.4,4.45){{\sc $\mathcal{I}_2$}}
\put(2.1,3.1){{\sc $\mathcal{I}_4$}}
\put(4.6,4.3){{\sc $\mathcal{I}_1$}}
\put(3.7,2.2){{\sc $\mathcal{I}_3$}}
\put(5.6,0.55){{\sc $S_{in}$}}
\put(8.5,4.8){{\sc $(c)$}}
\put(6.6,4.7){{\sc $D$}}
\put(7.5,4.7){{\sc ${\color{green}\Gamma_{\sLP}}$}}
\put(9.95,4.65){{\sc ${\color{magenta}\Gamma_{\sH}}$}}
\put(11,2){{\sc ${\color{red}\Gamma_u}$}}
\put(7.4,2.6){{\sc $\mathcal{I}_0$}}
\put(8.1,4.3){{\sc $\mathcal{I}_2$}}
\put(9.6,3.3){{\sc $\mathcal{I}_4$}}
\put(7.55,1){{\sc $\mathcal{I}_1$}}
\put(9.6,1.2){{\sc $\mathcal{I}_3$}}
\put(11,0.55){{\sc $S_{in}$}}
\end{picture}
\end{center}
\vspace{-0.9cm}
\caption{MATCONT: (a) operating diagram of (\ref{ModFlocGen}) in the case considered in Section \ref{SubSec-DOCasFig6}. (b-c) Magnifications of (a) showing the curve $\Gamma_{\sH}$.}\label{FigDO-Fig6Siam-Num}
\end{figure}

Fig. \ref{FigDO-Fig6Siam-Num} illustrates the operating diagram obtained numerically using MATCONT. It is identical to the operating diagram obtained theoretically in Fig. \ref{Fig-CasFig6Siam}.
However, MATCONT detects the nature of bifurcations at the intersection points between the curves $\Gamma_u$ and $\Gamma_{\sLP}$ which is of Cusp (CP) type and
between the curves $\Gamma_u$ and $\Gamma_H$ with the $S_{in}$-axis which are Bogdanov-Takens points (BT).
The critical operating parameters, the state, and the normal form coefficient for BT and CP bifurcations are summarized in Table \ref{TabCriParmStatNFC}.
\begin{table}[!h]
\caption{Operating parameters, state, normal form coefficient values $(\theta_1, \theta_2)$ [resp. $\theta_3$] for BT [resp. CP], at the bifurcation points in Fig. \ref{FigDO-Fig6Siam-Num}. The abbreviation BP [resp. CP] means a Bogdanov-Takens point [resp. Cusp] bifurcation.}  \label{TabCriParmStatNFC}
\vspace{-0.2cm}
\begin{center}\begin{tabular}{ @{\hspace{1mm}}l@{\hspace{2mm}} @{\hspace{2mm}}l@{\hspace{2mm}} @{\hspace{2mm}}l@{\hspace{1mm}} @{\hspace{2mm}}l@{\hspace{1mm}} }			
\hline
Bifurcation       & Parameter $(S_{in},D)$   &  State $(S,u,v)$  &  Normal form coefficient
\\ \hline
BT                & (3.714,0)       & (3.714,0,0)      &   $(\theta_1, \theta_2)=\left(3.23\,10^{-8}, -0.306\right)$
\\
CP                & (3.819,0.032)   & (3.819,0,0)      &  $\theta_3=-0.483$
\\
BT                & (3.823,0)       & (3.714,0,0)     &   $(\theta_1, \theta_2)=\left(2.38\,10^{-6}, -0.306\right)$
\\
BT                & (14.588,1.147)  & (14.588,0,0)     &   $(\theta_1, \theta_2)=(2.019, -1.127)$
\\ \hline
\end{tabular}\end{center}
\end{table}
\subsection{Bifurcation diagram with respect to $S_{in}$, corresponding to $D=0.1$ in the case considered in Section \ref{SubSec-DOCasFig6}}                \label{Sec-AnalBif}
In what follows, we will analyze the various types of bifurcation by crossing one region to another in the operating diagram of Fig. \ref{Fig-CasFig6Siam} or Fig. \ref{FigDO-Fig6Siam-Num}.
Using Prop. \ref{PropRegionDO}, the nature of all the bifurcations by passing through the various curves $\Gamma_i$ defined in Table \ref{TabCurvDO} is described in the following result.
\begin{Proposition}                                  \label{PropNatDB}
Let $f_1$ and $f_2$ be the specific growth rates defined in (\ref{SpeciFunc}). Let the set of the biological parameter values be in Table \ref{TabParamVal} (line 3).
The nature of all the bifurcations of model (\ref{ModFlocGen}) by crossing the different regions of the operating diagram in Fig. \ref{Fig-CasFig6Siam} is provided in Table \ref{TabDOBif}.
\end{Proposition}
\begin{table}[ht]
\caption{Nature of all the bifurcations of system (\ref{ModFlocGen}) by passing the different curves $\Gamma_i$, $i=\left\{ {\rm{\scM{BP}}},{\rm{\scM{LP}}},{\rm{\scM{H}}} \right\}$ defined in Table \ref{TabCurvDO}.} \label{TabDOBif}
\vspace{-0.2cm}
\begin{center}	
\begin{tabular}{@{\hspace{2mm}}l@{\hspace{2mm}} @{\hspace{2mm}}l@{\hspace{2mm}} @{\hspace{2mm}}l@{\hspace{1mm}} @{\hspace{2mm}}l@{\hspace{2mm}}}			
            Transition              &     Curve       &  Bifurcation  &  Steady states             \\ \hline
$\mathcal{I}_0$ to $\mathcal{I}_2$  & $\Gamma_{\sLP}$ &      LP       & $E_1^1=E_1^2$   \\
$\mathcal{I}_0$ to $\mathcal{I}_1$  &   $\Gamma_u$   &      BP        & $E_0=E_1^1$    \\
$\mathcal{I}_2$ to $\mathcal{I}_1$  & $\Gamma_{\rm{\scM{BP}}}$& BP      & $E_1^2=E_0$    \\
$\mathcal{I}_2$ to $\mathcal{I}_4$  &   $\Gamma_{\sH}$   &      H         & $E_1^1$  \\
$\mathcal{I}_1$ to $\mathcal{I}_3$  &   $\Gamma_{\sH}$   &      H         & $E_1^1$  \\
$\mathcal{I}_4$ to $\mathcal{I}_3$  &   $\Gamma_u$   &      BP        & $E_1^2=E_0$
\end{tabular}
\end{center}
\end{table}

\begin{figure}[tbhp]
\setlength{\unitlength}{1.0cm}
\begin{center}
\begin{picture}(7,5.5)(0,0)
\put(-5.5,0){\rotatebox{0}{\includegraphics[width=6.8cm,height=12cm]{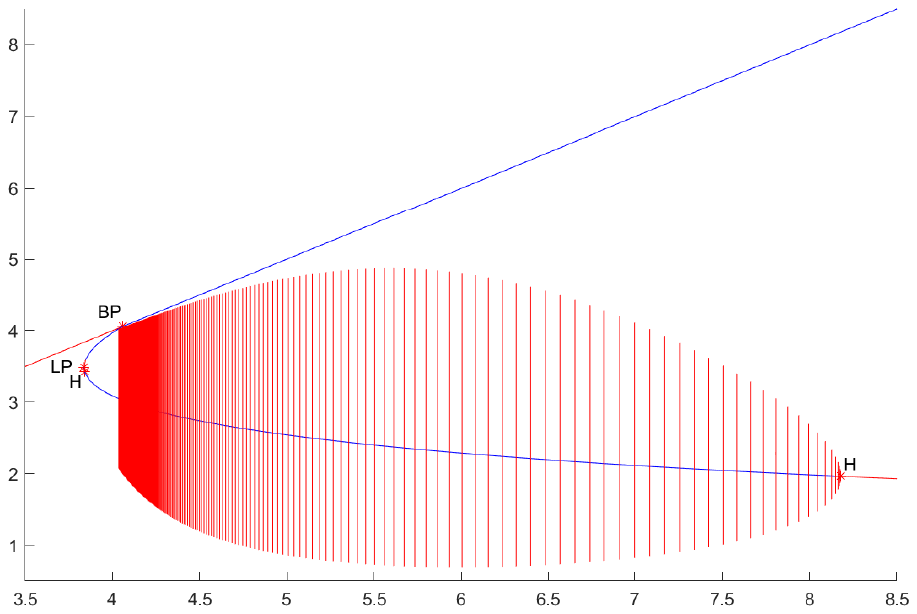}}}
\put(0,0){\rotatebox{0}{\includegraphics[width=6.8cm,height=12cm]   {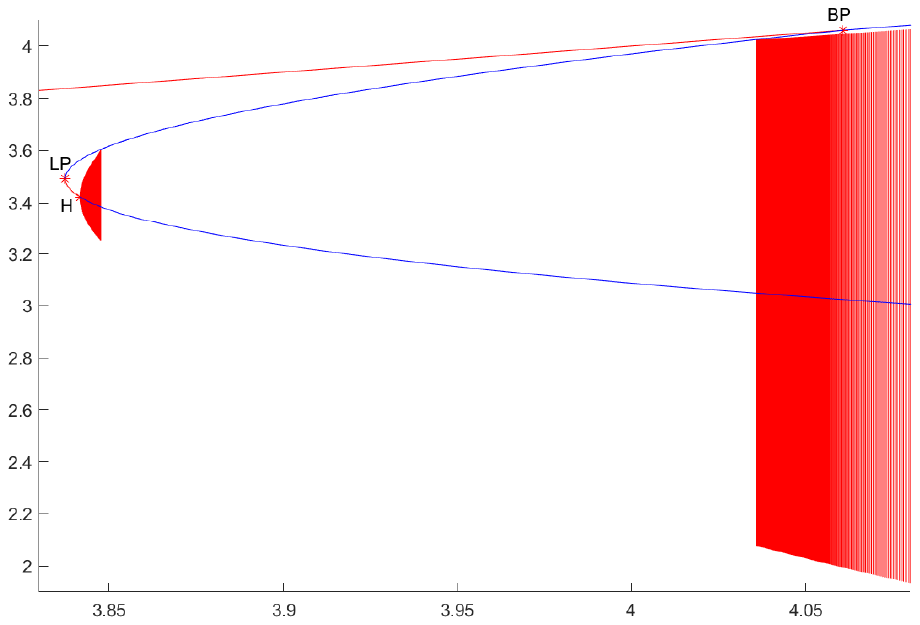}}}
\put(5.5,0){\rotatebox{0}{\includegraphics[width=6.8cm,height=12cm] {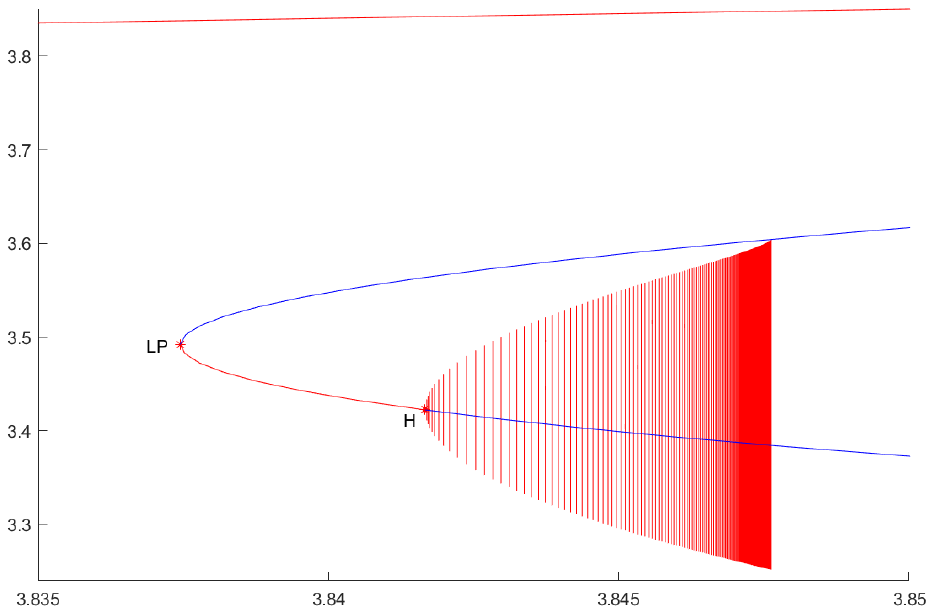}}}
\put(-2.3,5.2){{\sc  $(a)$}}
\put(-4.25,5.2){{\sc  $S$}}
\put(-4.3,2.85){{\sc ${\color{red}E_0}$}}
\put(-1.5,4.4){{\sc ${\color{blue}E_0}$}}
\put(0.15,1.9){{\sc ${\color{red}E_1^1}$}}
\put(0.4,1.1){{\sc  $S_{in}$}}
\put(3.5,5.2){{\sc  $(b)$}}
\put(1.25,5.2){{\sc  $S$}}
\put(2.5,5){{\sc ${\color{red}E_0}$}}
\put(5.6,5.3){{\sc ${\color{blue}E_0}$}}
\put(3,4.55){{\sc ${\color{blue}E_1^2}$}}
\put(3,3.3){{\sc ${\color{blue}E_1^1}$}}
\put(5.9,1.1){{\sc  $S_{in}$}}
\put(9,5.35){{\sc  $(c)$}}
\put(6.7,5.3){{\sc $S$}}
\put(8.5,5){{\sc ${\color{red}E_0}$}}
\put(9,3.45){{\sc ${\color{blue}E_1^2}$}}
\put(7.9,2.3){{\sc ${\color{red}E_1^1}$}}
\put(10.8,2.1){{\sc ${\color{blue}E_1^1}$}}
\put(11.4,1.1){{\sc  $S_{in}$}}
\end{picture}
\end{center}
\vspace{-1.3cm}
\caption{MATCONT: one parameter bifurcation diagram of (\ref{ModFlocGen}) in variable $S$ in the case considered in Section \ref{SubSec-DOCasFig6}.
(b) a magnification of two homoclinic bifurcations when $S_{in}\in [3.83,4.08]$; (c) a magnification of supercritical Hopf bifurcation when $S_{in}\in [3.835,3.85]$.} \label{FigDB-CasFig6Siam-MC}
\end{figure}

Let $D$ be fixed at $D=D^\ast=0.1$.
Next, we analyze the one-parameter bifurcation diagram with respect to $S_{in}$ as the bifurcating parameter to show the nature of bifurcations by crossing various boundaries between the different regions in the operating diagram. Note that the one-parameter bifurcation diagram in $D$ can be obtained in the same way.
Using MATCONT, we illustrate in Fig. \ref{FigDB-CasFig6Siam-MC} the one-parameter bifurcation diagram in $S_{in}$, with $S$ on the $y$-axis.
Similarly, we can obtain the one-parameter bifurcation diagram for the concentrations of isolated and attached bacteria, $u$ and $v$, respectively.
Note that the two-parameter bifurcation diagram does not show the disappearance of the limit cycle like the one-parameter bifurcation diagram.

In the following, we present the step-by-step approach to obtain the one-parameter bifurcation diagram in the variable $S$ using MATCONT.
Increasing $S_{in}$ from zero, the bifurcation diagram in Fig. \ref{FigDB-CasFig6Siam-MC}(a) illustrates the BP bifurcation occurring at $S_{in}=\sigma_5\approx 4.061$ between $E_0$ and $E_1^2$.
Increasing $S_{in}$ further, the washout steady state $E_0$ changes stability and becomes unstable (see Fig. \ref{FigDB-CasFig6Siam-MC}(a-b)).
Starting from this BP bifurcation and counting backward, $E_1^2$ emerges at $S_{in}=\sigma_5$ unstable by decreasing $S_{in}$.
A first close-up is illustrated in Fig. \ref{FigDB-CasFig6Siam-MC}(b) and a second close-up is illustrated in Fig. \ref{FigDB-CasFig6Siam-MC}(c).

Next, there is a LP between $E_1^1$ and $E_1^2$ at $S_{in}=\sigma_1\approx 3.837$ when $S\approx3.492$, $u\approx996\,10^{-5}$ and $v\approx107\,10^{-5}$ so that these two interior steady states disappear by decreasing $S_{in}$ further. Inversely, increasing $S_{in}$ from LP, $E_1^1$ and $E_1^2$ appear LES and unstable, respectively.
After that, increasing $S_{in}$ further, a Hopf (H) bifurcation occurs at $E_1^1$ when $S_{in}=\sigma_2\approx 3.842$, $S\approx3.422$, $u\approx 0.012$ and $v\approx 1.5\,10^{-3}$.
A stable limit cycle emerges through a supercritical Hopf bifurcation where the first Lyapunov coefficient is given by $-0.430$.
Moreover, $E_1^1$ changes stability and becomes unstable.
Increasing $S_{in}$ further, once again a Hopf bifurcation occurs at $E_1^1$ when $S_{in}=\sigma_6\approx 8.179$, $S\approx1.963$, $u\approx0.140$ and $v\approx0.139$. A stable limit cycle disappears through a supercritical Hopf bifurcation where the first Lyapunov coefficient is given by $-34\,10^{-3}$. Moreover, $E_1^1$ changes stability and becomes LES.
Fig. \ref{Fig3DHom2} shows the stable limit cycles in the three-dimensional space $(S,u,v)$ for different values of $S_{in}$ between $\sigma_4$ and $\sigma_6$.
Starting from the first Hopf bifurcation at $\sigma_2$ and increasing $S_{in}$, the radius of the stable limit cycle increases until his disappearance through a homoclinic bifurcation when $S_{in}=\sigma_3\approx 3.8477$. Fig. \ref{FigT-HomCasFig6Siam}(a) shows the period of the cycle tends to infinity when $S_{in}$ tends to $\sigma_3$.
Starting from the second Hopf bifurcation at $\sigma_6$ and decreasing $S_{in}$, the radius of the stable limit cycle first increases and then decreases until his disappearance through a homoclinic bifurcation when $S_{in}=\sigma_4\approx 4.03468$. Fig. \ref{FigT-HomCasFig6Siam}(b) shows the period of the cycle tends to infinity when $S_{in}$ tends to $\sigma_4$ confirming the homoclinic bifurcation.
The analysis of the one-parameter bifurcation diagram in $S_{in}$ from the operating diagram in Fig. \ref{Fig-CasFig6Siam} is summarized in the following result.
\begin{Proposition}                                     \label{PropBD}
For the specific growth rates $f$ and $g$ defined in (\ref{SpeciFunc}) and the set of the biological parameter values in Table \ref{TabParamVal} (see line 3),
the existence and the local stability of all steady states of (\ref{ModFlocGen}) according to $S_{in}$ are described in Table \ref{TableExiStabBD} when $D=0.1$ is fixed.
The critical values $\sigma_i$, $i=1,\ldots,6$ of different bifurcations according to the parameter $S_{in}$ and the corresponding nature are defined in Table \ref{TableSigmai}.
\end{Proposition}
\begin{table}[ht]
\caption{Existence and stability of all steady states of (\ref{ModFlocGen}) according to $S_{in}$ for the set of parameter in Fig. \ref{Fig-CasFig6Siam} when $D=0.1$. The critical values $\sigma_i$, $i=1,\ldots,6$ are defined in Table \ref{TableSigmai}.}\label{TableExiStabBD}
\vspace{-0.2cm}
\begin{center}
\begin{tabular}{l|llllll}
Interval of $S_{in}$   & $E_0$ & $E_1^1$ &  $E_1^2$ \\ \hline
$(0,\sigma_1)$         &     S           &                 &                 \\
$(\sigma_1,\sigma_2)$  &     S           &       S         &        U          \\
$(\sigma_2,\sigma_3)$  &     S           &       U         &        U           \\
$(\sigma_3,\sigma_4)$  &     S           &       U         &        U          \\
$(\sigma_4,\sigma_5)$  &     S           &       U         &        U           \\
$(\sigma_5,\sigma_6)$  &     U           &       U         &                    \\
$(\sigma_6,+\infty)$   &     U           &       S         &
\end{tabular}
\end{center}
\vspace{-0.6cm}
\end{table}
\begin{table}[ht]
\caption{Definitions of the critical values $\sigma_i$, $i=1,\ldots,6$ of $D$ and their corresponding nature of bifurcation when $D=0.1$ is fixed. The abbreviations Hom and T mean homoclinic and $T$ period of solutions, respectively.}\label{TableSigmai}
\vspace{-0.2cm}
\begin{center}
\begin{tabular}{lll}
Definition                                                           &  Value    &   Bifurcation                      \\
\hline
$\sigma_1=\lambda_{\sLP}(D)$                                         &  3.837    &    LP           \\
$\sigma_2$ is the first solution of equation  $c_4(S_{in})=0$        &  3.842    &    H            \\
$\sigma_3$: $\ds\lim_{S_{in}\rightarrow \sigma_3}T(S_{in})= +\infty$ &  3.84770  &    Hom         \\
$\sigma_4$  $\ds\lim_{S_{in}\rightarrow \sigma_4}T(S_{in})= +\infty$ &  4.03468  &    Hom         \\
$\sigma_5=\lambda_{\sBP}(D)$                            &  4.061    &    BP           \\
$\sigma_6$ is the second solution of equation $c_4(S_{in})=0$        &  8.179    &    H
\end{tabular}
\end{center}
\vspace{-0.6cm}
\end{table}
\section{Effect of flocculation on the operating diagram}      \label{Sec-EffectFloc}
\begin{figure}[!h]
\setlength{\unitlength}{1.0cm}
\begin{center}
\begin{picture}(5.2,5.1)(0,0)
\put(-4.9,0){\rotatebox{0}{\includegraphics[width=5cm,height=6cm]{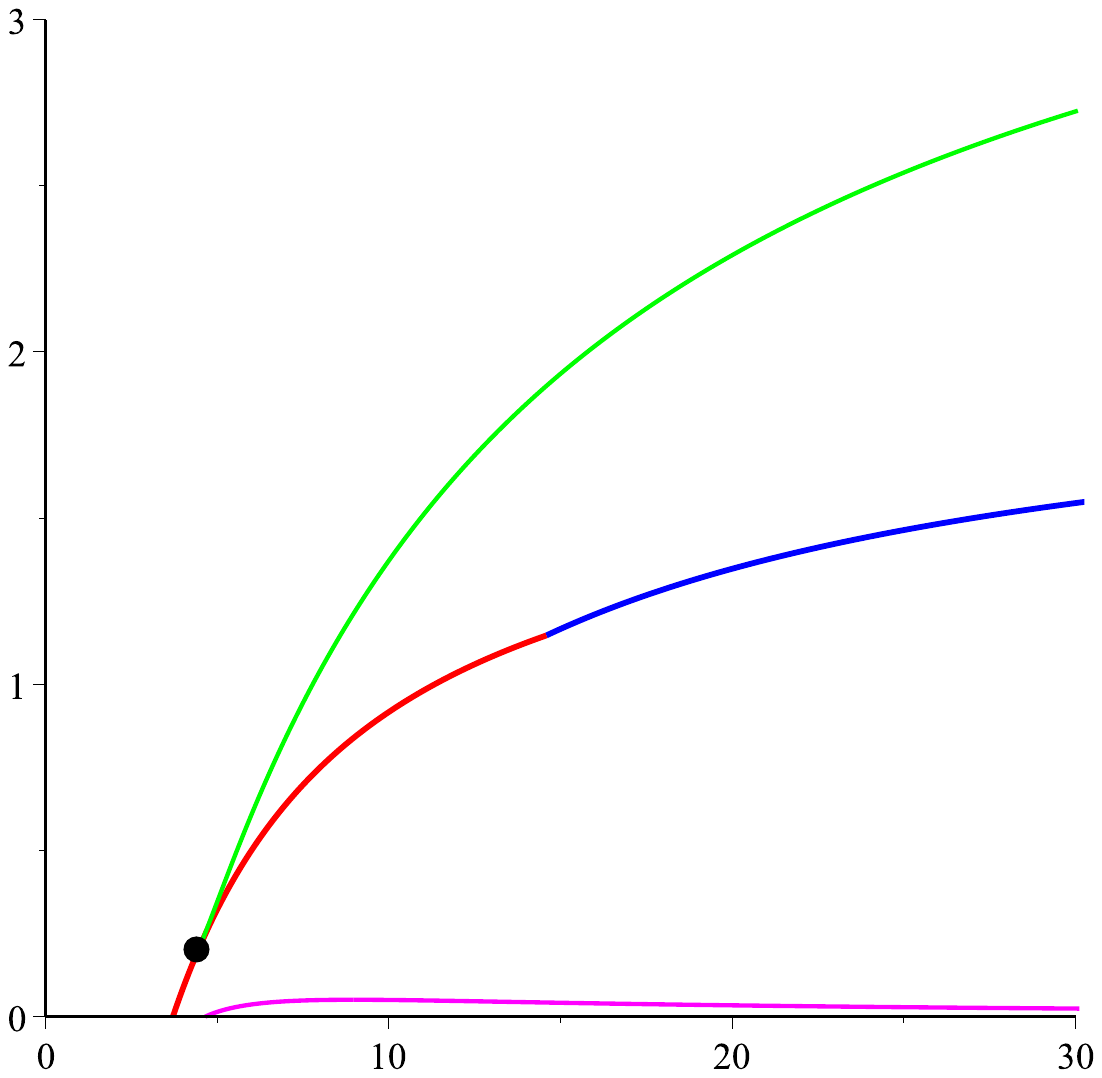}}}
\put(0,0){\rotatebox{0}{\includegraphics[width=5cm,height=6cm]{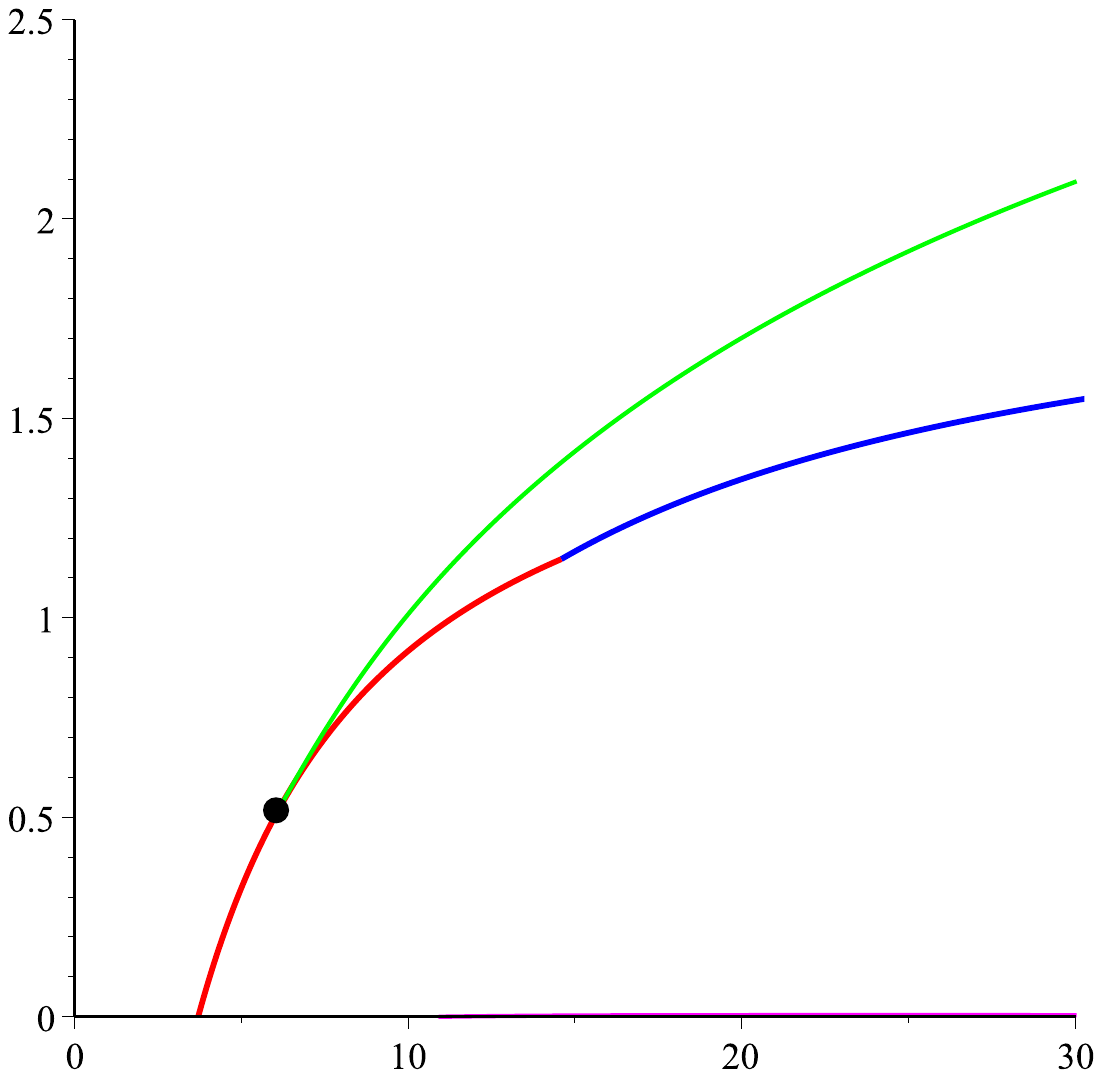}}}
\put(2.7,1.3){\rotatebox{0}{\includegraphics[width=2.2cm,height=2cm]{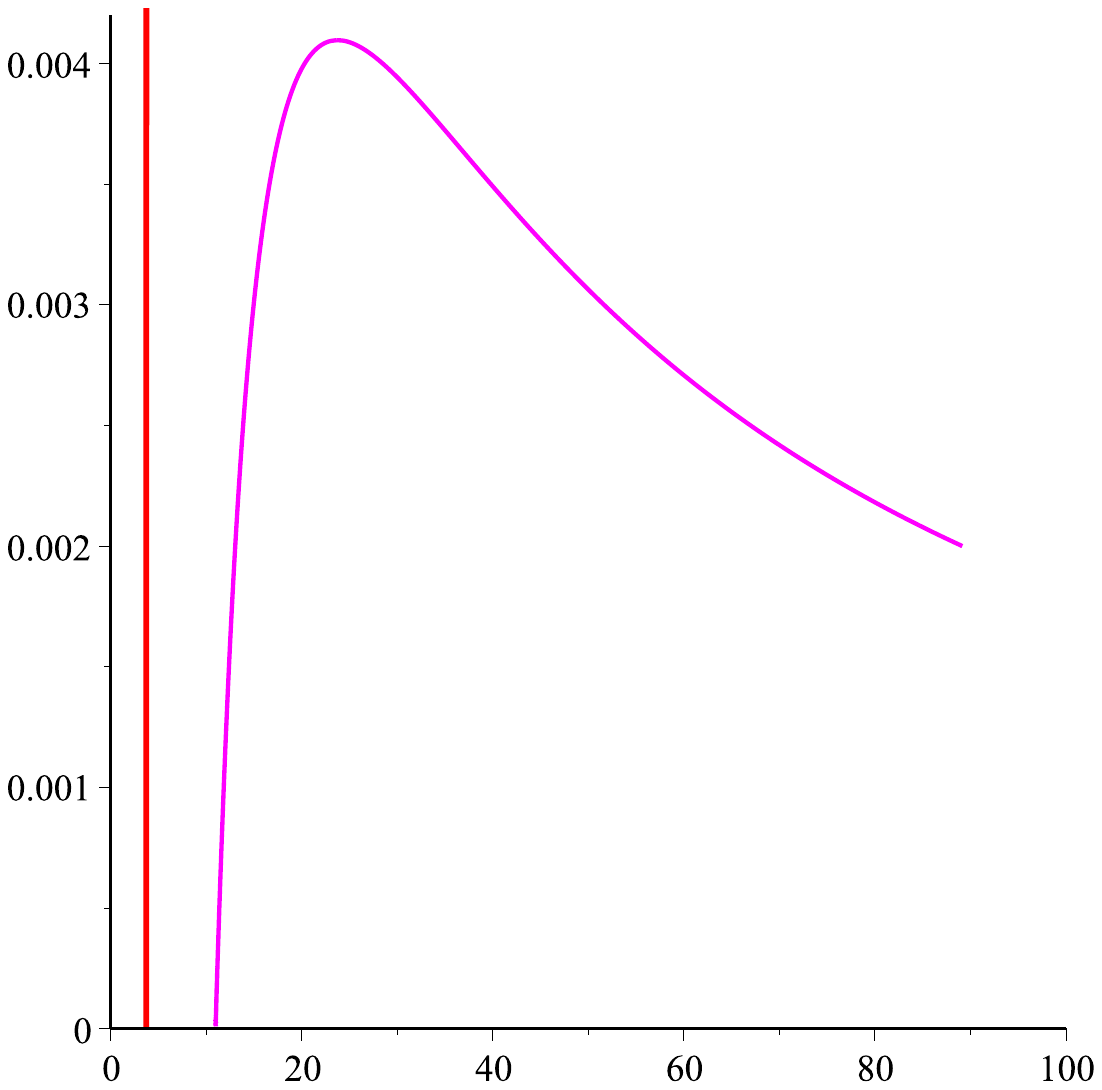}}}
\put(4.9,0){\rotatebox{0}{\includegraphics[width=5cm,height=6cm]{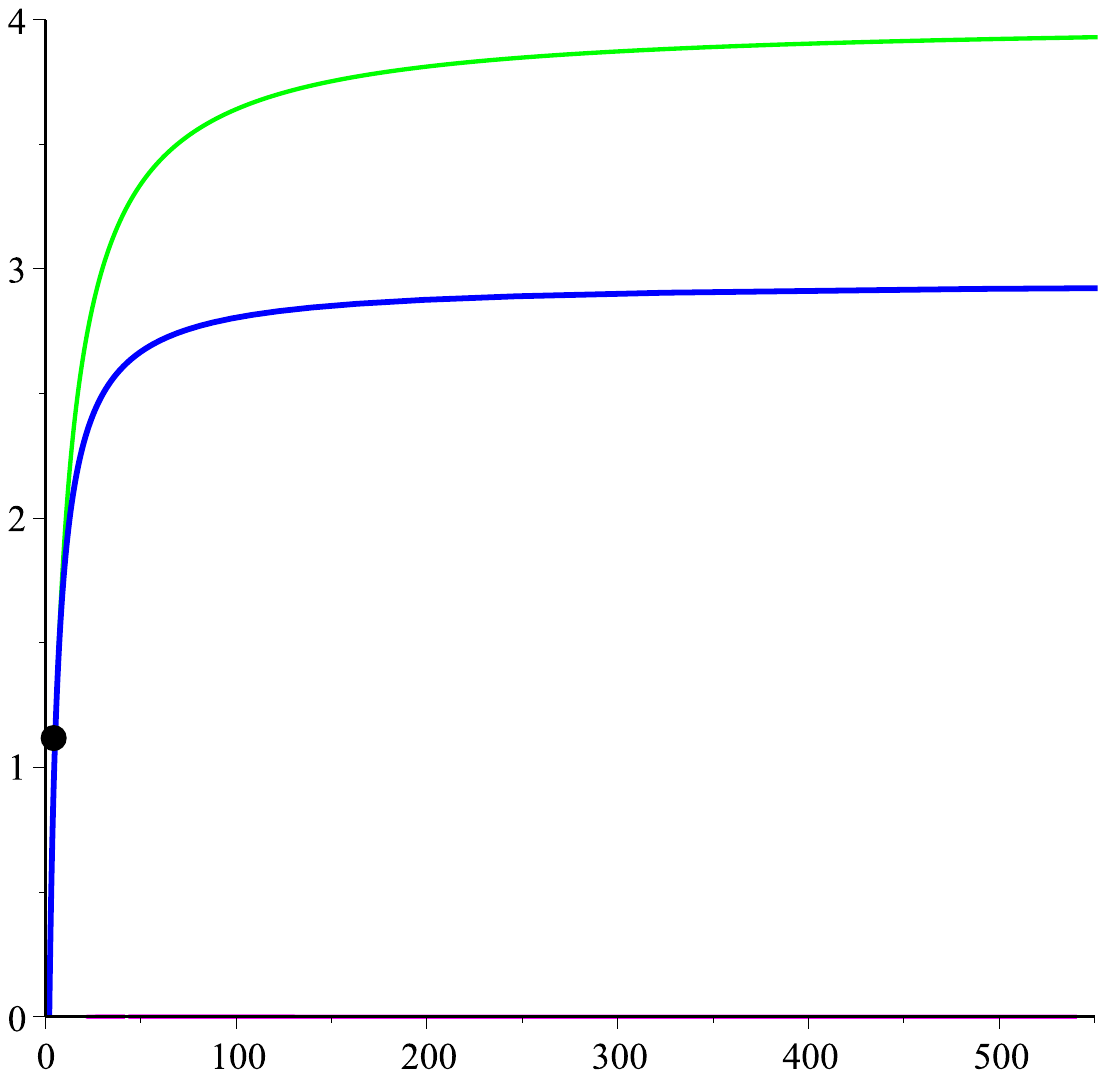}}}
\put(6.7,1.2){\rotatebox{0}{\includegraphics[width=3cm,height=2.7cm]{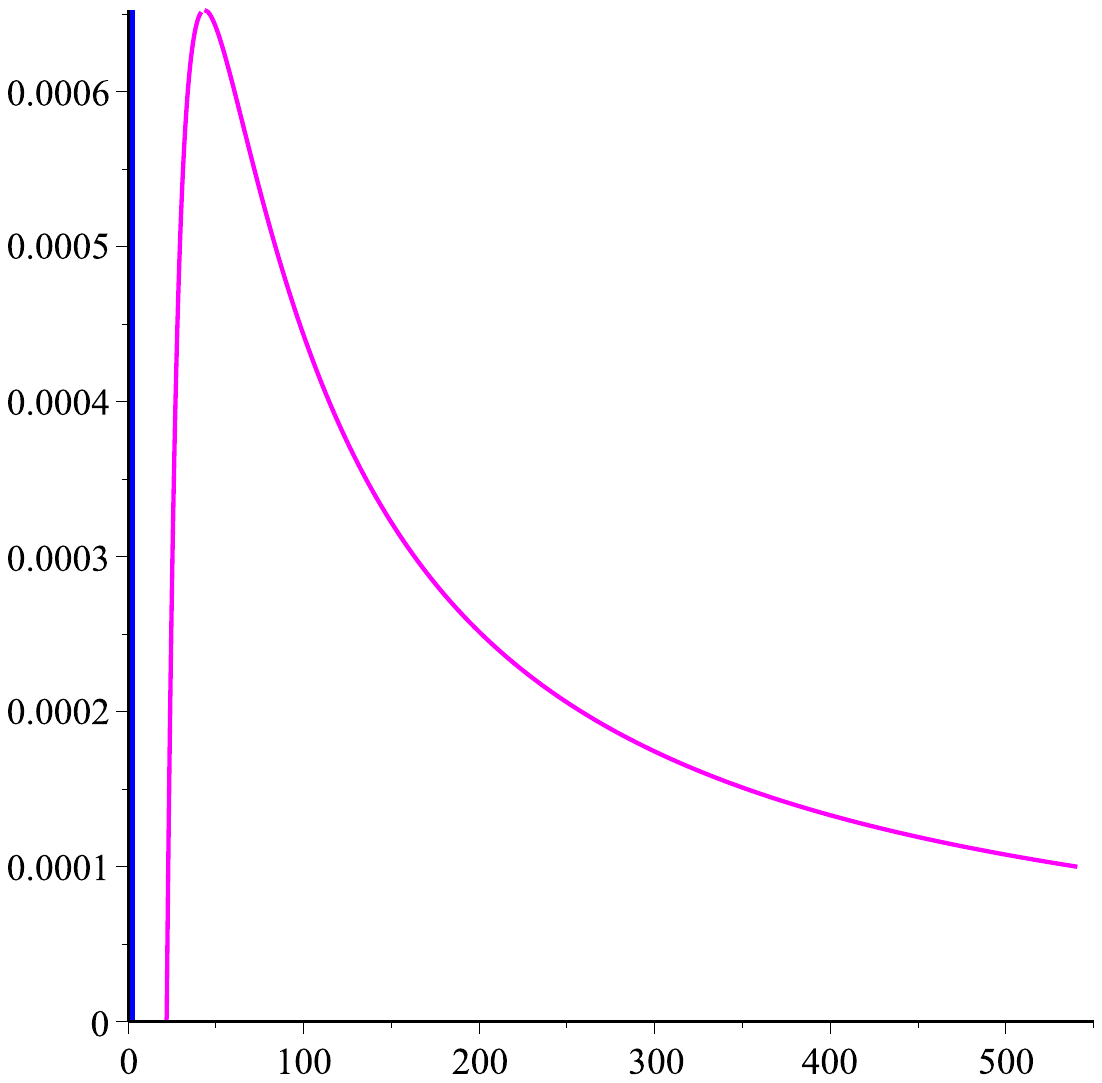}}}
\put(-4.9,-4.2){\rotatebox{0}{\includegraphics[width=5cm,height=5.8cm]{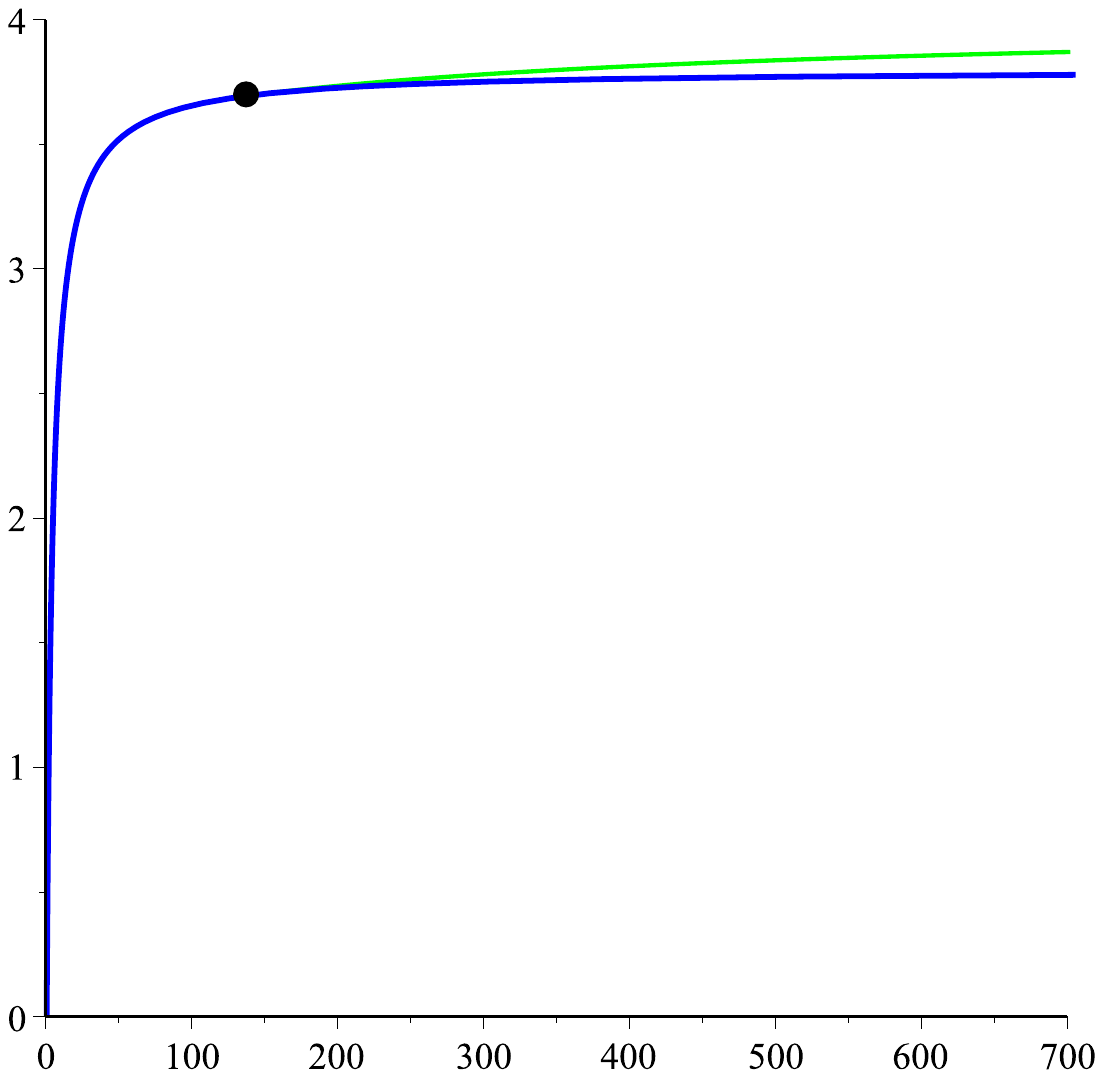}}}
\put(0,-4.2){\rotatebox{0}{\includegraphics[width=5cm,height=5.8cm]{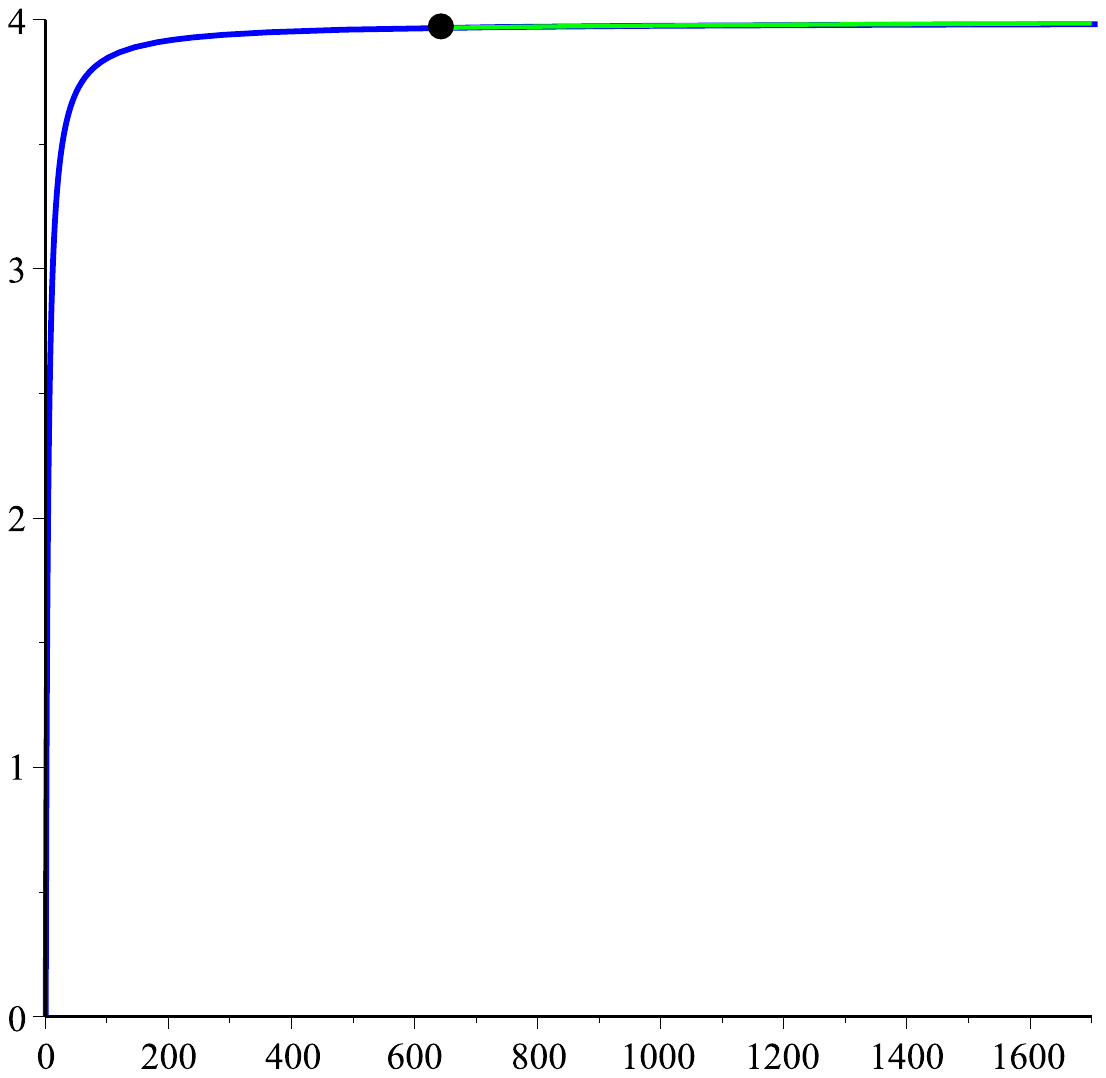}}}
\put(4.9,-4.2){\rotatebox{0}{\includegraphics[width=5cm,height=5.8cm]{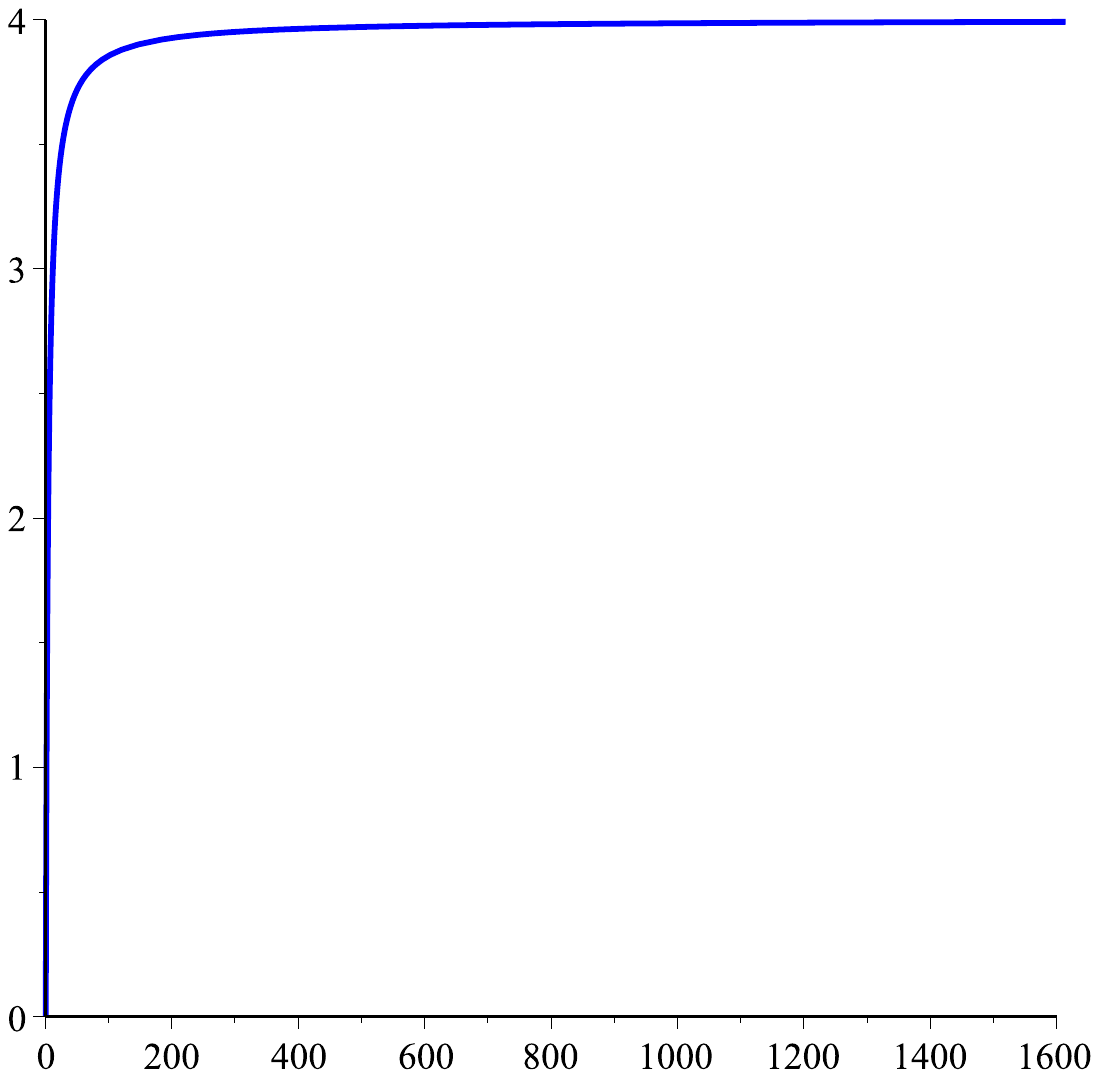}}}
\put(-2.3,5){{\sc $(a)$}}
\put(-4.35,5){{\sc $D$}}
\put(-0.3,4.6){{\sc ${\color{green}\Gamma_{\sLP}}$}}
\put(-0.3,3.1){{\sc ${\color{blue}\Gamma_b}$}}
\put(-3,2.65){{\sc ${\color{red}\Gamma_u}$}}
\put(-3,1.47){{\sc ${\color{magenta}\Gamma_{\sH}}$}}
\put(-0.3,1.3){{\sc $S_{in}$}}
\put(2.6,5){{\sc $(b)$}}
\put(0.65,5){{\sc $D$}}
\put(4.6,4.3){{\sc ${\color{green}\Gamma_{\sLP}}$}}
\put(4.6,3.5){{\sc ${\color{blue}\Gamma_b}$}}
\put(2.2,2.6){{\sc ${\color{red}\Gamma_u}$}}
\put(3.3,2.15){{\sc ${\color{magenta}\Gamma_{\sH}}$}}
\put(3.35,2.1){\tiny \vector(0,-1){0.8}}	
\put(4.6,1.3){{\sc $S_{in}$}}
\put(7.5,5){{\sc $(c)$}}
\put(5.45,5){{\sc $D$}}
\put(9.6,4.8){{\sc ${\color{green}\Gamma_{\sLP}}$}}
\put(9.6,3.9){{\sc ${\color{blue}\Gamma_b}$}}
\put(7.4,2.2){{\sc ${\color{magenta}\Gamma_{\sH}}$}}
\put(7.45,2.15){\tiny \vector(0,-1){0.8}}	
\put(9.6,1.3){{\sc $S_{in}$}}
\put(-2.4,0.7){{\sc $(d)$}}
\put(-4.35,0.55){{\sc $D$}}
\put(-0.3,0.5){{\sc ${\color{green}\Gamma_{\sLP}}$}}
\put(-0.3,0.2){{\sc ${\color{blue}\Gamma_b}$}}
\put(-0.3,-2.9){{\sc $S_{in}$}}
\put(2.6,0.7){{\sc $(e)$}}
\put(0.55,0.55){{\sc $D$}}
\put(4.65,0.5){{\sc ${\color{green}\Gamma_{\sLP}}$}}
\put(1.3,0.62){{\sc ${\color{blue}\Gamma_b}$}}
\put(4.7,-2.9){{\sc $S_{in}$}}
\put(7.2,0.7){{\sc $(f)$}}
\put(5.5,0.55){{\sc $D$}}
\put(9.45,0.5){{\sc ${\color{blue}\Gamma_b}$}}
\put(9.45,-2.9){{\sc $S_{in}$}}
\end{picture}
\end{center}
\vspace{2.5cm}
\caption{Operating diagram in the case considered in Section \ref{SubSec-DOCasFig6} when (a) $a=0.5$ and $b=2$; (b) $a=0.1$ and $b=2$; (c) $a=0.5$ and $b=1.05$; (d) $a=0.01$ and $b=0.2$; (e) $a=0.01$ and $b=0.01$; (f) $a=0$ and $b=0$.}\label{FigDO-LP-H-a05b2}
\end{figure}

\begin{figure}[!h]
\setlength{\unitlength}{1.0cm}
\begin{center}
\begin{picture}(4,5.5)(0,0)
\put(-5.8,0){\rotatebox{0}{\includegraphics[width=4cm,height=6.5cm]{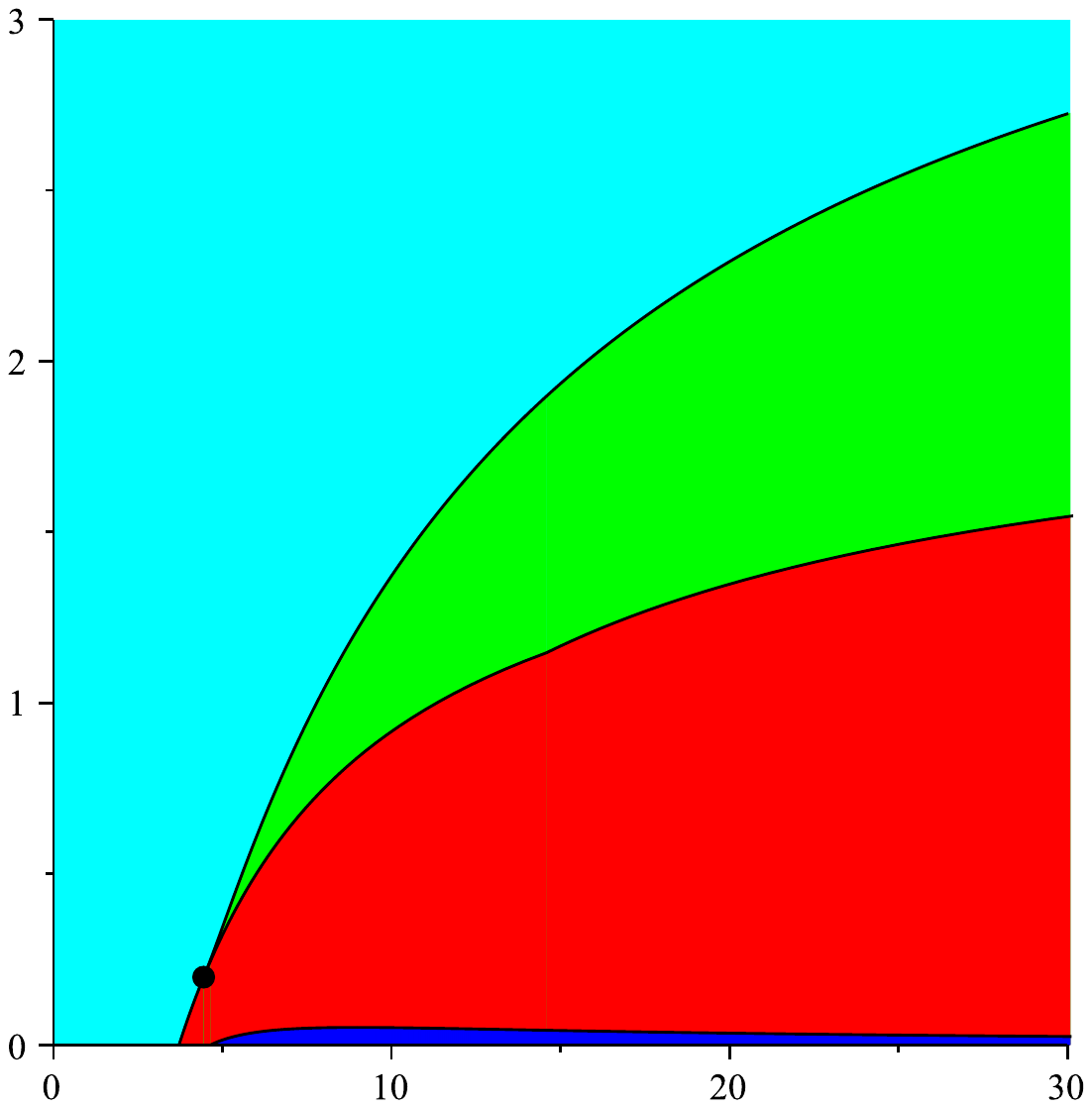}}}
\put(-1.9,0){\rotatebox{0}{\includegraphics[width=4cm,height=6.5cm]{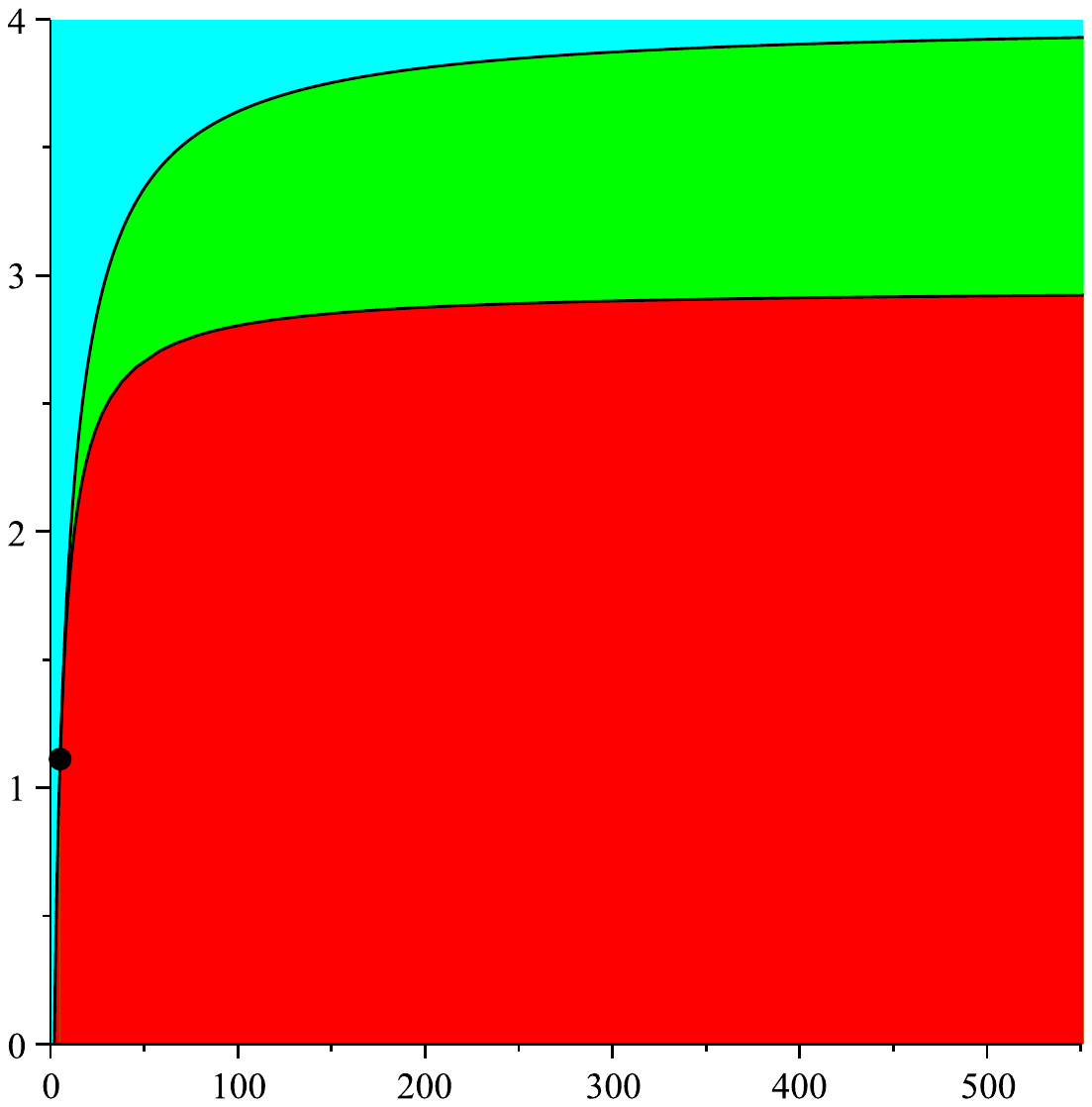}}}
\put(1.9,0){\rotatebox{0}{\includegraphics[width=4cm,height=6.5cm]{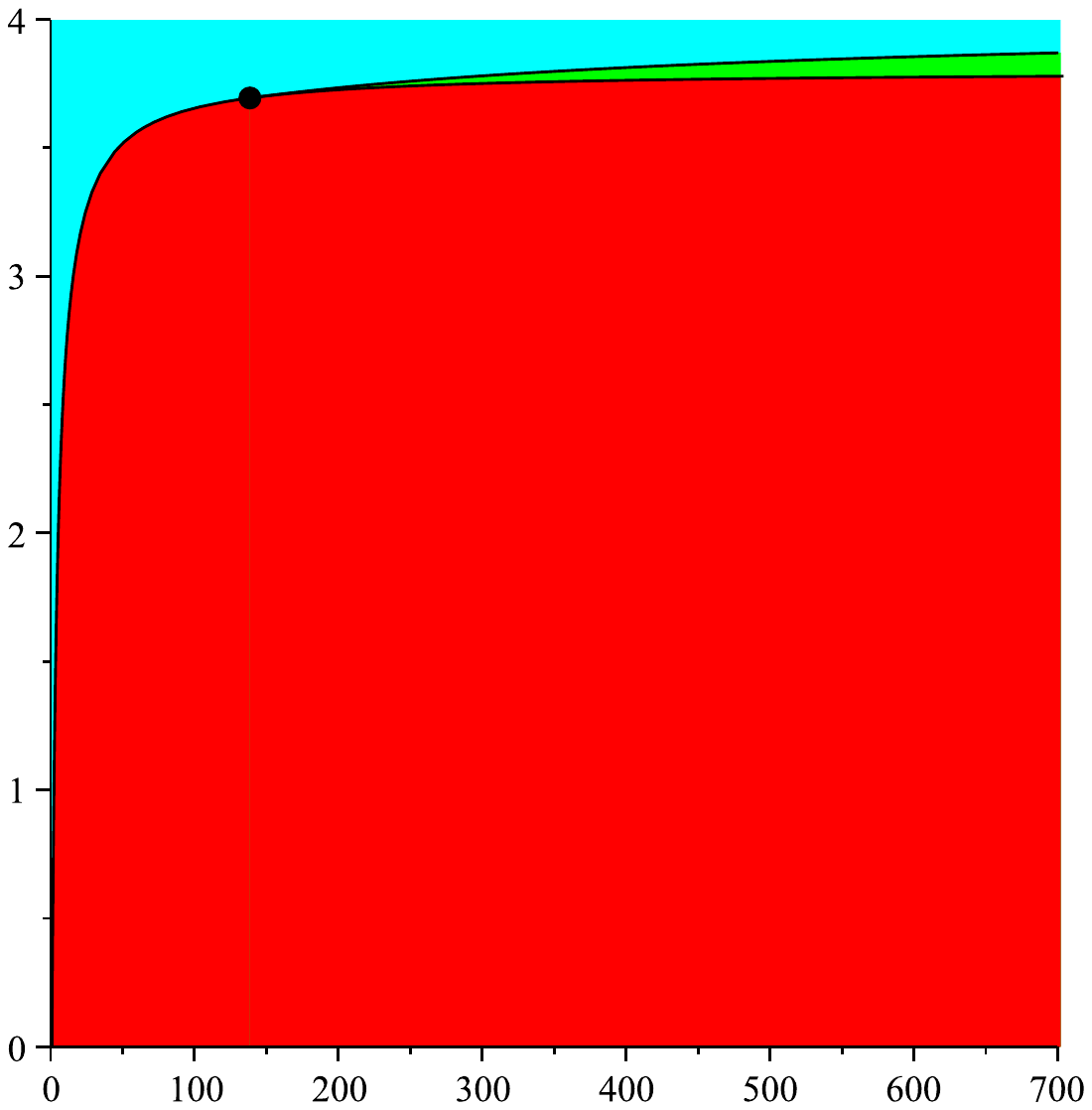}}}
\put(5.8,0){\rotatebox{0}{\includegraphics[width=4cm,height=6.5cm]{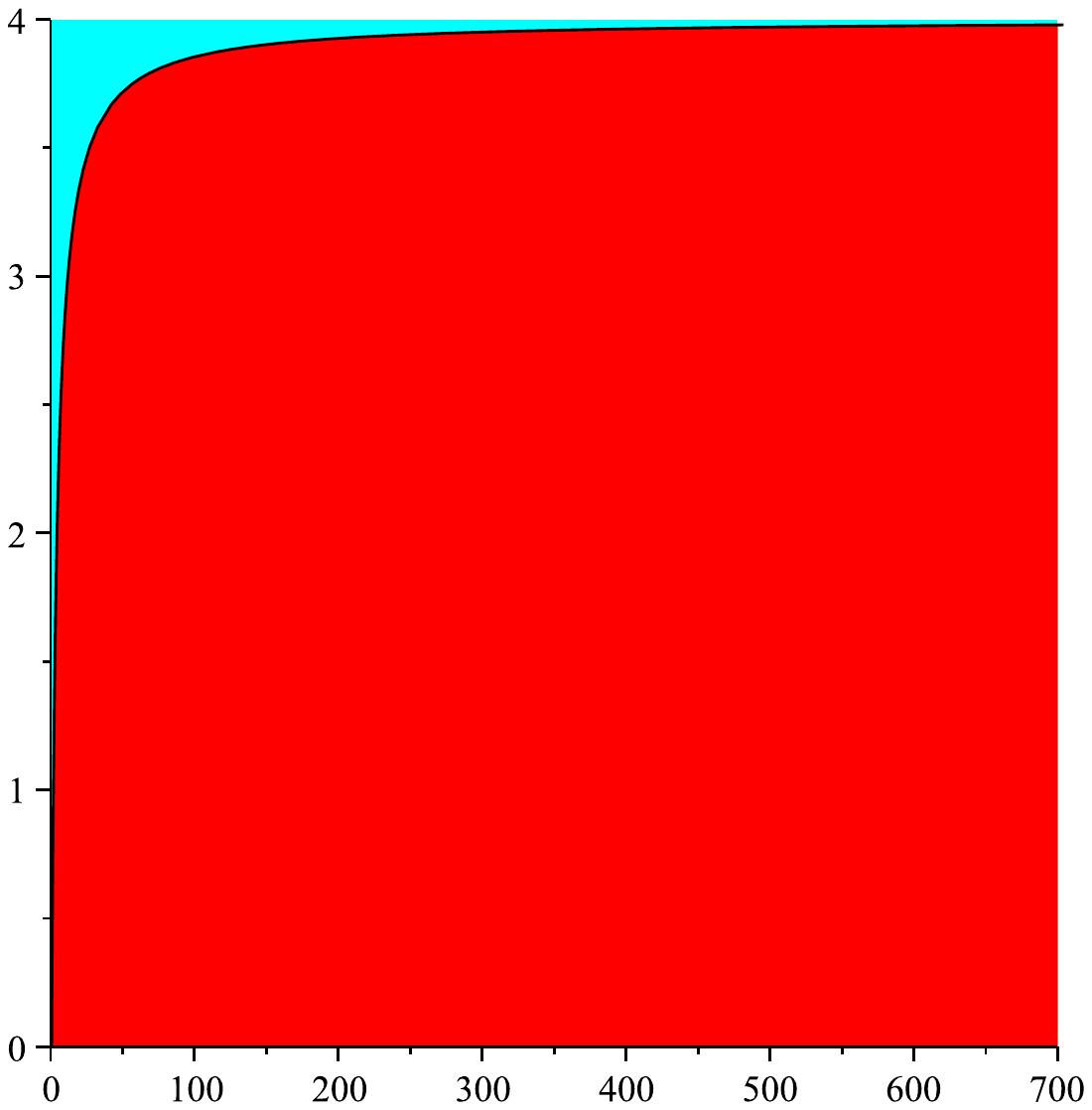}}}
\put(-4,5.5){{\sc $(a)$}}
\put(-5.35,5.4){{\sc $D$}}
\put(-2.1,4.9){{\sc $\Gamma_{\sLP}$}}
\put(-2.1,3.3){{\sc $\Gamma_b$}}
\put(-2.7,1.35){{\sc ${\color{white}\Gamma_{\sH}}$}}
\put(-4.7,4){{\sc $\mathcal{I}_0$}}
\put(-3,3.8){{\sc $\mathcal{I}_2$}}
\put(-4.2,1.7){{\sc ${\color{white}\mathcal{I}_3}$}}
\put(-4.15,1.65){\tiny {\color{white}\vector(0,-1){0.3}}}	
\put(-3.5,2.1){{\sc ${\color{white}\mathcal{I}_1}$}}
\put(-2.1,1.3){{\sc $S_{in}$}}
\put(0,5.5){{\sc $(b)$}}
\put(-1.45,5.4){{\sc $D$}}
\put(1.85,5.15){{\sc $\Gamma_{\sLP}$}}
\put(1.85,4.2){{\sc $\Gamma_b$}}
\put(-1.35,5){{\sc $\mathcal{I}_0$}}
\put(0.2,4.6){{\sc $\mathcal{I}_2$}}
\put(0.1,3){{\sc ${\color{white}\mathcal{I}_1}$}}
\put(1.85,1.3){{\sc $S_{in}$}}
\put(3.8,5.5){{\sc $(c)$}}
\put(2.35,5.4){{\sc $D$}}
\put(5.55,5.2){{\sc $\Gamma_{\sLP}$}}
\put(5.55,4.95){{\sc $\Gamma_b$}}
\put(2.5,5.1){{\sc $\mathcal{I}_0$}}
\put(5,5.15){{\sc $\mathcal{I}_2$}}
\put(4,3){{\sc ${\color{white}\mathcal{I}_1}$}}
\put(5.6,1.3){{\sc $S_{in}$}}
\put(7.6,5.5){{\sc $(d)$}}
\put(6.25,5.4){{\sc $D$}}
\put(9.45,5.25){{\sc $\Gamma_b$}}
\put(6.3,5.15){{\sc $\mathcal{I}_0$}}
\put(7.5,3){{\sc ${\color{white}\mathcal{I}_1}$}}
\put(9.5,1.3){{\sc $S_{in}$}}
\end{picture}
\end{center}
\vspace{-1.8cm}
\caption{Operating diagram in the case considered in Section \ref{SubSec-DOCasFig6} when (a) $a=0.5$ and $b=2$; (b) $a=0.5$ and $b=1.05$; (c) $a=0.01$ and $b=0.01$; (d) $a=0$ and $b=0$.}\label{FigDO-EffectFloc}
\end{figure}

In the following, we consider the same parameter values as in Section \ref{SubSec-DOCasFig6} (or in \cite[Fig. 6]{FekihSIADS2019}) except for the parameters $a$ and $b$ which are variables to see the effects of the attachment and detachment rates on the asymptotic behavior of the process.
Fig. \ref{FigDO-LP-H-a05b2}(a-c) illustrates the reduction in the size of the coexistence region $\mathcal{I}_3$ by decreasing the rates of attachment $a$ and/or detachment $b$.
In Fig. \ref{FigDO-LP-H-a05b2}(d-f), region $\mathcal{I}_3$ has disappeared and region $\mathcal{I}_2$ is reduced to disappearance in the limiting case $a=b=0$ where we obtain the operating diagram of the classic chemostat model.
Fig. \ref{FigDO-EffectFloc} illustrates the operating diagrams with the various colors of regions by decreasing the rates of attachment and detachment where there is a reduction in the size of regions $\mathcal{I}_2$ and $\mathcal{I}_3$ until their disappearance.
\section{Conclusion}                                  \label{Sec-Conc}
In this work, we have extended our mathematical study in \cite{FekihSIADS2019} by considering distinct yields in the flocculation model (\ref{ModFlocGen}) involving the attachment and detachment dynamics of isolated and aggregated bacteria in the presence of a single resource in a chemostat.
Considering distinct removal rates and without ignoring the yield coefficients, we have provided a complete analysis of the existence and local asymptotic stability of all steady states for general monotonic growth rates.
Using the necessary and sufficient conditions of existence and stability, we have studied theoretically and numerically the operating diagrams of flocculation model (\ref{ModFlocGen}) according to the operating parameters which are the dilution rate $D$ and the input concentration of the substrate $S_{in}$.


To have a better understanding of the theoretical study of the operating diagram, we start with a simple case where the positive steady state is unique and stable if it exists (see Section \ref{SubSec-DOCasFig13}).
It is revealed that there can only be two regions: the region $\mathcal{I}_0$ of the washout ($E_0$ is the only steady state) or the region $\mathcal{I}_1$ of coexistence of isolated and attached bacteria around the positive steady state (the only steady states are $E_0$ which is unstable and $E_1$ which is LES).


Next, we have considered a case with the emergence of two positive steady states and the destabilization of one positive steady state via a Hopf bifurcation (see Section \ref{SubSec-DOCasFig12Siam}).
The operating diagram shows the emergence of the green region $\mathcal{I}_2$ corresponding to the bistability between $E_0$ and $E_1^1$, and the blue region $\mathcal{I}_3$
corresponding to the destabilization of the positive steady state $E_1^1$ where there can be coexistence around a stable limit cycle.
In \ref{Sec-AppedixCasArima}, we have considered a similar case for the set of parameter values in \cite{FekihArima2020}.
Indeed, we obtain the operating diagram in Fig. \ref{Fig-CasFig13Siam}(b-c) which is similar to one in Fig. \ref{FigDO-Arima}. However, the region $\mathcal{I}_3$ has not been detected numerically in \cite{FekihArima2020} because of its size where the maximum value of $D$ is around $10^{-5}$.
To detect it with good accuracy, we have changed the default value of ``Digits'' in MAPLE to 20.


Then, we have considered another case with the emergence of two positive steady states and the destabilization of one positive steady state via a Hopf bifurcation (see Section \ref{SubSec-DOCasFig6}).
The operating diagram is divided into five regions where there can be one more behavior (yellow region $\mathcal{I}_4$): the bistability with either coexistence around a stable limit cycle or the washout of the isolated and attached bacteria according to the initial condition.


Using the software MATCONT \cite{MATCONT}, we found numerically the operating diagram obtained theoretically in the case considered in Sections \ref{SubSec-DOCasFig12Siam} and \ref{SubSec-DOCasFig6}.
However, we have also detected new bifurcations with two parameters like those of type Bogdanov-Takens (BT) or Cusp (CP). Moreover, the one bifurcation diagram shows the various types of bifurcations by crossing the different regions in the two-dimensional plane $(S_{in},D)$.
Then, the study of the operating diagram with the two control parameters $(S_{in},D)$ using MATCONT gives a more general vision of the asymptotic behavior of the system compared to the study of the bifurcation diagram as a function of $S_{in}$.


Finally, we analyze the effect of flocculation and deflocculation on the size and shape of various regions in the operating diagram.
Decreasing the rates of attachment and/or detachment, the regions $\mathcal{I}_2$ and $\mathcal{I}_3$ are reduced until their disappearance.
In the limiting case $a=b=0$, we obtain the operating diagram of the classic chemostat model where the CEP asserts that generically at most one species can survive the competition. Thus, the flocculation process promotes the coexistence of isolated and attached bacteria of a microbial species around a limit cycle or positive steady state.
This flocculation mechanism also favors bistability, where the asymptotic behavior of the solutions depends on the initial condition.


The behavior of the process in the various regions of the operating diagram of the model with $n$ species including the mechanism of flocculation is a question of major interest and importance from the biological and ecological point of view. This question deserves further attention and will be the object of future work.
\appendix
\section{Case of Section \ref{SubSec-DOCasFig13}: positivity of the stability condition $c_4$}                      \label{Sec-AppedixCasFig13}
\begin{figure}[!h]
\setlength{\unitlength}{1.0cm}
\begin{center}
\begin{picture}(8.1,6)(0,0)
\put(-3.5,0){\rotatebox{0}{\includegraphics[width=7cm,height=7cm]{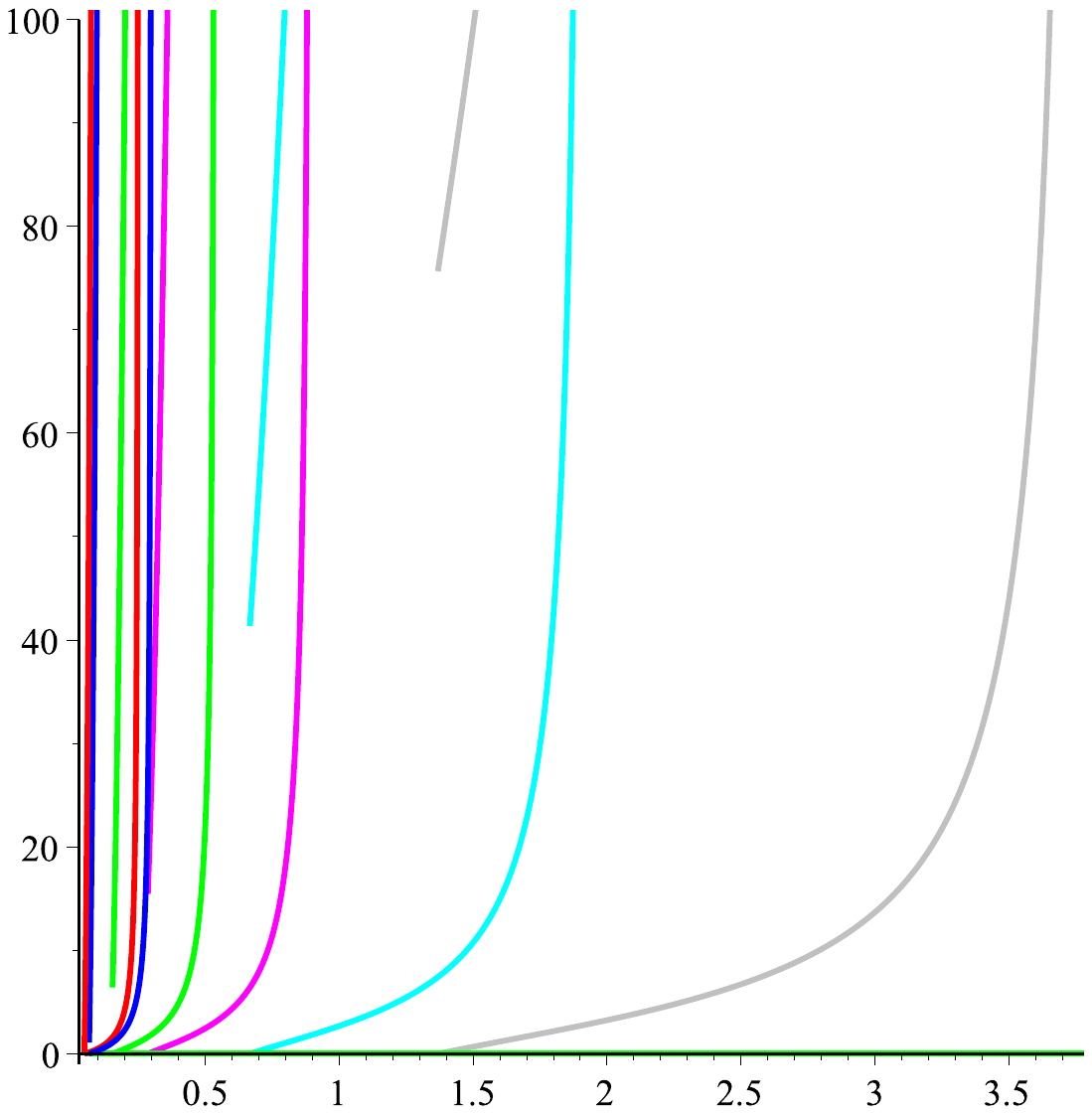}}}
\put(3.5,0){\rotatebox{0}{\includegraphics[width=7cm,height=7cm]{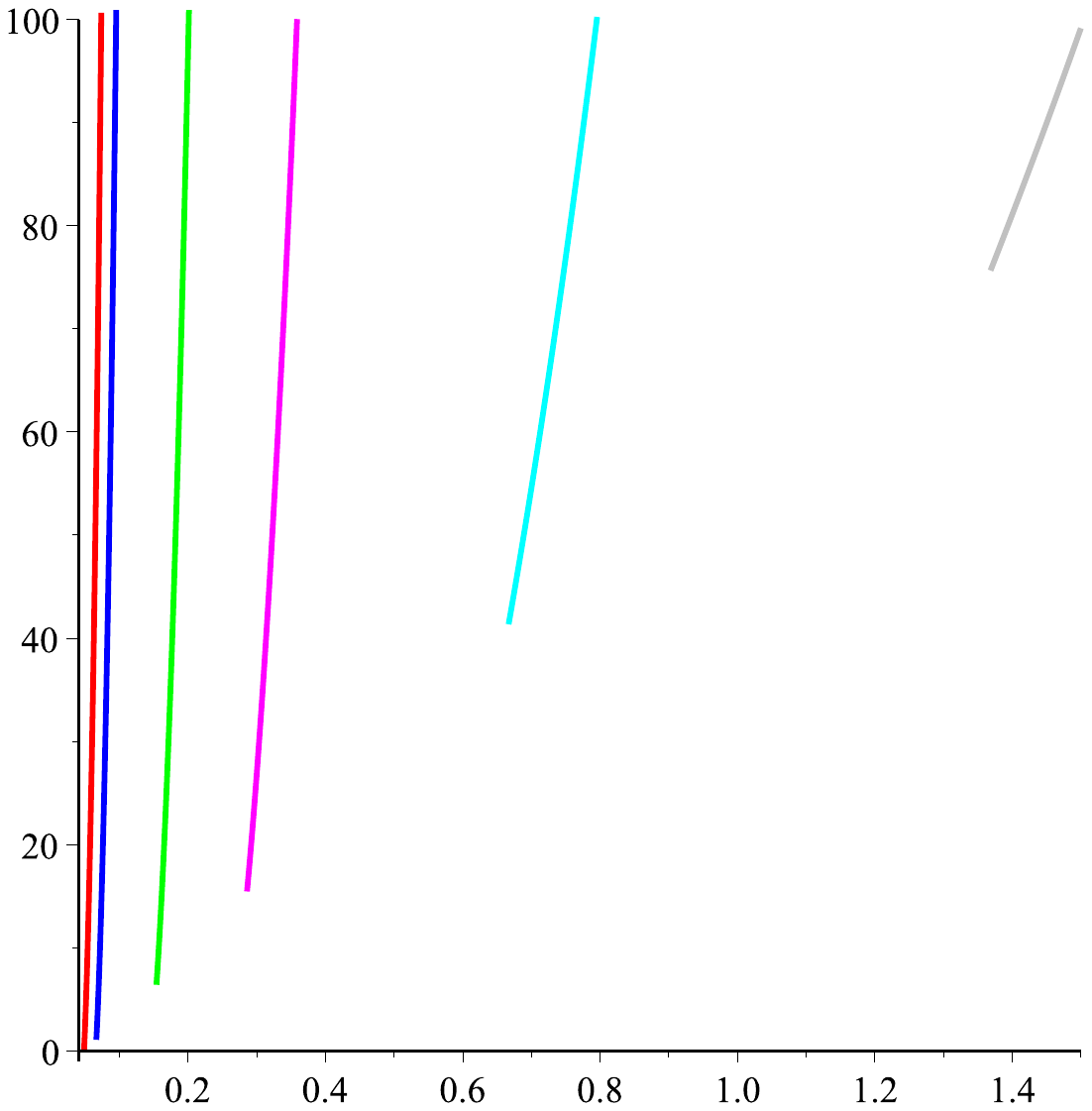}}}
\put(0.4,5.9){{\sc $(a)$}}
\put(-0.35,5.65){{\sc ${\color{gray}c_4}$}}
\put(2.8,5.6){{\sc ${\color{gray}H}$}}
\put(3,1.35){{\sc $S$}}
\put(7.5,5.9){{\sc $(b)$}}
\put(10,5.55){{\sc ${\color{gray}c_4}$}}
\put(10,1.35){{\sc $S$}}
\end{picture}
\end{center}
\vspace{-1.8cm}
\caption{Case of Section \ref{SubSec-DOCasFig13}: (a) curves of the function $H(S)$ and the corresponding curves of $c_4(S)$ in their existence domain when $D\in \{0.01, 0.1, 0.5, 1, 2, 3 \}$ corresponding to red, blue, green, magenta, cyan, and grey colors, respectively. (b) Only the curves of the function $c_4(S)$.}\label{Fig-CasFig13SiamHC4}
\vspace{-0.4cm}
\end{figure}

In the following, we show that the stability condition $c_4>0$ holds for the positive steady state $E_1$ in the case considered in Section \ref{SubSec-DOCasFig13} so that the curve $\Gamma_{\sH}$ corresponding to $c_4=0$ does not exist in the operating diagram of Fig. \ref{Fig-CasFig13Siam}(a) for model (\ref{ModFlocGen}). Fig. \ref{Fig-CasFig13SiamHC4} illustrates the positivity of the function $c_4(S)$ for several values of $D$ and the corresponding curves of $H(S)$.
\section{Case of Section \ref{SubSec-DOCasFig12Siam}: sign of $c_4$}                      \label{Sec-AppedixCasFig12}
\begin{figure}[!h]
\setlength{\unitlength}{1.0cm}
\begin{center}
\begin{picture}(7.5,6.5)(0,0)
\put(-3.7,0){\rotatebox{0}{\includegraphics[width=7.5cm,height=7.5cm]{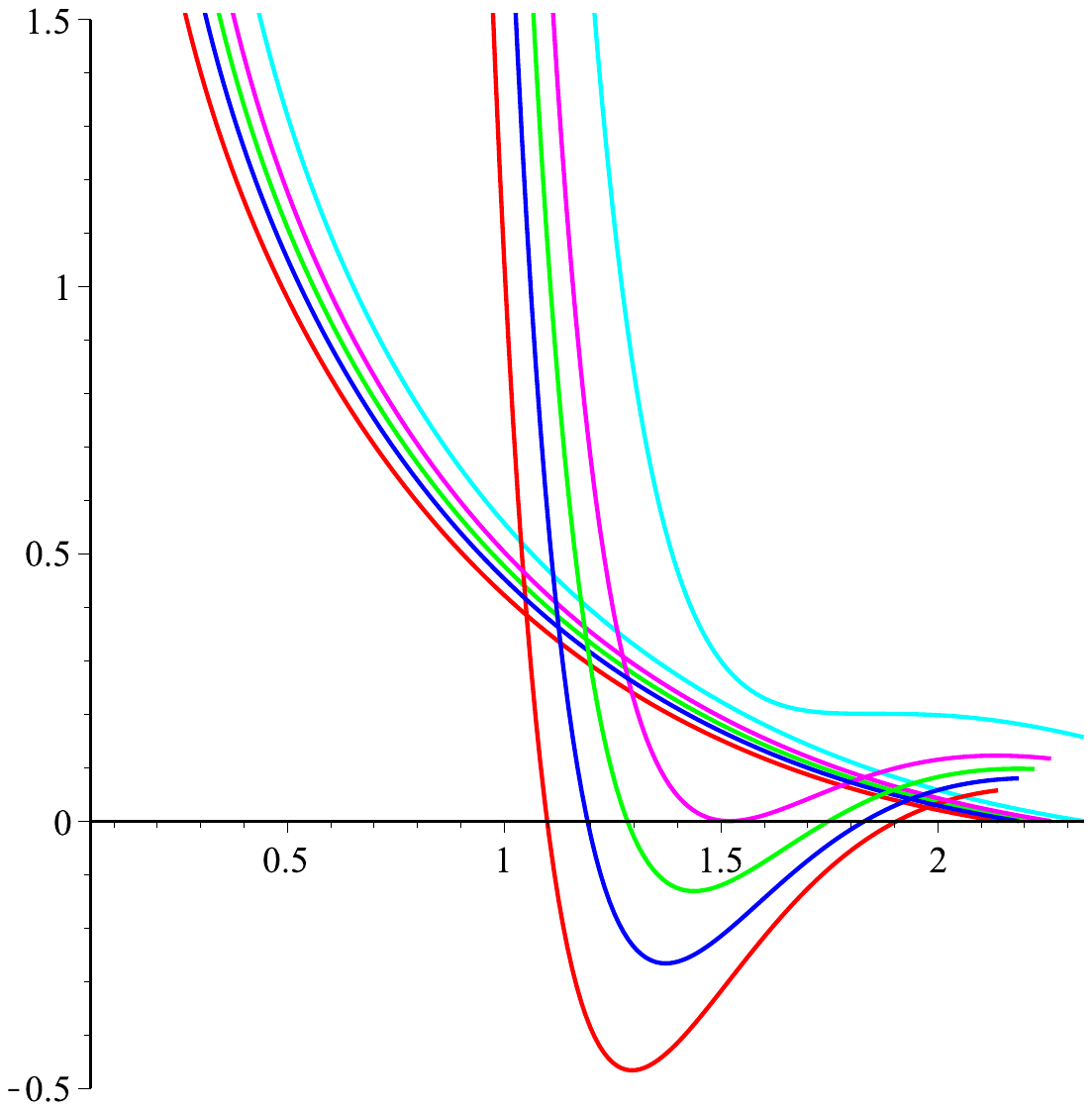}}}
\put(3.7,0){\rotatebox{0}{\includegraphics[width=7.5cm,height=7.5cm]{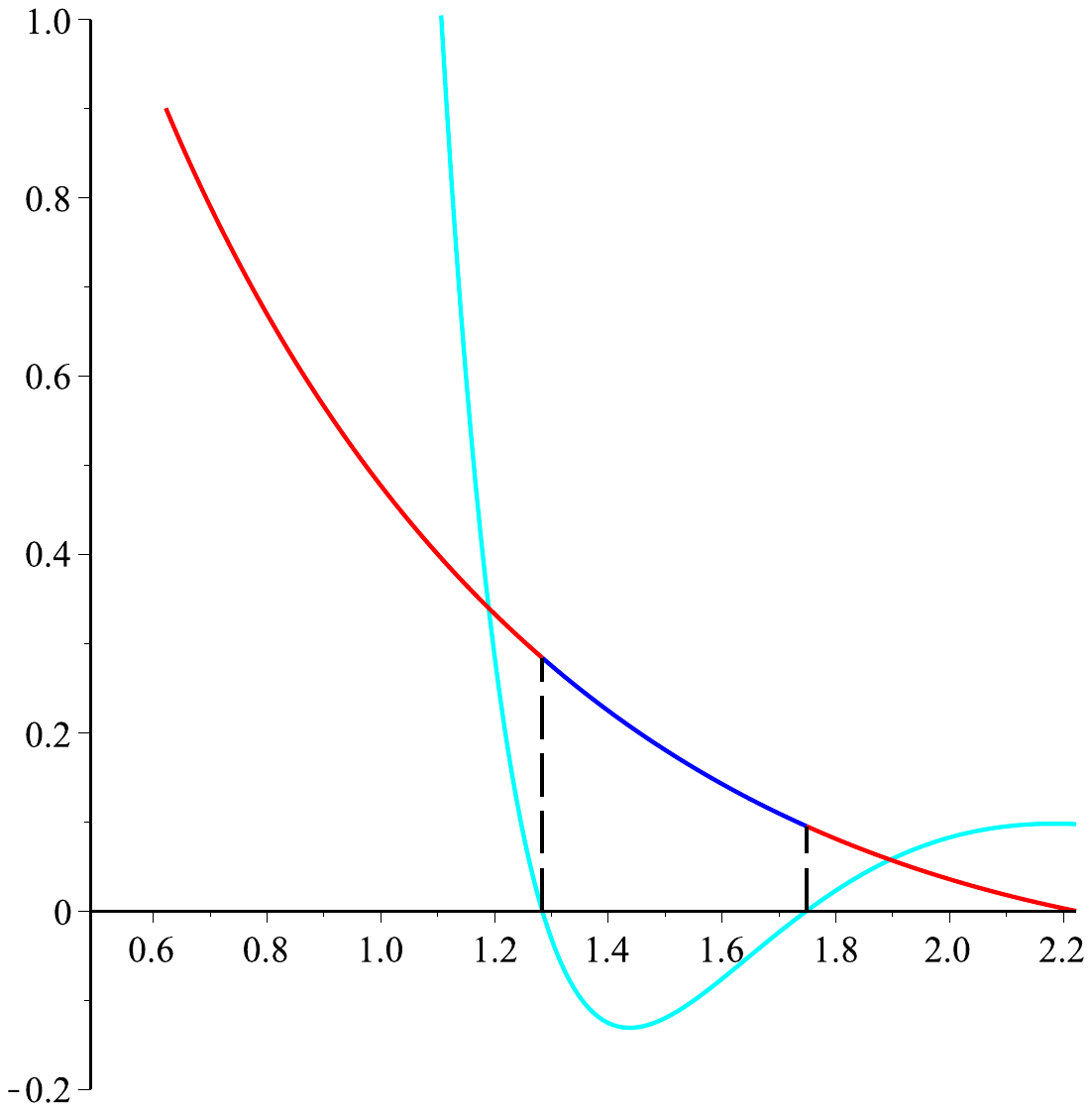}}}
\put(0.2,6.4){{\sc $(a)$}}
\put(-1.5,6){{\sc ${\color{cyan}c_4}$}}
\put(0.4,6){{\sc ${\color{cyan}H}$}}
\put(3.3,2.45){{\sc $S$}}
\put(7.5,6.4){{\sc $(b)$}}
\put(6.9,5.5){{\sc ${\color{cyan}c_4}$}}
\put(6,4.65){{\sc ${\color{red}H}$}}
\put(8,2.95){{\sc ${\color{blue}H}$}}
\put(9.95,2.25){{\sc ${\color{red}H}$}}
\put(10.65,2.1){{\sc $S$}}
\end{picture}
\end{center}
\vspace{-1.8cm}
\caption{Case of Section \ref{SubSec-DOCasFig12Siam}: (a) the functions $H(S)$ and $c_4(S)$ on the existence domain of $E_1$ when $D\in \left\{0.09, 0.12, 0.142, D_{\sH}^{max}, 0.21 \right\}$ $\left(D_{\sH}^{max}\approx 0.165\right)$ corresponding to red, blue, green, magenta, and cyan colors, respectively. (b) The function $H(S)$ in red [resp. in blue] when the function $c_4(S)$ is positive [resp. negative], for $D=D^\ast=0.142$.}\label{Fig-CasFig12Siam-Hc4}
\end{figure}

In the section, we show that the stability condition $c_4>0$ of the positive steady state $E_1$ is not always verified, so that $c_4(S)$ changes sign in the interval $I(D)=\left]\lambda_v(D),\lambda_{\sBP}(D)\right]$ of the existence of $E_1$.
Fig. \ref{Fig-CasFig12Siam-Hc4}(a) shows that the equation $c_4=0$ has two roots $S_{\sH}^1(D)$ and $S_{\sH}^2(D)$ as defined in Section \ref{SubSec-DOCasFig12Siam} for all $D<D_{\sH}^{max}$.
Let $D$ be fixed at $D=D^\ast=0.142$ (the green curve in \ref{Fig-CasFig12Siam-Hc4}(a)).
Fig. \ref{Fig-CasFig12Siam-Hc4}(b) shows the curve of the function $S\mapsto H(S)$ in red [resp. in blue] when the function $S\mapsto c_4(S)$ is positive [resp. negative].
More precisely, $c_4(S)$ is positive for all $S\in\left]\lambda_v,S_{\sH}^2\right[\cup\left]S_{\sH}^1,\lambda_{\sBP}\right]$ and negative for all $S\in\left]S_{\sH}^2,S_{\sH}^1\right[$ where $\lambda_v\left(D^\ast\right) \approx 0.078$, $\lambda_{\sBP}(D^\ast)=\lambda_b\approx2.222$ and the critical values according to $S$ and corresponding to Hopf bifurcation are given by
$$
 S_{\sH}^2\left(D^\ast\right) \approx 1.284, \quad S_{\sH}^1\left(D^\ast\right) \approx 1.748.
$$
These critical values are equivalent to the following critical values according to $S_{in}$,
$$
S_{in}^{H2}\left(D^\ast\right) \approx 3.674, \quad S_{in}^{H1}\left(D^\ast\right) \approx 2.640, \mbox{ respectively}.
$$
By increasing the value of $S_{in}$ from zero to $\lambda_{\sBP}(D^\ast)$, $E_1$ emerges LES via a Branch Point (BP) with $E_0$ when $S=S_{in}=\lambda_{\sBP}(D^\ast)$. Increasing $S_{in}$ further, $E_1$ destabilizes trough the first Hopf bifurcation at $S_{in}^{H1}$ and remains unstable up to the value of $S_{in}^{H2}$.
Finally, $E_1$ returns LES for all $S_{in}>S_{in}^{H2}$ via a second Hopf bifurcation.
\section{Case of Section \ref{SubSec-DOCasFig6}}         \label{Sec-AppedixCasFig6}
\begin{figure}[!h]
\setlength{\unitlength}{1.0cm}
\begin{center}
\begin{picture}(6.7,5.3)(0,0)
\put(-3.7,0){\rotatebox{0}{\includegraphics[width=7cm,height=6.5cm]{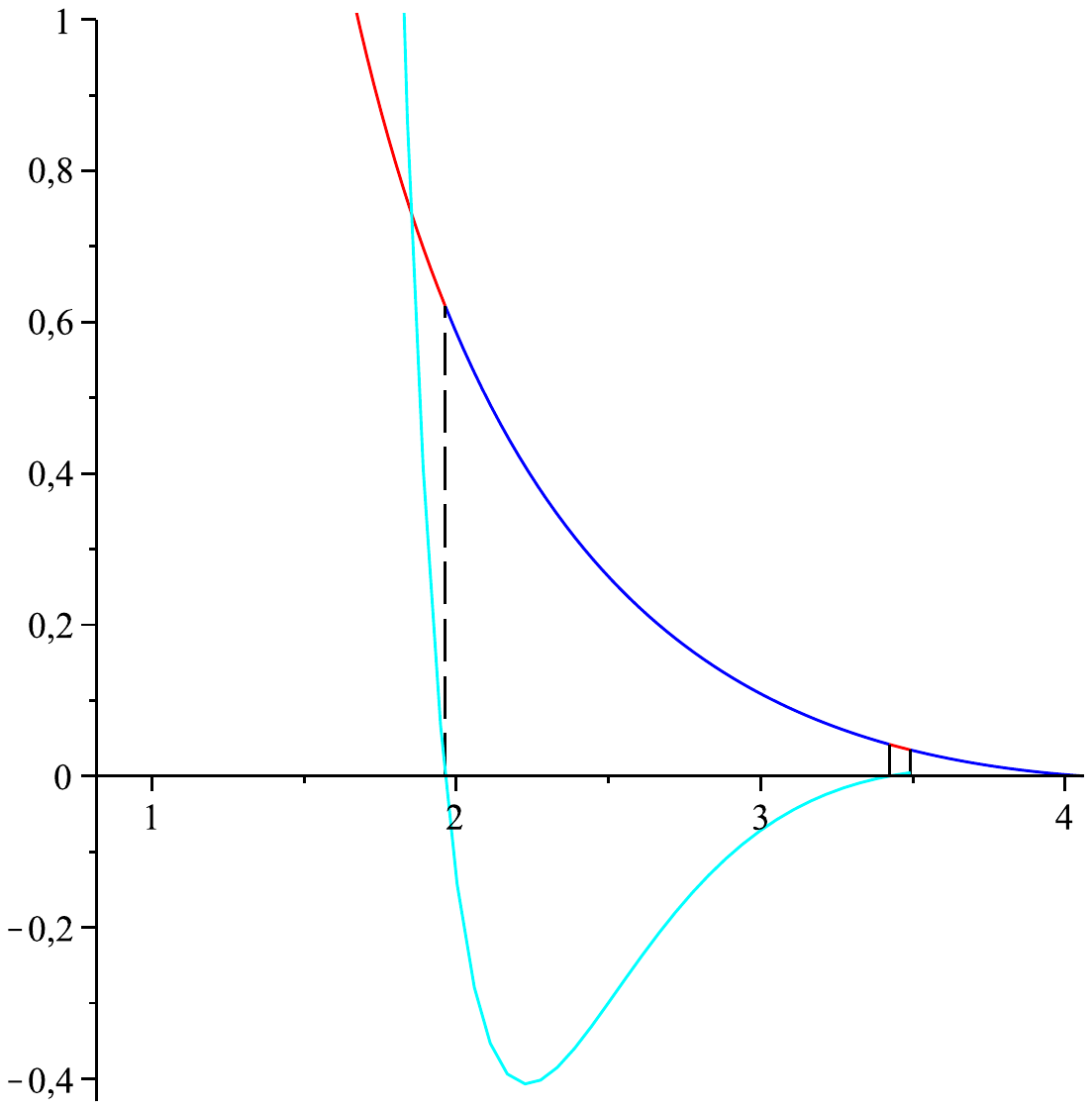}}}
\put(3.7,0){\rotatebox{0}{\includegraphics[width=7cm,height=6.5cm]{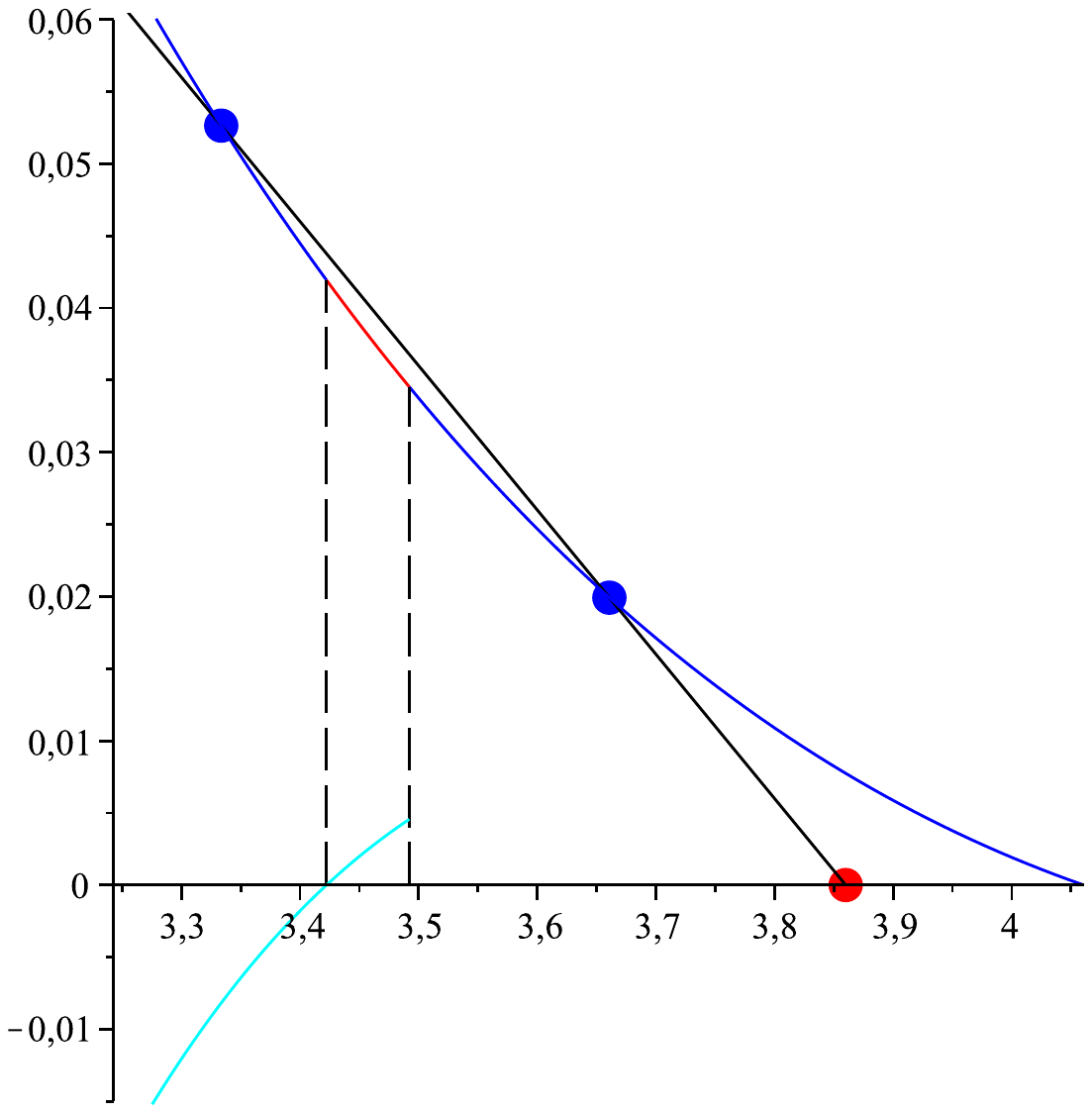}}}
\put(0.3,5.2){{\sc $(a)$}}
\put(-0.85,4.9){{\sc ${\color{cyan}c_4}$}}
\put(0.2,3.05){{\sc $H$}}
\put(-1,1.9){{\sc  $S_{\sH}^2$}}
\put(1.6,1.9){{\sc  $S_{\sH}^1$}}
\put(2,1.9){{\sc  $S_{\sLP}$}}
\put(3,2.2){{\sc  $S$}}
\put(7.5,5.2){{\sc $(b)$}}
\put(5.55,4.7){{\sc ${\color{blue}E_1^1}$}}
\put(7.6,3){{\sc ${\color{blue}E_1^2}$}}
\put(8.9,1.5){{\sc ${\color{red}E_0}$}}
\put(8.9,2.3){{\sc $H$}}
\put(6.6,3.8){{\sc $\delta$}}
\put(5.1,1.3){{\sc ${\color{cyan}c_4}$}}
\put(5.9,1.4){{\sc  $S_{\sH}^1$}}
\put(6.5,1.4){{\sc  $S_{\sLP}$}}
\put(10.45,1.8){{\sc  $S$}}
\end{picture}
\end{center}
\vspace{-1.6cm}
\caption{Case of Section \ref{SubSec-DOCasFig6}: (a) the function $H(S)$ and $c_4(S)$ when $D=D^\ast=0.1$ showing the changes of the sign of $c_4(S)$; (b) a magnification for $3.25<S<\lambda_u=4.061$ where $S_{in}=3.86$.} \label{Fig-Hc4-CasFig6Siam}
\end{figure}

In this appendix, we give numerical evidence of the change of sign of $c_4(S)$ and the appearance of a stable limit cycle as well as its disappearance by a homoclinic bifurcation for fixed $D$.
Fig. \ref{Fig-Hc4-CasFig6Siam} shows the curve of the function $S\mapsto H(S)$ in red [resp. in blue] when the function $S\mapsto c_4(S)$ is positive [resp. negative] and $D$ is fixed at $D^\ast=0.1$.
More precisely, $c_4(S)$ is positive for all $S\in\left]\lambda_v,S_{\sH}^2\right[\cup\left]S_{\sH}^1,S_{\sLP}\right]$ and negative for all $S\in\left]S_{\sH}^2,S_{\sH}^1\right[$ where $\lambda_v\left(D^\ast\right) \approx 0.846$, and the critical values according to $S$ and corresponding to Hopf bifurcation are given by
$$
 S_{\sH}^2\left(D^\ast\right) \approx 1.963, \quad S_{\sH}^1\left(D^\ast\right) \approx 3.422, \quad S_{\sLP}\left(D^\ast\right) \approx 3.492.
$$
These critical values are equivalent to the following critical values according to $S_{in}$,
$$
S_{in}^{H2}\left(D^\ast\right) \approx 8.179, \quad S_{in}^{H1}\left(D^\ast\right) \approx 3.842, \quad \lambda_{\sLP}\left(D^\ast\right)\approx 3.837, \mbox{ respectively}.
$$
By increasing the value of $S_{in}$ from zero to $\lambda_{\sLP}\left(D^\ast\right)$, the two positive steady states $E_1^1$ and $E_1^2$ emerge LES and unstable, respectively, via a Limit Points (LP) bifurcation where $S=S_{\sLP}\left(D^\ast\right)$. Increasing $S_{in}$ further, $E_1^1$ destabilizes trough the first Hopf bifurcation at $S_{in}^{H1}$ and remains unstable up to the value of $S_{in}^{H2}$.
Finally, $E_1^1$ returns LES for all $S_{in}>S_{in}^{H2}$ via a second Hopf bifurcation.
In Fig. \ref{Fig-Hc4-CasFig6Siam}(b), we have chosen the red color for LES steady states and the blue color for unstable steady states.

\begin{figure}[tbhp]
\setlength{\unitlength}{1.0cm}
\begin{center}
\begin{picture}(5.5,6)(0,0)
\put(-5,0){\rotatebox{0}{\includegraphics[width=7.5cm,height=14cm]{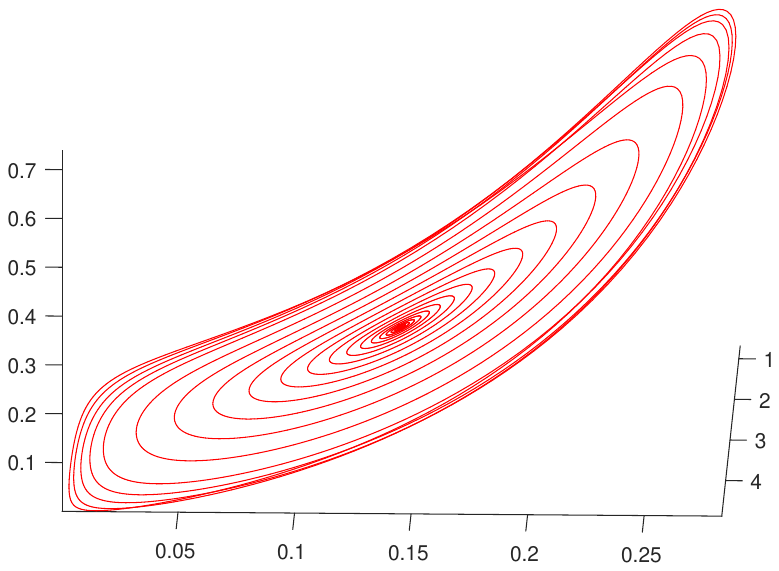}}}
\put(0,0){\rotatebox{0}{\includegraphics[width=7.5cm,height=14cm]{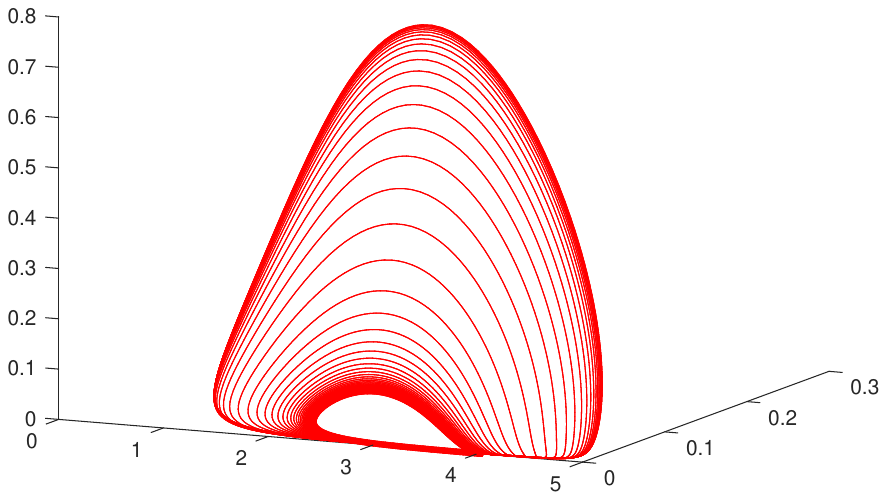}}}
\put(5.7,0){\rotatebox{0}{\includegraphics[width=7.5cm,height=14cm]{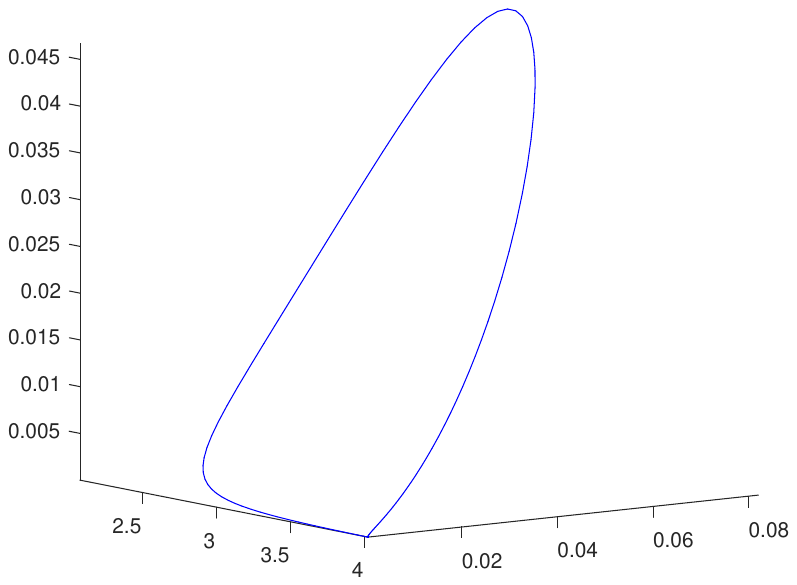}}}
\put(-2,5){{\sc  $(a)$}}
\put(-2.3,2.25){{\sc H}}
\put(-2.4,2.15){{\sc ${\color{black}\bullet}$}}
\put(-4.3,3.8){{\sc $v$}}
\put(-2.7,0.2){{\sc $u$}}
\put(-0.1,1.3){{\sc $S$}}
\put(3,5){{\sc  $(b)$}}
\put(0.8,4.4){{\sc $v$}}
\put(4.75,0.6){{\sc $u$}}
\put(2.1,0.4){{\sc $S$}}
\put(8.7,5.2){{\sc  $(c)$}}
\put(6.4,4.8){{\sc $v$}}
\put(9.5,0.3){{\sc $u$}}
\put(6.85,0.4){{\sc $S$}}
\end{picture}
\end{center}
\vspace{-0.8cm}
\caption{Case of Section \ref{SubSec-DOCasFig6}: the three-dimensional space $(S,u,v)$ in MATCONT when $D=0.1$: (a) a family of limit cycles, starting from a Hopf point (H) at $\sigma_6=8.179$ and decreasing $S_{in}$ until the maximum radius. (b) A family of limit cycles, starting from the maximum radius and approaching a homoclinic orbit by decreasing $S_{in}$. (c) A homoclinic orbit for $S_{in}=\sigma_4\approx 4.03468$.} \label{Fig3DHom2}
\end{figure}

\begin{figure}[!h]
\setlength{\unitlength}{1.0cm}
\begin{center}
\begin{picture}(8.5,5.3)(0,0)
\put(-4,0){\rotatebox{0}{\includegraphics[width=12cm,height=12cm]{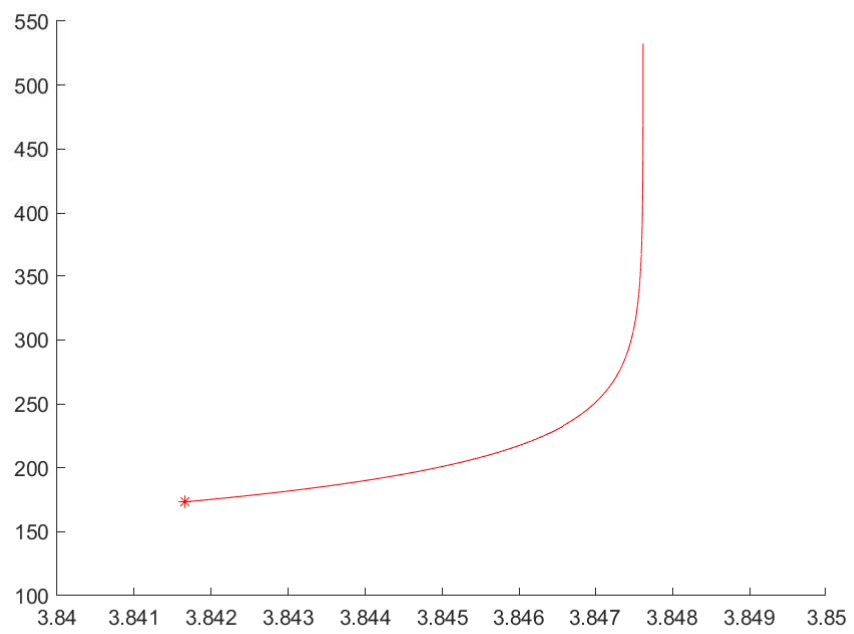}}}
\put(4,0){\rotatebox{0}{\includegraphics[width=10cm,height=10cm]{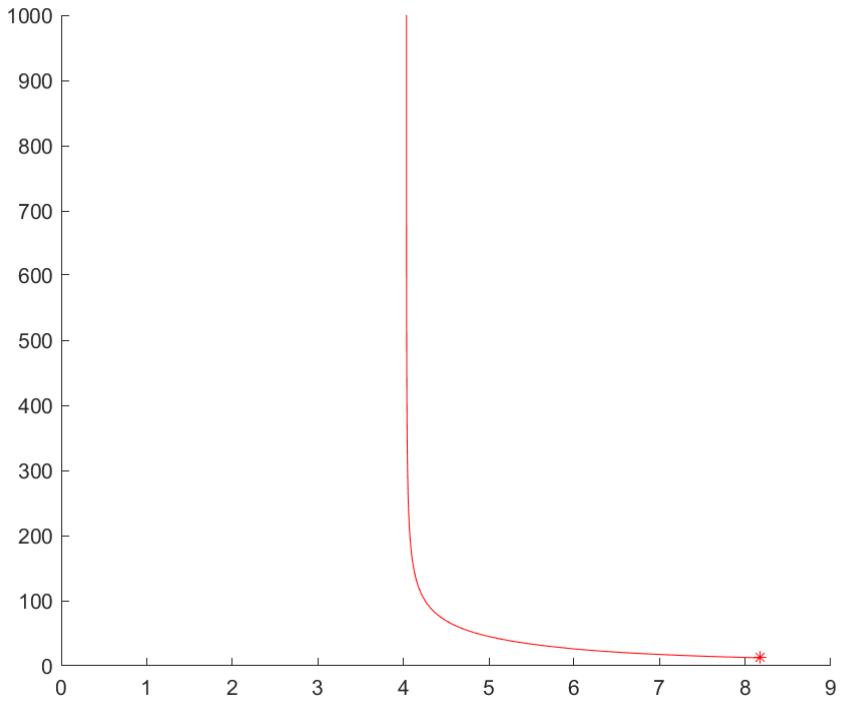}}}
\put(0.5,4.9){{\sc $(a)$}}
\put(-2.65,4.7){{\sc $T$}}
\put(4.05,0.55){{\sc $S_{in}$}}
\put(8,4.9){{\sc $(b)$}}
\put(5.2,4.5){{\sc $T$}}
\put(11.1,0.55){{\sc $S_{in}$}}
\end{picture}
\end{center}
\vspace{-1cm}
\caption{Case of Section \ref{SubSec-DOCasFig6}: A plot in MATCONT of time period $T$ of the limit cycle solutions of model (\ref{ModFlocGen}) for $D=0.1$ starting from the Hopf bifurcation at (a) $\sigma_2$ [(b) resp. $\sigma_6$]; homoclinic bifurcation at $S_{in}=\sigma_3\approx 3.8477$ [resp. $S_{in}=\sigma_4\approx 4.03468$].}\label{FigT-HomCasFig6Siam}
\end{figure}

Fig. \ref{Fig3DHom2} shows the limit cycles in the three-dimensional phase plot $(S,u,v)$ for various values of $S_{in}$ between $\sigma_4$ and $\sigma_6$ (defined in Table \ref{TableSigmai}) until their disappear by homoclinic bifurcation at $\sigma_4$.
Decreasing $S_{in}$ from the Hopf bifurcation at $\sigma_6$, the radius of the limit cycle increases (see Fig. \ref{Fig3DHom2}(a)).
Decreasing again $S_{in}$, the radius of the limit cycle decreases until his disappearance by approaching a homoclinic orbit when $S_{in}=\sigma_4\approx 4.03468$ (see Fig. \ref{Fig3DHom2}(b-c)).
Fig. \ref{FigT-HomCasFig6Siam} reveals the homoclinic bifurcations at $S_{in}=\sigma_3$ and $S_{in}=\sigma_4$ where the time period $T$ of the limit cycle solutions of model (\ref{ModFlocGen}) tends to $+\infty$ as $S_{in}$ tends to these critical values.
\section{Case of parameter set in Table \ref{TabParamVal} (line 5)}    \label{Sec-AppedixCasArima}
\begin{figure}[!h]
\setlength{\unitlength}{1.0cm}
\begin{center}
\begin{picture}(8.4,6.5)(0,0)
\put(-4.1,0){\rotatebox{0}{\includegraphics[width=8cm,height=7.5cm]{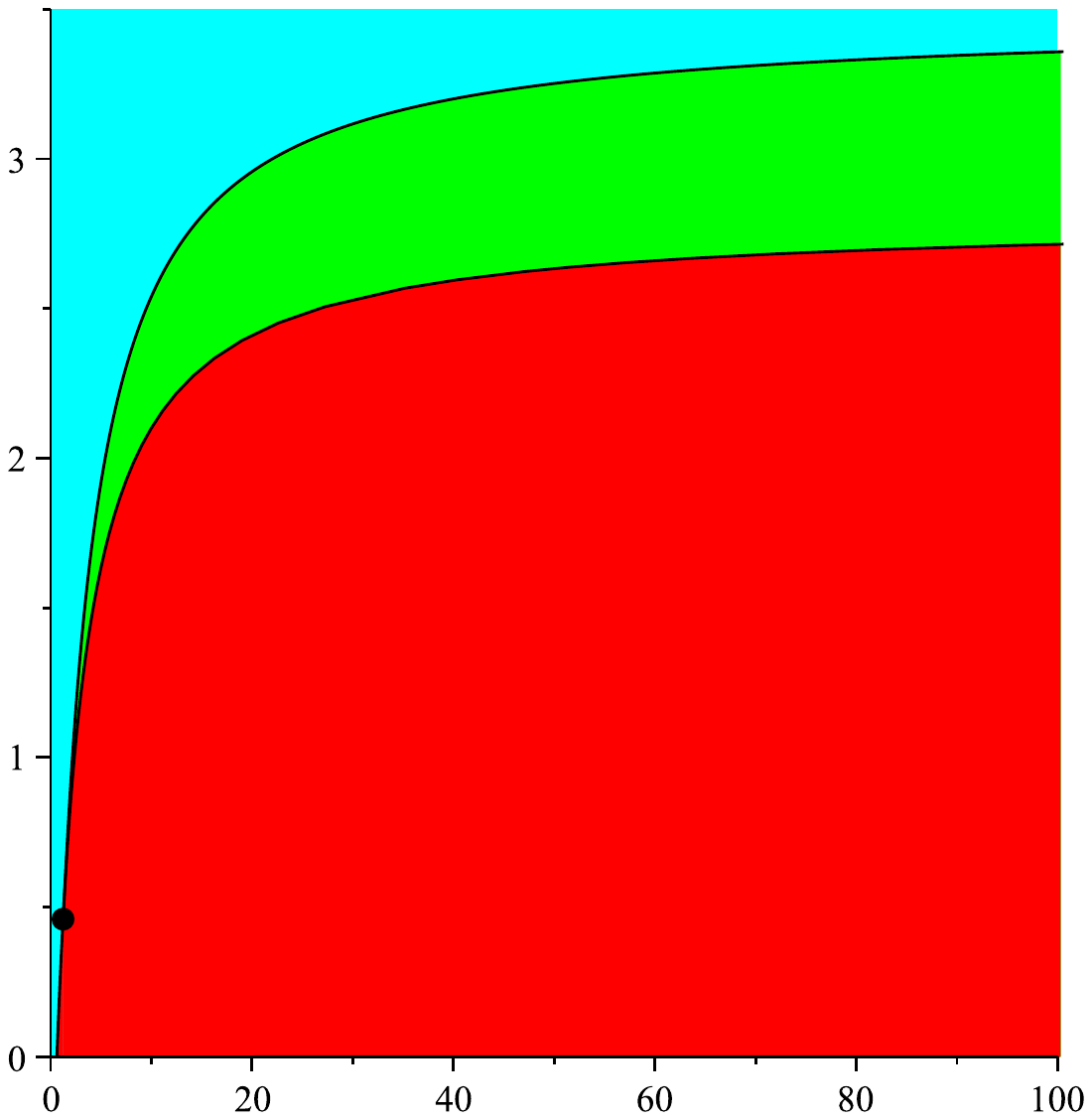}}}
\put(4.1,0){\rotatebox{0}{\includegraphics[width=8cm,height=7.5cm]{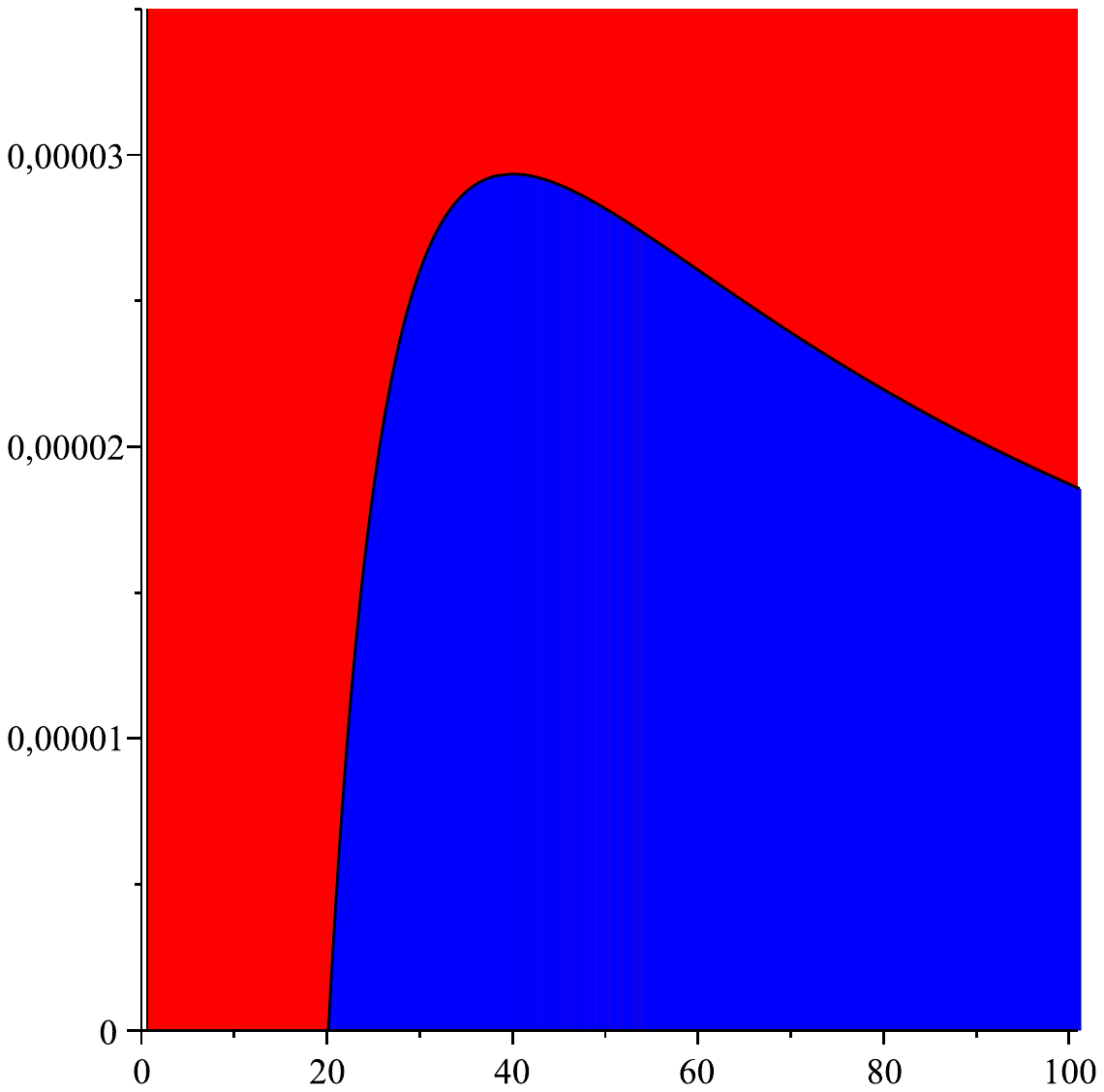}}}
\put(0,6.4){{\sc $(a)$}}
\put(-3.2,6.35){{\sc $D$}}
\put(3.2,6){{\sc $\Gamma_{\sLP}$}}
\put(3.2,5.2){{\sc $\Gamma_u$}}
\put(-2.6,5.8){{\sc $\mathcal{I}_0$}}
\put(-0.7,5.4){{\sc $\mathcal{I}_2$}}
\put(-0.2,3.3){{\sc ${\color{white}\mathcal{I}_1}$}}
\put(3.2,1.6){{\sc $S_{in}$}}
\put(8,6.4){{\sc $(b)$}}
\put(5.35,6.35){{\sc $D$}}
\put(11.45,4){{\sc $\Gamma_{\sH}$}}
\put(7,5.8){{\sc ${\color{white}\mathcal{I}_1}$}}
\put(8.6,3.3){{\sc ${\color{white}\mathcal{I}_3}$}}
\put(6.6,4){{\sc ${\color{white}\Gamma_{\sH}^1}$}}
\put(9.2,5){{\sc ${\color{white}\Gamma_{\sH}^2}$}}
\put(11.5,1.5){{\sc $S_{in}$}}
\end{picture}
\end{center}
\vspace{-1.8cm}
\caption{MAPLE: (a) operating diagram of (\ref{ModFlocGen}) in case of  Table \ref{TabParamVal} (line 5). (b) Magnification on the region $\mathcal{I}_3$ and the curve $\Gamma_{\sH}=\Gamma_{\sH}^1 \cup \Gamma_{\sH}^2$ when $D\in [0,0.000035]$.}\label{FigDO-Arima}
\end{figure}

The main purpose of this appendix is to show that the region of destabilization of the positive steady state is omitted in the construction of the operating diagram in \cite{FekihArima2020}.
With the same set of parameters in \cite{FekihArima2020}, see Table \ref{TabParamVal} (line 5), we find the operating diagram in Fig. \ref{FigDO-Arima} which is similar to that in Fig. \ref{Fig-CasFig13Siam}(b-c) but where $\lambda_{\sBP}(D)=\lambda_u(D)$. Thus, the existence and the local stability of all steady states of model (\ref{ModFlocGen}) in the four regions $\mathcal{I}_i$, $i=0,\ldots,3$ of the operating diagram shown in Fig. \ref{FigDO-Arima} can be obtained from Table \ref{Tab-DO-LP-Fig12Siam}.
Note that the numbering of the $\mathcal{I}_1$ and $\mathcal{I}_2$ regions is reversed in \cite{FekihArima2020}. Similarly, for the regions $\mathcal{I}_3$ and $\mathcal{I}_4$.
Next, we will show that the region $\mathcal{I}_3$ corresponds to the emergence of the stable limit cycle via Hopf bifurcations.

Since the order of $D$ can reach $10^{-10}$ to plot the curve $\Gamma_{\sH}$, we modified the value of ``Digits'' in MAPLE to 20 instead of the default value 10 to avoid the introduction of round-off error.
This allows in particular to have precision in the tracing of the $\Gamma_{\sH}$ curve and that of the function $c_4(S)$ for fairly small $D$.

\begin{figure}[!h]
\setlength{\unitlength}{1.0cm}
\begin{center}
\begin{picture}(7,6.3)(0,0)
\put(-4,0){\rotatebox{0}{\includegraphics[width=7cm,height=7cm]{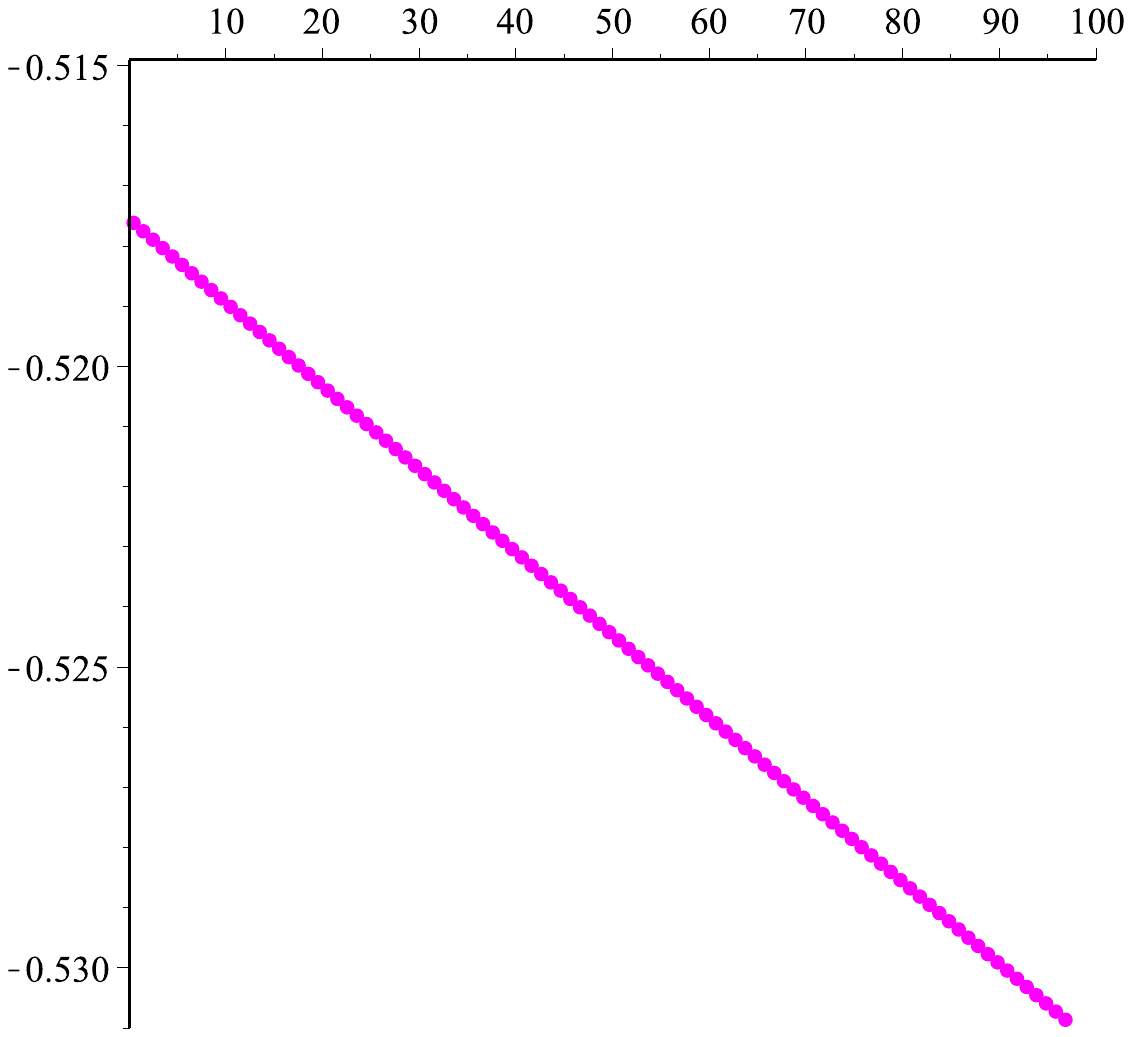}}}
\put(4,0){\rotatebox{0}{\includegraphics[width=7cm,height=7cm]{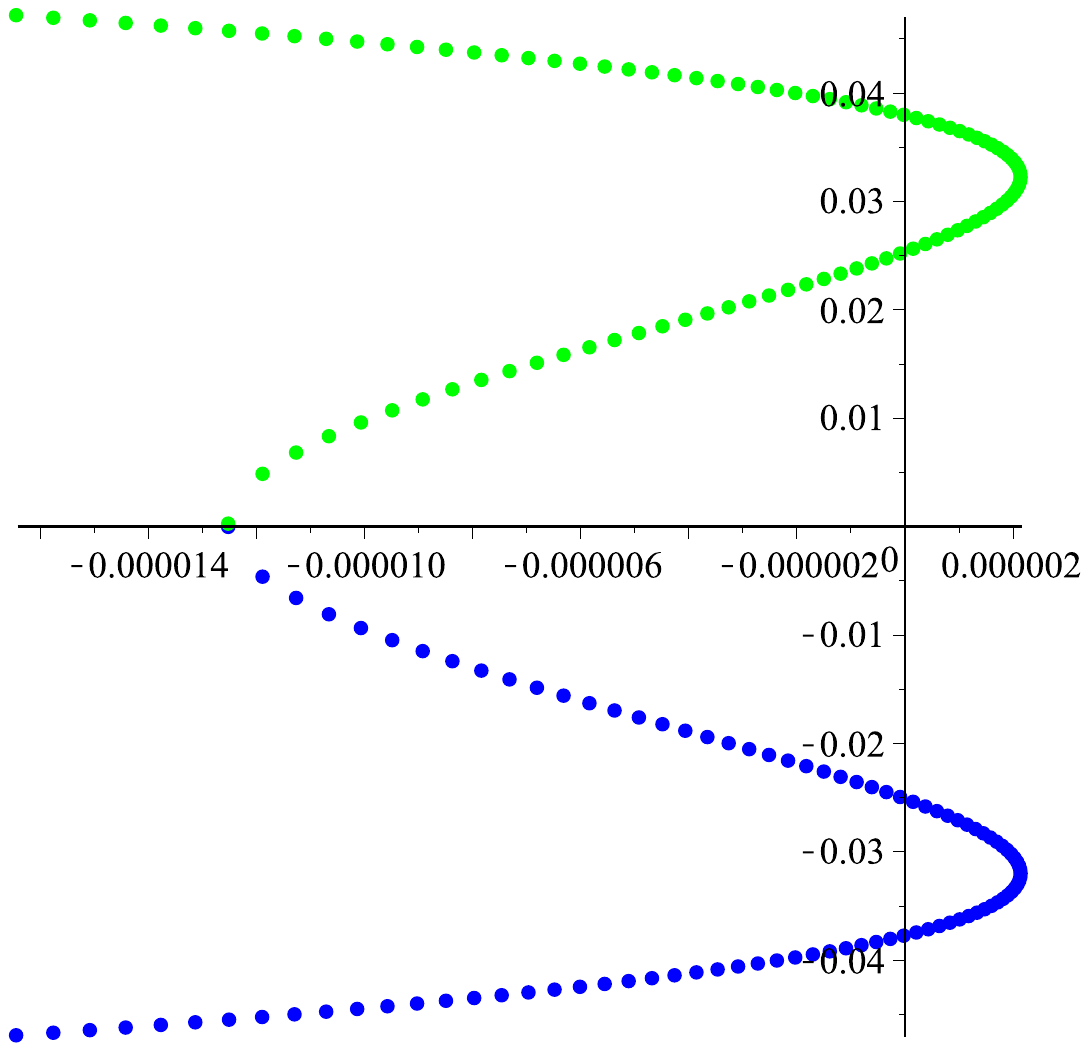}}}
\put(-0.3,6){\sc (a)}	
\put(-2.9,5.4){$\sc \eta^\star$}
\put(2.1,5){$\sc \eta_f$}	
\put(2.6,5.2){$\sc S_{in}$}	
\put(0,2.7){\sc {\color{magenta} $\lambda_1(S_{in})$}}
\put(7,6){{\sc (b)}}
\put(9.6,5.5){$\sc \beta(S_{in})$}	
\put(10.2,3.3){\sc $\alpha(S_{in})$}
\put(5.6,4){\sc {\color{green} $\lambda_2(S_{in})$}}
\put(6.2,3.9){\sc {\color{green} \vector(2,1){0.5}}}
\put(5.6,2.35){\sc {\color{blue} $\lambda_3(S_{in})$}}
\put(6.2,2.6){\sc {\color{blue} \vector(2,-1){0.5}}}
\end{picture}
\end{center}
\vspace{-1.6cm}
\caption{Case of the parameter set in Table \ref{TabParamVal} (line 5): variation of $S_{in}$ from $\eta^\star$ to $\eta_f$ when $D=D^\ast$; (a) the real eigenvalue $\lambda_1(S_{in})$. (b) The pair of complex-conjugate eigenvalues $\lambda_{2,3}(S_{in})$.} \label{FigEigenValSinCasArima}
\end{figure}

\begin{figure}[!h]
\setlength{\unitlength}{1.0cm}
\begin{center}
\begin{picture}(8.5,7.3)(0,0)
\put(-3.7,0){\rotatebox{0}{\includegraphics[width=10cm,height=17cm]{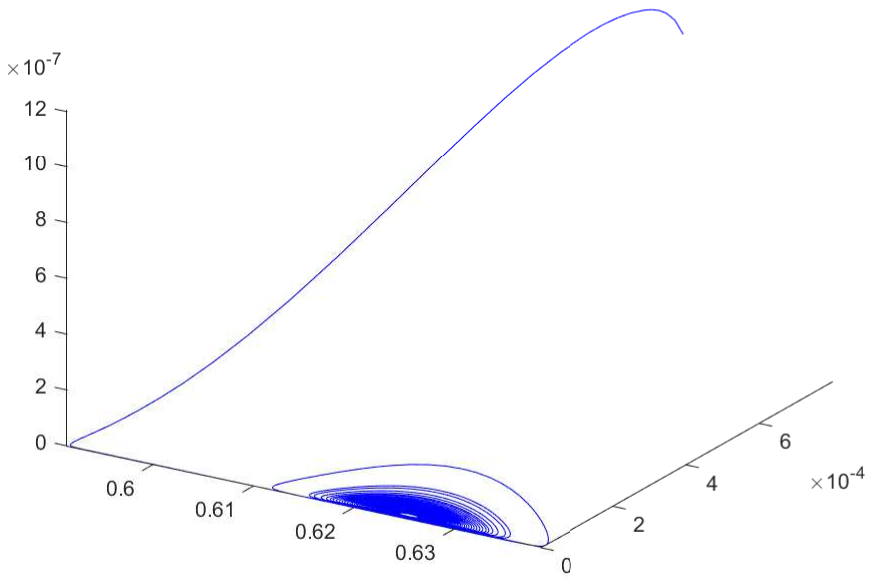}}}
\put(3.7,0){\rotatebox{0}{\includegraphics[width=10cm,height=16cm]{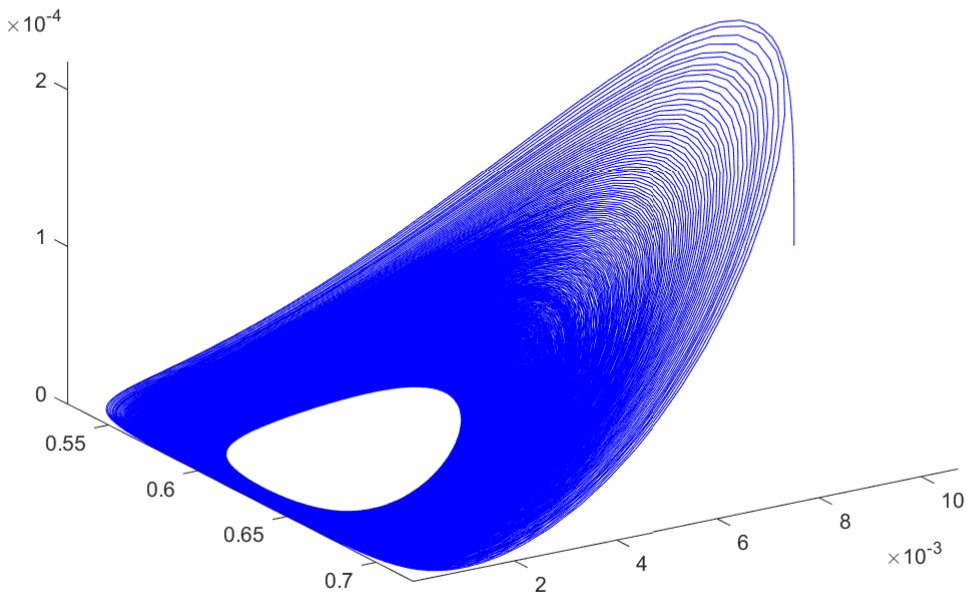}}}
\put(-0.3,6){\sc (a)}	
\put(-3.5,3.6){$\sc v$}
\put(-1.5,0.7){$\sc S$}	
\put(2.7,1){$\sc u$}
\put(-0.1,1.1){{\sl {\color{red}$\bullet$  }}}	
\put(7.7,6){\sc (b)}	
\put(4.3,3.9){$\sc v$}
\put(6,1.3){$\sc S$}	
\put(9.5,1){$\sc u$}	
\end{picture}
\end{center}
\vspace{-1.3cm}
\caption{Case of the parameter set in Table \ref{TabParamVal} (line 5): the three-dimensional space $(S,u,v)$ in MATCONT; (a) convergence to the positive steady state $E_1$ when $S_{in}=1$; (b) convergence to the stable limit cycle when $S_{in}=48$.} \label{Fig3DCasArima}
\end{figure}

To give numerical evidence of the Hopf bifurcation occurring through the transition through the curve $\Gamma_{\sH}$ from region $\mathcal{I}_1$ to region $\mathcal{I}_3$,
we determine numerically the eigenvalues of the Jacobian matrix of system (\ref{ModFlocGen}) at $E_1^1$ by computing the roots of the characteristic polynomial as we vary
the parameter $S_{in}$.
Let $D$ be fixed such that $D=D^\ast=2.5\,10^{-5}$.
Fig. \ref{FigEigenValSinCasArima}(a) shows that one eigenvalue denoted by $\lambda_1(S_{in})$ is real and remains negative for all $S_{in}\in [\eta^\star,\eta_f]$ where
 $\eta^\star=\lambda_u(D^\ast) \approx 0.625$ denotes the value of $S_{in}$ at which the positive steady state appears and $\eta_f=100$ denotes the final value of the variation of $S_{in}$.
Fig. \ref{FigEigenValSinCasArima}(b) shows that the two other eigenvalues $\lambda_2(S_{in})$ and $\lambda_3(S_{in})$ defined by
$$
 \lambda_{2,3}(S_{in})=\alpha_{2,3}(S_{in}) \pm i \beta_{2,3}(S_{in}),\quad \mbox{for all} \quad S_{in}\in \left[\eta^\star,\eta_f\right]
$$
are complex-conjugate so that the real part $\alpha_{2,3}(S_{in})$ is negative for all $S_{in}\in \left[\eta^\star,\eta_1\right)\cup(\eta_2,\eta_f]$ and positive for all $S_{in}\in (\eta_1,\eta_2)$.
When $S_{in} = \eta_i$, $i=1,2$, the pair $\lambda_{2,3}(\eta_i)$ is purely imaginary such that $\alpha_{2,3}(\eta_i)=0$, with $\beta_{2,3}(\eta_i)\neq0$.
Moreover, the following transversality condition is checked numerically
\begin{equation}                                   \label{IneqDiffRel}
\frac{d \alpha_{2,3}}{d S_{in}}(\eta_1) > 0 \quad \mbox{and} \quad  \frac{d \alpha_{2,3}}{d S_{in}}(\eta_2) < 0.
\end{equation}
that is, the two complex-conjugate eigenvalues cross the imaginary axis with non-zero speed.
Thus, the positive steady state $E_1$ is destabilized via two Hopf bifurcations with the occurrence or disappearance of a stable limit cycle
when $S_{in}$ increases and crosses the critical values $\eta_1$ and $\eta_2$. This result is consistent with the numerical simulation in Fig. \ref{Fig3DCasArima}(b) showing the emergence of a stable limit cycle where the oscillations are sustained.

Recall that $D$ is fixed at $D=D^\ast=2.5\,10^{-5}$. Fig. \ref{Fig3DCasArima}(a) illustrates the convergence towards $E_1$ in the three-dimensional phase space $(S,u,v)$ when $S_{in}=1$ where the pair of complex-conjugate eigenvalues have negative real parts. In this case, the point $(S_{in},D)$ belongs to region $\mathcal{I}_1$ where there are only two steady states: $E_0$ is unstable while $E_1$ is LES.

Fig. \ref{Fig3DCasArima}(b) illustrates the convergence towards a stable limit cycle when $S_{in}=48$ where the pair of complex-conjugate eigenvalues have positive real parts. In this case, the point $(S_{in},D)$ belongs to region $\mathcal{I}_3$ where there are only two steady states $E_0$ and $E_1$ which are unstable.
To solve the problem of the calculation time of the solution of (\ref{ModFlocGen}) until convergence to the limit cycle where $D$ is small enough, we have changed the default solver ``{\it{ode45}}'' to ``{\it{ode23}} '' in MATCONT.

\begin{figure}[!h]
\setlength{\unitlength}{1.0cm}
\begin{center}
\begin{picture}(6.1,6)(0,0)
\put(-5.2,0){\rotatebox{0}{\includegraphics[width=5.7cm,height=7cm]{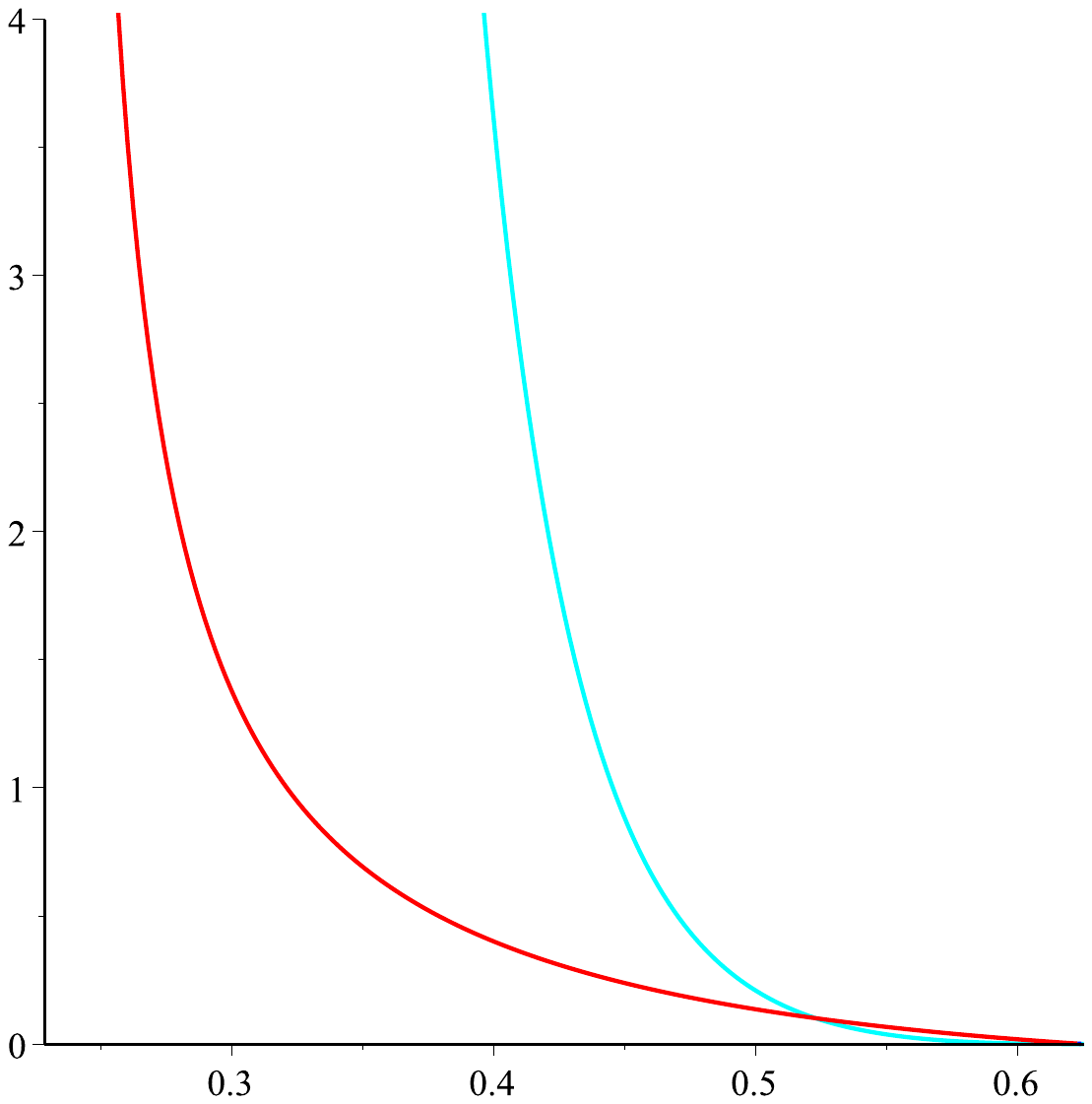}}}
\put(0,0){\rotatebox{0}{\includegraphics[width=5.7cm,height=7cm]{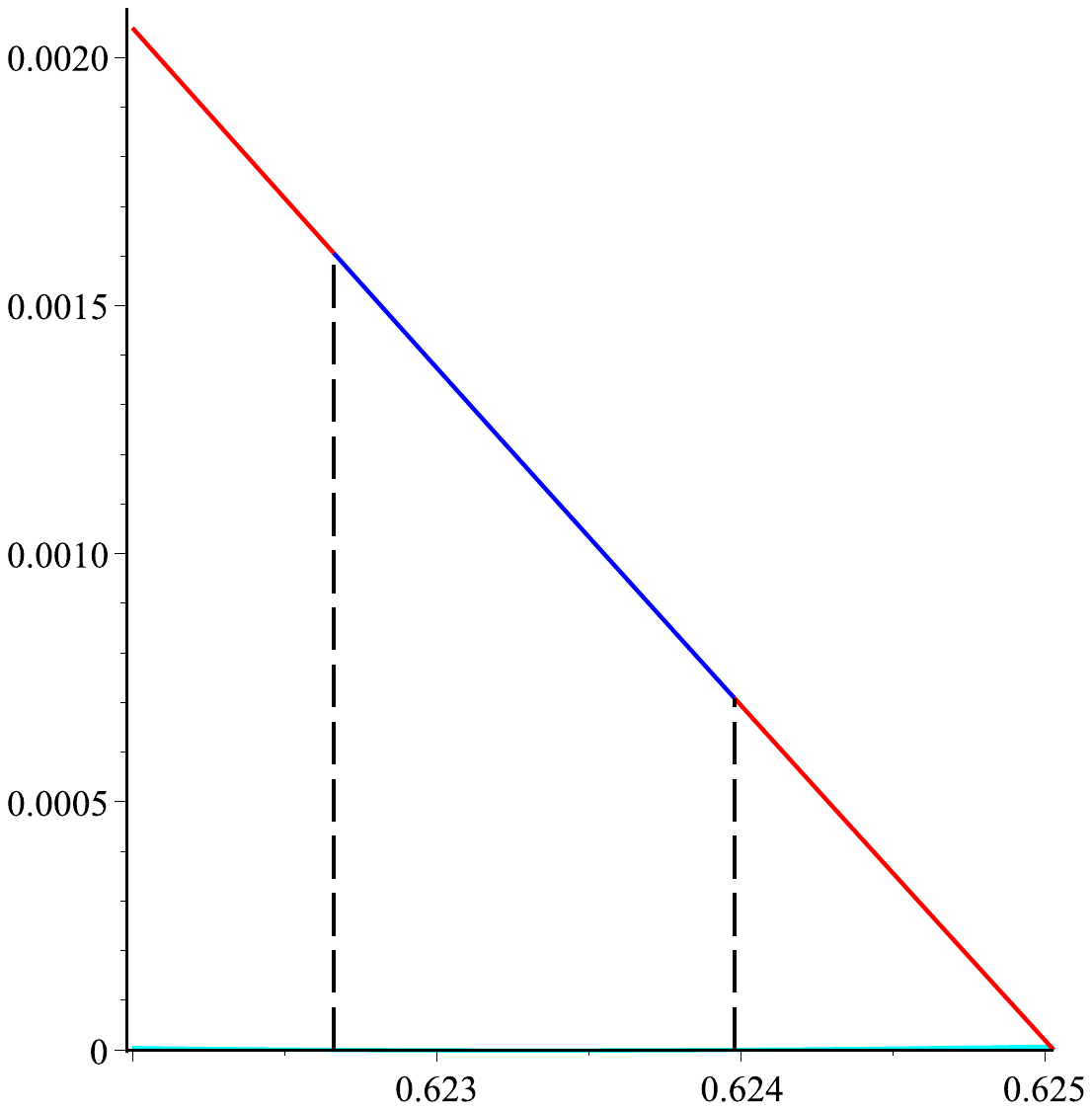}}}
\put(5.2,0){\rotatebox{0}{\includegraphics[width=5.7cm,height=7cm]{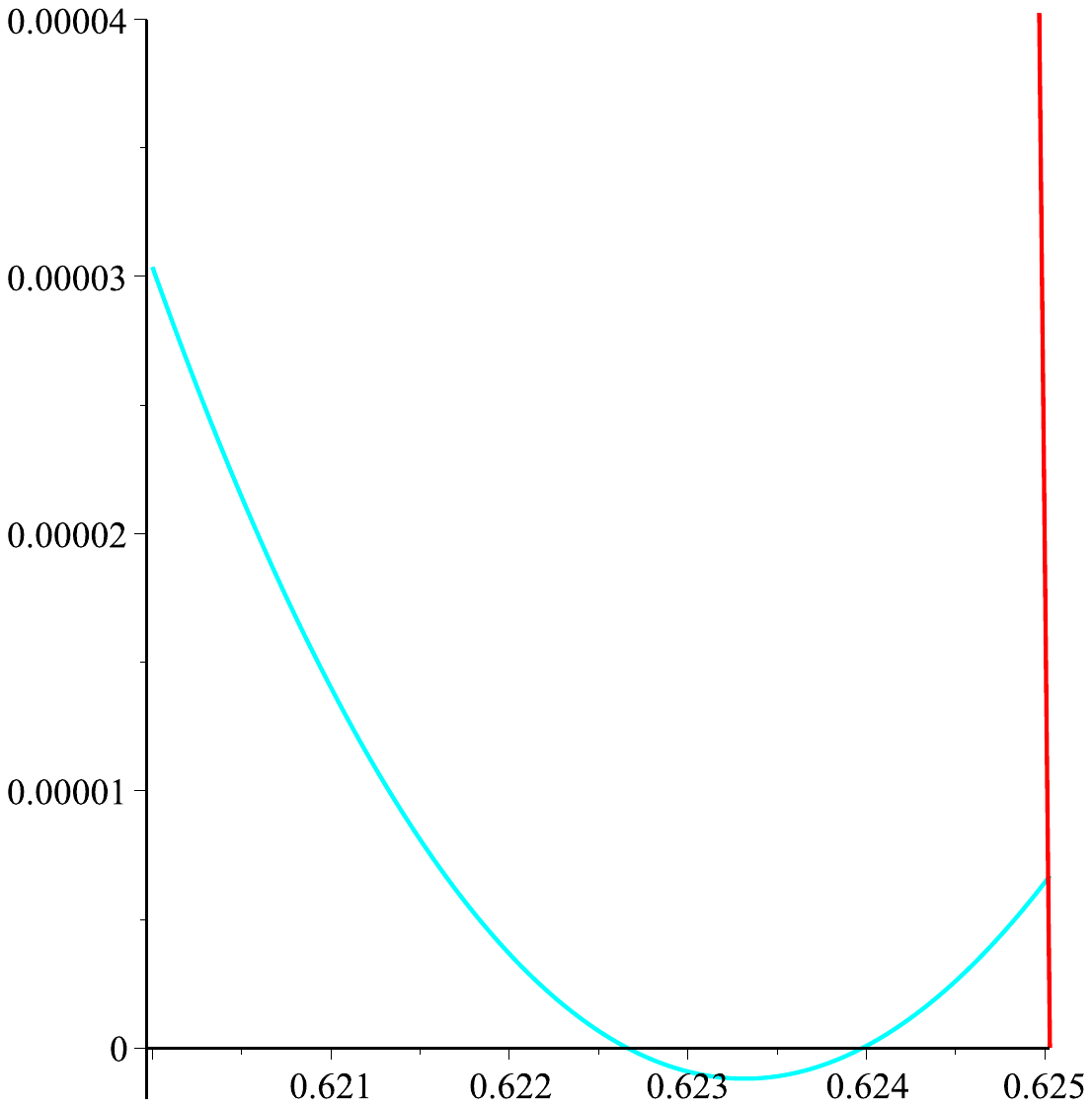}}}
\put(-2.1,5.8){{\sc $(a)$}}
\put(-2.55,5.4){{\sc ${\color{cyan}c_4}$}}
\put(-4.2,5.4){{\sc ${\color{red}H}$}}
\put(0,1.1){{\sc  $\lambda_u$}}
\put(0.15,1.35){{\sc  $S$}}
\put(2.9,5.8){{\sc $(b)$}}
\put(1.35,5.4){{\sc ${\color{red}H}$}}
\put(1.1,1.45){{\sc ${\color{cyan}c_4>0}$}}
\put(2.5,1.45){{\sc ${\color{cyan}c_4<0}$}}
\put(4.1,1.45){{\sc ${\color{cyan}c_4>0}$}}
\put(1.7,1.1){{\sc  $S_{\sH}^2$}}
\put(3.35,1.48){{\sc  $S_{\sH}^1$}}
\put(5,0.95){{\sc  $\lambda_u$}}
\put(5.2,1.35){{\sc  $S$}}
\put(8.1,5.8){{\sc $(c)$}}
\put(6.5,4.3){{\sc ${\color{cyan}c_4}$}}
\put(8.35,1.45){{\sc  $S_{\sH}^2$}}
\put(9.35,1.45){{\sc  $S_{\sH}^1$}}
\put(10.3,5.4){{\sc ${\color{red}H}$}}
\put(10.2,0.95){{\sc  $\lambda_u$}}
\put(10.4,1.35){{\sc  $S$}}
\end{picture}
\end{center}
\vspace{-1.6cm}
\caption{Case of the parameter set in Table \ref{TabParamVal} (line 5): curve of the function $H(S)$ where $c_4(S)$ is positive when $S\in\left[0.622,S_{\sH}^2\right[\cup\left]S_{\sH}^1,\lambda_u\right]$ and negative when $S\in\left]S_{\sH}^2,S_{\sH}^1\right[$. Magnifications when (b) $(S,H(S))\in[0.622,\lambda_u]\times\left[-5\,10^{-6},2.1\,10^{-3}\right]$ and (c) $(S,H(S))\in[0.620,\lambda_u]\times\left[-2\,10^{-6},4\,10^{-5}\right]$.}\label{Fig-HS_C4_CasArimaDfix}
\end{figure}

Fig. \ref{Fig-HS_C4_CasArimaDfix} illustrates the curve of the function $S\mapsto H(S)$ in red [resp. in blue] when the function $S\mapsto c_4(S)$ is positive [resp. negative] and $D$ is fixed at $D^\ast=2.5\,10^{-5} \in\left]0,D_{\sH}^{max}\right[$.
The solutions $S_{\sH}^1$ and $S_{\sH}^2$ of the equation $c_4(S)=0$ correspond to the critical values $S_{in}^{H1}$ and $S_{in}^{H2}$.
They are the intersections of the horizontal line of equation $D=D^\ast$ with the curves $\Gamma_{\sH}^1$ and $\Gamma_{\sH}^2$, respectively,
in the $(S_{in},D)$-plane of the operating diagram in Fig. \ref{FigDO-Arima}(b).
By increasing the value of $S_{in}$ from zero to $\eta^\star=\lambda_u\left(D^\ast\right)\approx 0.625$, $E_0$ becomes unstable by a Branch Point (BP) with $E_1$ that
appears stable until the first Hopf bifurcation at $S_{in}^{H1}\approx 28.990$ $\left(\mbox{or equivalently } S=S_{\sH}^1\approx 0.62398\right)$.
Then, $E_1$ remains unstable up to the value of $S_{in}^{H2}\approx 64.878$ $\left(\mbox{or equivalently } S=S_{\sH}^2\approx 0.62267\right)$, that is, for all $S\in \left]S_{\sH}^2,S_{\sH}^1\right[$. Finally, for $S_{in}>S_{in}^{H2}$ $\left(\mbox{or equivalently } S \in \left]\lambda_v(D^\ast),S_{\sH}^2\right[\right)$, $E_1$ returns stable via a second Hopf bifurcation.
\section{Parameter values used in numerical simulations} \label{Sec-AppedixParamVal}
All the parameter values used in the numerical simulations are provided in Table \ref{TabParamVal}.
\begin{table}[!h]
{\footnotesize
\caption{Parameter values used for system (\ref{ModFlocGen}) when the specific growth rates $f$ and $g$ are given by (\ref{SpeciFunc}). The abbreviation Var means Variable.} \label{TabParamVal}
\vspace{-0.1cm}
\begin{center}
\begin{tabular}{@{\hspace{2mm}}l@{\hspace{2mm}} @{\hspace{2mm}}l@{\hspace{2mm}} @{\hspace{2mm}}l@{\hspace{2mm}} @{\hspace{2mm}}l@{\hspace{2mm}}
                @{\hspace{2mm}}l@{\hspace{2mm}} @{\hspace{2mm}}l@{\hspace{2mm}} @{\hspace{2mm}}l@{\hspace{2mm}} @{\hspace{2mm}}l@{\hspace{2mm}}
                @{\hspace{2mm}}l@{\hspace{2mm}} @{\hspace{2mm}}l@{\hspace{2mm}} @{\hspace{2mm}}l@{\hspace{2mm}} @{\hspace{2mm}}l@{\hspace{2mm}}
                @{\hspace{2mm}}l@{\hspace{2mm}} }
 Parameter                 &   $m_1$      &  $k_1$   &   $m_2$     &  $k_2$   &   $a$      &  $b$        &   $\alpha$  &  $\beta$    &   $m_u$      &  $m_v$       &  $y_{u,v}$  &  $\overline{D}$\\
&   $\left(h^{-1}\right)$ &  $(g/l)$ &  $\left(h^{-1}\right)$ &  $(g/l)$ &  $(l/h/g)$ &  $\left(h^{-1}\right)$ &  &  &  $\left(h^{-1}\right)$  &  $\left(h^{-1}\right)$  &       &  $\left(h^{-1}\right)$   \\ \hline
Figs. \ref{Fig-CasFig13Siam}(a), \ref{Fig-CasFig13SiamHC4} &  4.5    &   1   &     3    &   2.7    &    2    &    3   &    0.8   &   0.5   &  0.2   &    0.25   & 1 & \\ \hline
\begin{tabular}{@{\hspace{0mm}}l@{\hspace{0mm}}}
Figs. \ref{Fig-CasFig13Siam}(b-c), \ref{FigDO-Fig12Siam-Num}, \ref{Fig-CasFig12Siam-Hc4}  \\
Figs. \ref{Fig-CasFig6Siam}, \ref{FigDO-Fig6Siam-Num}, \ref{FigDB-CasFig6Siam-MC}, \ref{Fig-Hc4-CasFig6Siam}-\ref{FigT-HomCasFig6Siam}  \\
Figs. \ref{FigDO-LP-H-a05b2}, \ref{FigDO-EffectFloc}
\end{tabular}
                           &   5    &   2   &     5    &   3    &
\begin{tabular}{@{\hspace{0mm}}l@{\hspace{0mm}}}
       4\\
       4 \\
       Var
\end{tabular}
                  &
 \begin{tabular}{@{\hspace{0mm}}l@{\hspace{0mm}}}
       2\\
       2 \\
       Var
\end{tabular}
              &    1   &
\begin{tabular}{@{\hspace{0mm}}l@{\hspace{0mm}}}
               0.9\\
               1 \\
               1
\end{tabular}                  &  3.25  &
\begin{tabular}{@{\hspace{0mm}}l@{\hspace{0mm}}}
               0\\
               1 \\
               1
\end{tabular}
& 1 &
\begin{tabular}{@{\hspace{0mm}}l@{\hspace{0mm}}}
               0.130\\
               0.032\\
               Var
\end{tabular}
\\ \hline
Figs. \ref{FigDO-Arima}-\ref{Fig-HS_C4_CasArimaDfix} &  3.5    &   2.5   &     3    &   1.5    &    1    &    1   &    1   &   0.75   &  0.7  &    0.4   & 1 & 0.460
\end{tabular}
\end{center}
}\vspace{-0.4cm}
\end{table}
\section*{Acknowledgments}
The first author thanks the financial support of Cimpa "research in pairs", the European Mathematical Society EMS-Simons for Africa, and also the I-SITE Excellence Program of the University of Montpellier for projects ``Support for international mobility - EXPLORE\#4''. We thank the Euro-Mediterranean research network
\href{http://www.inra.fr/treasure}{TREASURE}.
 
\end{document}